\documentclass[fleqn,10pt]{article}
\usepackage{latexsym, graphicx, epsfig, amsmath, amssymb,amsfonts}
\usepackage{natbib,amsthm,version}
\usepackage{amsbsy,bm,multirow,enumerate}
\usepackage[titletoc,page]{appendix}
\usepackage[mathscr]{eucal}
\usepackage{mathtools}
\usepackage{color}
\usepackage[utf8]{inputenc}
\usepackage[english]{babel}
\usepackage{amsthm}
\usepackage{enumerate}
\usepackage[hidelinks]{hyperref}
\usepackage{url}
\usepackage{subfigure}
\usepackage[lined,boxed]{algorithm2e}

\usepackage{hyperref}
\usepackage{amsmath,amssymb,amsthm}
\usepackage{graphicx}
\usepackage{graphics}
\usepackage {graphicx,fancyhdr}
\usepackage{color}
\usepackage{flafter}
\numberwithin{equation}{section}
\usepackage{multirow}
\usepackage{mathrsfs} 
\usepackage{subfigure}
\usepackage{tabularx}
\usepackage{lineno}
\usepackage{xfrac}

\usepackage{stmaryrd} 
\usepackage{bm} 
\usepackage{cases} 

\DeclareMathOperator{\esssup}{ess\,sup}
\newtheorem{myTheorem}{Theorem}[section]
\newtheorem{myLemma}[myTheorem]{Lemma}
\newtheorem{myRemark}[myTheorem]{Remark}

{\normalfont}

\newcommand{\bs}[1]{\boldsymbol{#1}}

\theoremstyle{definition}

\catcode`\@=11 \theoremstyle{plain}
\@addtoreset{equation}{section}   

\@addtoreset{figure}{section}
\renewcommand\thefigure{\thesection.\@arabic\c@figure}

\theoremstyle{remark}

\def \epsilon {{\varepsilon}}

\definecolor{bgblue}{rgb}{0.04,0.39,0.54}
\definecolor{lired}{rgb}{0.3, 0.0, 0.0}
\definecolor{ligreen}{rgb}{0.0, 0.3, 0.0}
\definecolor{liblue}{rgb}{0.9, 1.0, 1.0}
\definecolor{gray}{rgb}{0.6, 0.6, 0.6}
\definecolor{sky}{rgb}{0.3, 1.0, 1.0}
\definecolor{bunhong}{rgb}{1.0, 0.3, 1.0}
\definecolor{yellow}{rgb}{0.97, 1, 0.0}
\definecolor{liyellow}{rgb}{0.9, 0.8, 0.0}
\definecolor{cengse}{rgb}{0.00,0.40,0.29}

\renewcommand \wedge \times

\graphicspath{{./figs/}}

\title{
  gPAV-Based Unconditionally Energy-Stable Schemes 
  for the Cahn-Hilliard Equation: Stability and Error Analysis
} 
\author{
  Yanxia Qian$^{1,2}$, \ Zhiguo Yang$^{2}$, \ Fei Wang$^1$, \
  Suchuan Dong$^2$\thanks{Authors of correspondence.
    Email: sdong@purdue.edu } \\
  $^1$School of Mathematics and Statistics \\
  Xi'an Jiaotong University, China \\
  $^2$Center for Computational and Applied Mathematics \\
  Department of Mathematics \\
  Purdue University, USA \\
 } 

\date{(June 14, 2020)}
\begin{document}
\maketitle



\begin{abstract}

  We present several first-order and second-order
  numerical schemes for the Cahn-Hilliard equation with discrete
  unconditional energy stability. These schemes stem
  from the generalized Positive Auxiliary Variable (gPAV) idea,
  and require only the solution of linear algebraic systems with a
  constant coefficient matrix.
  More importantly, the computational complexity
  (operation count per time step) of these schemes
  is approximately a half of those of the gPAV and the scalar auxiliary
  variable (SAV) methods in previous works.
  We investigate the stability properties of the proposed schemes
  to establish stability bounds for the field function and the auxiliary
  variable, and also provide their error analyses.
  Numerical experiments are presented to verify the theoretical
  analyses and also demonstrate the stability of the schemes at
  large time step sizes.

\end{abstract}


\vspace{0.05cm}
Keywords: {\em 
  energy stability;
  auxiliary variable;
  generalized positive auxiliary variable;
  scalar auxiliary variable;
  Cahn-Hilliard equation
}

\section{Introduction}
\label{sec:intro}

%

The Cahn-Hilliard equation~\cite{CahnH1958} plays a key role
in the modeling of two-phase and multiphase flows based on the
phase field (or diffuse interface) approach~\cite{LowengrubT1998,LiuS2003,AbelsGG2012,Dong2017}.
Under appropriate
boundary conditions, the mass (or volume) of each fluid component described
by the Cahn-Hilliard equation is conserved.
Indeed, the Cahn-Hilliard equation can be derived from the mass
balance equations for individual fluid components in
a multicomponent fluid mixture by choosing an appropriate
form for the free energy density function;
see e.g.~\cite{AbelsGG2012,Dong2014,Dong2018}.
Developing effective and efficient numerical algorithms for
the Cahn-Hilliard equation can have important ramifications
on the modeling and simulation of two-phase and multiphase flows.
This problem has witnessed a sustained interest from
the research community, and we refer to~\cite{ShenY2015,TierraG2015,YanCWW2018,YangLD2019} for some examples.

The design of numerical schemes for the Cahn-Hilliard equation
confronts several challenges. The predominant issue among these
is posed by the nonlinear term. The nonlinear term in the Cahn-Hilliard equation
stems from
the potential energy (double-well) in the free energy density function.
The system described by the Cahn-Hilliard equation admits
an energy balance equation (energy law) on the continuous level.
To achieve discrete energy stability in the numerical scheme,
i.e.~retaining a corresponding
discrete energy law,  hinges on the numerical
treatment of the nonlinear term. 
Energy-stable schemes allow the use of larger time step sizes,
which can potentially accelerate dynamic simulations if
its computational cost per time step is manageable.
A downside 
about energy-stable schemes
is that their cost is typically markedly higher when
compared with semi-implicit type schemes (see e.g.~\cite{BadalassiCB2003,YueFLS2004,DongS2012,Dong2014obc,GonzalezT2013,Dong2017}),
which are only conditionally stable.
This is because the energy-stable schemes oftentimes
entail the solution of nonlinear algebraic equations,
or the solution of linear algebraic systems
(either coupled linear systems or a linear system for multiple times).
To achieve discrete energy stability and a low computational
complexity (or operation count) per time step in the numerical
scheme is the focus of the current work.


Numerical algorithms for the Cahn-Hilliard equation
available in the literature generally consist of
two classes: nonlinear schemes and linear
schemes. With nonlinear schemes one treats the nonlinear term
or a part of the nonlinear term implicitly, and this requires
the solution of nonlinear algebraic equations
upon discretization; see e.g.~\cite{ElliotFM1989,Furihata2001,FengP2004,KimKL2004,MelloF2005,Feng2006,Wise2010,HuaLLW2011},
among others.
Among the nonlinear schemes,
convex splitting of the potential energy~\cite{ElliotS1993,Eyre1998} and
its variants are a widely-used approach for treating
the nonlinear term.
Other approaches include the midpoint approximation~\cite{ElliotFM1989,LinLZ2007},
Taylor expansion approximation~\cite{KimKL2004}, and
special quadrature rules~\cite{GomezH2011}, among others.

Unlike nonlinear schemes,
the linear schemes (see e.g.~\cite{ShenY2010,GonzalezT2013,Yang2016,ShenXY2018})
require only the solution of linear algebraic systems
upon discretization, due to the explicit treatment of
the nonlinear term, while maintaining  energy stability.
Among the linear schemes, incorporation of a stabilizing zero term,
together with a modified potential energy with bounded second derivative,
is an often-used method~\cite{ShenY2010,GonzalezT2013}. 
Other researchers have also employed a Lagrange multiplier to
enforce  energy stability~\cite{BadiaGG2011,GonzalezT2013}.
In the past few years, the use of certain auxiliary functions or variables
proves to be effective in  devising linear energy-stable
schemes. 
The invariant energy quadratization (IEQ)~\cite{Yang2016} and
the scalar auxiliary variable (SAV)~\cite{ShenXY2018} are two
prominent examples of such methods;
see also~\cite{GongZYW2018,Shen2018Convergence,Zhaoetal2018,LiZW2019,YangLD2019,HouAX2019,YangD2019}, among others. 
%
The IEQ method introduces an auxiliary field function as an approximation
of the square root of the potential energy density function
together with a dynamic equation for this field function, and
allows one to ensure the energy stability relatively easily.
It gives rise to a system of linear algebraic equations
involving time-dependent coefficient matrices upon discretization.
The SAV method uses an auxiliary variable (a scalar number rather than
a field function) to approximate the square root of the potential energy
integral. It retains the ease to ensure the energy stability, and
moreover leads to linear algebraic systems with a constant coefficient
matrix, thus making the implementation considerably simpler~\cite{YangLD2019}.
The use of the square root function form in IEQ and SAV,
either for a field function and a scalar number,  is critical
to the proof of energy stability in these methods.
%
A recent further development in this area is the generalized Positive Auxiliary
Variable (gPAV) method~\cite{YangD2020}.
The gPAV method also employs a scalar-valued number
as the auxiliary variable to ensure the energy stability, and it gives rise
to a linear algebraic system with a constant coefficient
matrix. This method makes three advances  in terms of
the methodology: (i) gPAV allows the use of a general class of function forms
to define the auxiliary variable, not limited to the square root
function as in IEQ and SAV.
(ii) gPAV guarantees the positivity of the computed values for
the auxiliary variable.
(iii) gPAV applies to
general types of dissipative or conservative partial differential equations (PDE)
for the development
of energy-stable schemes, not limited to gradient type systems.


In the current paper we present several
unconditionally energy-stable linear
schemes with
first- and second-order accuracy
for solving the Cahn-Hilliard equation, and
provide analyses for their stability properties and errors.
These schemes stem from the gPAV idea~\cite{YangD2020}, and inherit
the useful properties of guaranteed positivity for
the computed auxiliary variable and constant coefficient matrix
for the resultant linear algebraic system upon discretization.
Two advances have been made algorithm-wise:
(i) Stability bounds for both the phase field function and
the auxiliary variable can be established with the current schemes.
In contrast, with the original gPAV scheme~\cite{YangD2020}
the stability property is only through the auxiliary variable.
(ii) The operation counts (or computational cost) per time step
of the current schemes are comparable to that of the semi-implicit
schemes (see e.g.~\cite{DongS2012}), and are about a half of
those of the gPAV scheme~\cite{YangD2020} and the
SAV scheme~\cite{ShenXY2018}.
This is because the linear system resulting from
the Cahn-Hilliard equation only needs to be solved once within
each time step with the current schemes. In contrast,
with gPAV~\cite{YangD2020} and SAV~\cite{ShenXY2018,YangLD2019}
the linear system needs to be solved twice
for the two copies of the field function therein
within each time step.
We provide the stability analyses and error estimates
for these schemes, and present numerical experiments
to verify the theoretical analyses.


The contributions of this work consist of two aspects:
(i) the unconditionally energy-stable schemes for
the Cahn-Hilliard equation, and
(ii) the stability and error analyses for the proposed schemes.


The rest of this paper is organized as follows.
In Section \ref{sec:method} we reformulate the Cahn-Hilliard equation
using the gPAV idea, and present two first-order and two
second-order schemes for numerically solving the reformulated system
of equations. We prove the unconditional energy stability properties
of these schemes and provide the error estimates.
In Section \ref{sec:tests} we present numerical examples to
verify the convergence rates and the unconditional stability 
of these schemes.
Section \ref{sec:summary} then concludes the presentations with
some closing remarks.
In the Appendix we provide proofs to several
theorems from the main text.


\section{gPAV-Based Unconditionally Energy-Stable Schemes}
\label{sec:method}

\subsection{Cahn-Hilliard Equation and gPAV Formulation}\label{Introduction}

Let $\Omega \in \mathbb{R}^d (d=2, 3)$ be a bounded domain with
a smooth boundary $\partial\Omega$.
We consider the following gradient flows
\begin{equation}\label{CH} 
\frac{\partial\phi}{\partial t} = \Delta\mu=\Delta\left(-\Delta \phi + \lambda \phi + h(\phi)\right),
\end{equation}
where $\phi(\bm{x},t)$ is the phase field function,
$\lambda\geqslant 0$ is a constant parameter,
$\Delta$ is the Laplace operator, and
$\bm{x}$ and $t$ and the spatial coordinate and time.
The nonlinear term
\begin{equation}\label{def_h}
h(\phi) = H'(\phi) = \phi^3-\phi, \quad\text{with} \
H(\phi)=\frac{1}{4}(\phi^2-1)^2
\end{equation}
being the double-well
potential function~\cite{CahnH1958,AllenC1979}.  As is well-known,
this is  the celebrated Cahn-Hilliard equation with $\lambda=0$. 

The equation \eqref{CH} is supplemented by the initial
condition \begin{equation}\label{initial_condition}
\phi(\bm{x},t=0) = \phi_{in}(\bm{x})
\end{equation}
where $\phi_{in}$ denote the initial phase field distribution,
and the boundary conditions of either
\begin{equation}\label{boundary}
\nabla\phi\cdot\bm{n}=\nabla\mu\cdot\bm{n}=0 \quad \text{on}\ \partial\Omega, 
\end{equation}
or the periodic boundary conditions for $\phi$.
Here $\bm n$ denotes the outward-pointing unit normal vector
of the boundary.

Taking the $L^2$ inner product between \eqref{CH} and $\phi$, using
the integration by parts and \eqref{boundary},
we derive the following free energy functional $E_{tot}$
for this system
\[
E_{tot}(t)=\int_{\Omega}\left(\frac{1}{2}|\nabla \phi|^2+\frac{\lambda}{2}\phi^2+H(\phi)\right)d\bm{x}.
\]
To facilitate energy-stable numerical approximations of the system \eqref{CH}, we define a shifted energy of the following form
\begin{equation}\label{energy_0}
E(t)=\bm{E}[\phi]=\int_{\Omega}\left(\frac{1}{2}|\nabla \phi|^2+\frac{\lambda}{2}\phi^2+H(\phi)\right)d\bm{x}+c_0,
\end{equation}
where $c_0$ is a chosen energy constant to ensure
that $E(t)>0$ for $0\leq t\leq T$, and $T$ denotes
the time interval on which the system to be solved.
It is straightforward to verify that
the system \eqref{CH}--\eqref{boundary} satisfies the energy dissipation law
\begin{equation}\label{energy}
\frac{dE}{dt}=\int_{\Omega}\left(\nabla \phi\cdot\nabla\phi_t+\lambda\phi\phi_t+h(\phi)\phi_t\right)d\bm{x}=\int_{\Omega}\phi_t\mu d\bm{x}=-\int_{\Omega}|\nabla\mu|^2d\bm{x}\leq 0,
\end{equation}
where $\phi_t$
denotes the time derivative of $\phi$.

Following the gPAV idea~\cite{YangD2020}, we introduce a
positive scalar variable $R(t)^2=E(t)$ (or $R(t)=\sqrt{E(t)}$).
$R(t)$ satisfies the following evolution equation
\begin{equation}\label{energy_1}
2R\frac{dR}{dt}=\frac{dE}{dt}=-\int_{\Omega}|\nabla\mu|^2d\bm{x}.
\end{equation}
 Noting that $\frac{R}{\sqrt{E}}=1$, we rewrite \eqref{energy_1} into
\begin{equation}\label{energy_2}
  \frac{dR}{dt}=-\frac{1}{2R}\int_{\Omega}|\nabla\mu|^2d\bm{x}
  =-\frac{1}{2\sqrt{E}}\int_{\Omega}|\nabla\mu|^2d\bm{x}
  =-\frac{1}{2\sqrt{E}}\frac{R}{\sqrt{E}}\int_{\Omega}|\nabla\mu|^2d\bm{x}
  =-\frac{R}{2E}\int_{\Omega}|\nabla\mu|^2d\bm{x}.
\end{equation}
Then, we rewrite \eqref{CH} into the following equivalent form
\begin{subequations}\label{CH_revise}
	\begin{align}
	\label{CH_eq1}
	&\phi_t = \Delta\mu,\\
	\label{CH_eq2}
	&\mu =-\Delta\phi+ \lambda\phi + \frac{R^2}{E}h(\phi),\\
	\label{CH_eq3}
	&\frac{dR}{dt} = -\frac{R}{2E}\int_{\Omega}|\nabla\mu|^2d\bm{x}.
	\end{align}
\end{subequations}
In this reformulated system, the dynamic variables are $\phi$, $\mu$ and $R$,
which are coupled in the equations \eqref{CH_revise},
together with the boundary conditions \eqref{boundary}, the initial condition \eqref{initial_condition} for $\phi$, and the following initial condition for $R(t)$ 
\begin{equation}\label{energy_3}
R(0)=\sqrt{\bm{E}[\phi_{in}]}, \qquad {\rm where}\ \bm{E}[\phi_{in}]=\int_{\Omega}\left(\frac{\lambda}{2}\phi_{in}^2+\frac{1}{2}|\nabla\phi_{in}|^2+H(\phi_{in})\right)dx+c_0. 
\end{equation}
Note that $R(t)$ in this system is determined by solving the coupled system of equations, not by using the relation $R^2(t)=\bm{E}(\phi)$.
Therefore $R^{2}(t)$ is an approximation of $E(t)$,
rather than $E(t)$ itself.

\subsection{Preliminaries}\label{Preliminaries}

We first outline the notation used herein
and recall some basic results, including the existence, uniqueness, and regularity results about the $H^{-1}$ gradient flows.

For the non-negative integers $p$, $k$ and an open Lipschitz subdomain $D \subset \overline{\Omega}$, let $L^p(D)$ denote the standard Banach space  with norm
$\|v\|_{0,p,D}=\left(\int_{D}|v|^pd\bm{x}\right)^{1/p}$
and $W^{k,p}(D)$ the standard Sobolev space
with the norm
$\|v\|_{k,p,D}=\left(\sum_{|\alpha|\leq k}\int_{D}\|D^{\alpha}v\|_{L^p}^pd\Omega\right)^{1/p}$. 
For simplicity, we take the Sobolev space $H^k(D)=W^{k,2}(D)$ with the norm $\|\cdot\|_{k,D}$ and semi-norm $|\cdot|_{k,D}$, and the space $H^0(D)=L^2(D)$ with the usual $L^2$-inner product $(\cdot,\cdot)_D$ and $L^2$-norm $\|\cdot\|_{0,D}$.  If $D$ is chosen as $\Omega$, we abbreviate them
by the norms $\|\cdot\|_{k}$, $\|\cdot\|_{0}$, the semi-norm $|\cdot|_{k}$ and the inner product $(\cdot,\cdot)$, respectively. Therefore, we introduce the space $L^p(0,T;V)$, $L^{\infty}(0,T;V)$ and $C(0,T;V)$ with the norms
\[
\|\varphi\|_{L^p(0,T;V)}=\left[\int_0^T\|\varphi(t)\|_{V}^pdt\right]^{1/p}, 
\quad \|\varphi\|_{L^{\infty}(0,T;V)}= \mathop{\esssup}\limits_{0\leq t\leq T}\|\varphi(t)\|_{V}\quad {\rm and} \quad
\|\varphi\|_{C(0,T;V)}=\max_{0\leq t\leq T}\|\varphi(t)\|_{V}.
\]

Assume that the nonlinear free energy potential $H(s) \in C^3(\mathbb{R})$. For some cases, in order to ensure the uniqueness, we assume the following condition: there exists a constant $c_1>0$ such that
\begin{equation}\label{bound}
H''(s)=h'(s)\geq -c_1.
\end{equation}

\begin{myLemma}[See \cite{Temam1997Infinite}]\label{sec2_Lemma1} (i) For the $H^{-1}$ gradient flow, assume that \eqref{bound} holds, $\phi_{in} \in L^2(\Omega)$ and there exists $p_0>0$ such that 
	\begin{equation}\label{Assume_H1}
	sh(s)\geq b |s|^{p_0}-c_2,
	\end{equation}
	where $b>0$ and $c_2$ are constants. Then, for any $T>0$, there exists a unique solution $\phi$ for \eqref{CH} such that
	\begin{equation}\label{Assume_u1}
	\phi\in C(0,T;L^2(\Omega))\cap L^2(0,T;H^2(\Omega))\cap L^{p_0}(0,T;L^{p_0}(\Omega)).
	\end{equation}
	(ii) For the $H^{-1}$ gradient flow, assume that $\phi_{in} \in H^2(\Omega)$ and 
	\begin{align}\label{Assume_H2}
	|h'(x)|<C(|x|^{p_1}+1), \qquad p_1>0 \ {\rm arbitrary} \ {\rm if}\ d=1,2; \quad 0<p_1<4 \ {\rm if}\  d=3,\\
	\label{Assume_H3}
	|h''(x)|<C(|x|^{p_2}+1), \qquad p_2>0 \ {\rm arbitrary} \ {\rm if}\ d=1,2; \quad 0<p_2<3 \ {\rm if}\  d=3,
	\end{align}
        where $d$ denotes the dimension in space.
	Then, for any $T>0$, there exists a unique solution $\phi$ for \eqref{CH} such that
	\begin{equation}\label{Assume_u2}
	\phi\in C(0,T;H^2(\Omega)) \cap L^2(0,T;H^4(\Omega)).
	\end{equation}
\end{myLemma}

\begin{myLemma}[See \cite{Shen2018Convergence}]\label{sec2_Lemma2} Assume that $\|\phi\|_1\leq M$. \\
	(i) Under the conditions of \eqref{Assume_H2} and $\phi \in H^3(\Omega)$, there exists $0<\sigma<1$ and a constant $C(M)$ such that 
	\begin{equation}\label{Assume_H4}
	\|\nabla h(\phi)\|_0^2\leq C(M)(1+\|\nabla\Delta\phi\|_0^{2\sigma}).
	\end{equation}
	(ii) Under the assumptions of \eqref{Assume_H2}, \eqref{Assume_H3} and $\phi \in H^4(\Omega)$, there exists $0<\sigma<1$ and a constant $C(M)$ such that
	\begin{align}\label{Assume_H5}
	\|\Delta h(\phi)\|_0^2\leq C(M)(1+\|\Delta^2\phi\|_0^{2\sigma}).
	\end{align}
\end{myLemma}

\begin{myLemma} [Discrete Gronwall lemma \cite{Heywood1990Finite}]\label{sec2_Lemma3} 
	Let $a_i$, $b_i$, $c_i$, $d_i$, $\Delta t$ and $C_0$, for integers $i\geq 0$, be non-negative numbers such that
	\[ a_n+\Delta t\sum_{i=0}^{n}b_i\leq \Delta t\sum_{i=0}^{n}d_ia_i+\Delta t\sum_{i=0}^{n}c_i+C_0,\quad\forall\,n\geq 0.\]
	Then, if $d_i\,\Delta t<1$ for all $i$,
	\[ a_n+\Delta t\sum_{i=0}^nb_i\leq \left(C_0 + \Delta t\sum_{i=0}^{n}c_i\right)\exp\left(\Delta t\sum_{i=0}^{n}\frac{d_i}{1-d_i\,\Delta t}\right),\quad\forall\, n\geq 0.
	\]
\end{myLemma}

\begin{myLemma} [Discrete Gronwall lemma \cite{Shen1990Long,He2003Two}]\label{sec2_Lemma4} 
	Let $a_i$, $b_i$, $c_i$, $d_i$, $\Delta t$ and $C_0$, for integers $i\geq 0$, be non-negative numbers such that
	\[ a_n+\Delta t\sum_{i=0}^{n}b_i\leq \Delta t\sum_{i=0}^{n-1}d_ia_i+\Delta t\sum_{i=0}^{n-1}c_i+C_0,\quad\forall\,n\geq 0.\]
	Then, 
	\[ a_n+\Delta t\sum_{i=0}^{n}b_i\leq \left(C_0 + \Delta t\sum_{i=0}^{n-1}c_i\right)\exp\left(\Delta t\sum_{i=0}^{n-1}d_i\right),\quad\forall\, n\geq 0.
	\]
\end{myLemma}

\subsection{First-order schemes}\label{gPAV_scheme}

We introduce several unconditionally energy-stable
schemes for solving the reformulated
Cahn-Hilliard equations~\eqref{CH_revise}, 
the boundary conditions~\eqref{boundary}, and
the initial conditions~\eqref{initial_condition}
and \eqref{energy_3}.
These schemes stem from the gPAV idea~\cite{YangD2020}, and
they inherit some of the attractive properties of gPAV.
For example, the auxiliary variable is computed
via a well-defined explicit form, and its
computed values are guaranteed to be positive.
The departure point lies in that all the schemes
presented herein require only the computation of
a single copy of the field functions per time step. In contrast,
the original gPAV method~\cite{YangD2020} entails
the computation of two copies of the field functions
within each time step. The amount of operations
involved in the current schemes is approximately a half of
that in the scheme of~\cite{YangD2020}.
The current schemes have a computational cost roughly the same as
the semi-implicit schemes for the Cahn-Hilliard
equation (see e.g.~\cite{DongS2012}).

We provide stability analysis and error estimates for these schemes
in what follows.
The first-order schemes are discussed in this section,
followed by the second-order schemes in the subsequent section.

\subsubsection{Scheme 1A}\label{gPAV}

Let $\Delta t>0$ be the time step size and
$n\geq 0$ denote the time step index,
and we set $t^n=n\Delta t$ for $0\leq n \leq N$ with $N= T/\Delta t$. For a generic
function $\chi(\bm{x},t)$ continuous in $t$, let $\chi^n$ denote
the approximation of $\chi(\bm{x},t^n)$ at time $t^n$,
where $\chi$ can be $\phi$, $\mu$ or $\psi$ (defined below).
Similarly, let $R^n$ denote the approximation of $R(t^n)$. Set
\begin{equation}\label{CH1_initial}
  \left\{
  \begin{split}
    &
    \phi^0=\phi_{in}, \ \ \mu^0=-\Delta \phi^0 + \lambda \phi^0 + h(\phi^0), \\
    &
    R^0=\sqrt{E(0)}=\sqrt{\int_{\Omega}\left(\frac{1}{2}|\nabla \phi^0|^2+\frac{\lambda}{2}|\phi^0|^2+H(\phi^0)\right)d\bm{x}+c_0}.
  \end{split}
  \right.
\end{equation}

Then given $\phi^n$, $\mu^n$ and $R^{n}$, we compute $\phi^{n+1}$, $\mu^{n+1}$ and $R^{n+1}$ through the following scheme,
\begin{subequations}\label{CH1_weak}
	\begin{align}
	\label{CH1_eq1}
	&\frac{\phi^{n+1}-\phi^{n}}{\Delta t} = \Delta\mu^{n+1},\\
	\label{CH1_eq2}
	&\mu^{n+1} =-\Delta\phi^{n+1}+ \lambda\phi^{n+1} + |\xi_{_{1A}}^{n+1}|^2h(\phi^{n}),\\
	\label{CH1_eq3}
	&\frac{R^{n+1}-R^{n}}{\Delta t} = -\frac{\xi_{_{1A}}^{n+1}}{2\sqrt{\bm{E}[\phi^{n}]}}\int_{\Omega}|\nabla\mu^{n}|^2d\bm{x},
	\end{align}
\end{subequations}
with the boundary conditions 
\begin{equation}\label{boundary_1}
  \nabla\phi^{n+1}\cdot\bm{n}=\nabla\mu^{n+1}\cdot\bm{n}=0,
  \quad \text{on}\ \partial\Omega,
\end{equation}
where 
\begin{equation}\label{xi_0}
\xi_{_{1A}}^{n+1} = \frac{R^{n+1}}{\sqrt{\bm{E}[\phi^{n}]}}.
\end{equation}
Note that here $\xi_{_{1A}}^{n+1}$ is an approximation
of the constant $\frac{R(t)}{\sqrt{E(t)}} = 1$.



\begin{myRemark}\label{rem_1A}
  In this scheme we have treated the nonlinear term $h(\phi)$
  in \eqref{CH1_eq2} and the $\left|\nabla\mu\right|^2$ term
  in \eqref{CH1_eq3} explicitly.
  Consequently, equation \eqref{CH1_eq3} for $R^{n+1}$ is
  not coupled with the equations~\eqref{CH1_eq1} and \eqref{CH1_eq2}
  for $\phi^{n+1}$ and $\mu^{n+1}$ on the discrete
  level.
  
\end{myRemark}

Substituting the $\xi_{_{1A}}^{n+1}$ expression in \eqref{xi_0} 
into equation \eqref{CH1_eq3}, we get
	\begin{equation}\label{xi}
	\xi_{_{1A}}^{n+1}=\frac{R^n}{\sqrt{\bm{E}[\phi^{n}]}+\frac{\Delta t}{2\sqrt{\bm{E}[\phi^{n}]}}\int_{\Omega}|\nabla \mu^n|^2d\bm{x}}.
	\end{equation}
Since $R^0>0$ according to equation \eqref{CH1_initial},
we conclude by induction that $\xi_{_{1A}}^{n}>0$ for all $n$.
Then $R^{n+1}$ is given by, in light of \eqref{xi_0},
\begin{equation}\label{R_1A}
  R^{n+1} = \xi_{_{1A}}^{n+1}\sqrt{E[\phi^n]}.
\end{equation}
We conclude that $R^{n+1}>0$ for all time steps
$n$.

The  time stepping with the current scheme is thus
as follows. Within a time step,
given $\phi^{n}$, $\mu^{n}$ and $R^{n}$,
we compute $\bm{E}[\phi^{n}]$ by \eqref{energy_0},
$\xi_{_{1A}}^{n+1}$ by \eqref{xi}, and $R^{n+1}$ by \eqref{R_1A}.
Then with $\xi_{_{1A}}^{n+1}$ known, we compute $\phi^{n+1}$ and $\mu^{n+1}$
by solving equations 
\eqref{CH1_eq1}-\eqref{CH1_eq2} together with
the boundary conditions~\eqref{boundary_1}.

It should be emphasized that the Cahn-Hilliard field equation
is only  solved once
per time step with the current scheme.
This is very different from the previous
gPAV and SAV-type
schemes (see e.g.~\cite{YangD2020,YangLD2019,ShenXY2018}),
which require solving the field equations twice per time step
(for two copies of the field variables therein).
Therefore, the operation count induced by the current scheme
is approximately a half of
those with the previous
SAV and gPAV schemes, and it
is comparable to that of
the semi-implicit type schemes (see e.g.~\cite{DongS2012}).
It can further be noted
that the auxiliary variables $R^{n+1}$ and $\xi_{_{1A}}^{n+1}$ are computed
by well-defined explicit forms, with their
values guaranteed to be positive.

\paragraph{Stability Properties}


The scheme given by equations \eqref{CH1_weak}--\eqref{xi_0}
is unconditionally stable. We summarize its stability properties
into several lemmas or theorems below.

\begin{myLemma}\label{sec3_Lemma1} The scheme \eqref{CH1_weak} is mass conserving
  in the sense that $(\phi^{n+1},1)=(\phi^n,1)$.
\end{myLemma}
\begin{proof}
  In light of the boundary conditions \eqref{boundary_1},
  the $L^2$ inner product between \eqref{CH1_eq1} and the constant one leads to
	\[
	(\phi^{n+1}-\phi^{n},1) = \Delta t(\Delta\mu^{n+1},1)=-\Delta t(\nabla\mu^{n+1},\nabla 1)=0.
	\]
	So the solution of \eqref{CH1_weak} satisfies $(\phi^{n},1)=(\phi^{0},1)$
        for any $n$.
\end{proof}

\begin{myLemma}\label{sec3_Lemma2} With the scheme \eqref{CH1_weak}
  for all time step $n$,
  \begin{equation}\label{sec3_eq1}
    0< R^{n+1} \leq R^n.
	\end{equation}
\end{myLemma}
\begin{proof} Multiplying $2\Delta tR^{n+1}$ to equation \eqref{CH1_eq3}
  and using equation \eqref{xi_0},
  we get
	\begin{equation}\label{sec3_eq2}
	  |R^{n+1}|^2 - |R^{n}|^2\leq|R^{n+1}|^2 - |R^{n}|^2 + |R^{n+1}-R^{n}|^2=- \frac{\Delta t |R^{n+1}|^2}{\bm{E}[\phi^{n}]}\int_{\Omega}|\nabla\mu^{n}|^2d\bm{x}
          \leq 0.
	\end{equation}
   We arrive at \eqref{sec3_eq2} by further noting that $R^n>0$ for all $n$.
\end{proof}

Lemma \ref{sec3_Lemma2} implies that
there exists a constant $M$, depending only on $\Omega$, $\phi_{in}$ and $c_0$,
such that for any $n$,
\begin{equation}\label{sec3_eq4}
R^{n}\leq M.
\end{equation}
Note that $c_0$ in \eqref{energy_0}
is a chosen constant to ensure $E(t)>0$ for $0\leq t\leq T$.
We can choose $c_0$ such that $E(t)\geq C_0$, for some constant $C_0>0$.
It then follows from equation \eqref{xi} that $\xi_{_{1A}}^{n+1}$ is bounded
from above, since
\begin{equation}\label{sec3_eq5}
\xi_{_{1A}}^{n+1}=\frac{R^n}{\sqrt{\bm{E}[\phi^{n}]}+\frac{\Delta t}{2\sqrt{\bm{E}[\phi^{n}]}}\int_{\Omega}|\nabla \mu^n|^2d\bm{x}}\leq \frac{R^n}{\sqrt{\bm{E}[\phi^{n}]}}\leq \frac{M}{\sqrt{C_0}}.
\end{equation}

\begin{myTheorem}\label{sec3_Lemma3}
  Suppose $\phi_{in}\in H^3(\Omega)$ and the condition \eqref{Assume_H2} holds.
  The following inequality holds for all $n$
  with the scheme~\eqref{CH1_weak}, 
        \begin{align*}
        \|\nabla\phi^{n+1}\|_0^2+\frac{\lambda}{2}\|\phi^{n+1}\|_0^2
        +\frac{\Delta t}{2}\|\nabla\Delta\phi^{n+1}\|_0^2+\lambda\Delta t\sum_{k=0}^{n}\|\Delta\phi^{k+1}\|_0^2+\frac{\Delta t}{2}\sum_{k=0}^{n}\|\nabla\mu^{k+1}\|_0^2\leq \widehat{C}_1,
        \end{align*}
        where $\widehat{C}_1=\exp\left(C(M)T\right)\left(\|\nabla\phi^{0}\|_0^2+\frac{\lambda}{2}\|\phi^{0}\|_0^2+\frac{\Delta t}{2}\|\nabla\Delta\phi^{0}\|_0^2\right)$,
        and $C(M)$ is a constant depending on $M$.
\end{myTheorem}

\begin{proof}
  Taking the $L^2$ inner product between \eqref{CH1_eq1} and $\Delta t \mu^{n+1}$
  and between \eqref{CH1_eq2}  and $-(\phi^{n+1}-\phi^{n})$, and summing up
  the two resultant equations, we have
	\begin{align}\label{sec3_eq6}
	&\frac{1}{2}(\|\nabla\phi^{n+1}\|_0^2-\|\nabla\phi^{n}\|_0^2+\|\nabla\phi^{n+1}-\nabla\phi^{n}\|_0^2)+\frac{\lambda}{2}(\|\phi^{n+1}\|_0^2-\|\phi^{n}\|_0^2
	\nonumber\\
	&\qquad+\|\phi^{n+1}-\phi^{n}\|_0^2)+\Delta t\|\nabla\mu^{n+1}\|_0^2=-|\xi_{_{1A}}^{n+1}|^2 (h(\phi^{n}),\phi^{n+1}-\phi^{n}),
	\end{align}
        where the boundary condition \eqref{boundary_1} has been used.
     In light of \eqref{CH1_eq1}, we have
	\begin{align}\label{eq_a}
	  -|\xi_{_{1A}}^{n+1}|^2 (h(\phi^{n}),\phi^{n+1}-\phi^{n})=-|\xi_{_{1A}}^{n+1}|^2\Delta t(h(\phi^{n}),\Delta\mu^{n+1})=|\xi_{_{1A}}^{n+1}|^2\Delta t(\nabla h(\phi^{n}),\nabla\mu^{n+1})
	  \nonumber\\
          \leq\frac{\Delta t}{2}\|\nabla\mu^{n+1}\|_0^2+\frac{|\xi_{_{1A}}^{n+1}|^4\Delta t}{2}\|\nabla h(\phi^{n})\|_0^2.
	\end{align}
        Taking the $L^2$ inner product between  \eqref{CH1_eq1} and $\Delta\phi^{n+1}$
        and between \eqref{CH1_eq2} and $\Delta^2\phi^{n+1}$, and
        summing up the resultant equations, we arrive at
	\begin{align}\label{sec3_eq7}
	&\frac{1}{2}(\|\nabla\phi^{n+1}\|_0^2-\|\nabla\phi^{n}\|_0^2+\|\nabla(\phi^{n+1}-\phi^{n})\|_0^2)+\Delta t\|\nabla\Delta\phi^{n+1}\|_0^2+\lambda\Delta t\|\Delta\phi^{n+1}\|_0^2
	\nonumber\\
	&\qquad =\Delta t|\xi_{_{1A}}^{n+1}|^2(\nabla h(\phi^{n}),\nabla\Delta\phi^{n+1})\leq \frac{\Delta t}{2}\|\nabla\Delta\phi^{n+1}\|_0^2+\frac{|\xi_{_{1A}}^{n+1}|^4\Delta t}{2}\|\nabla h(\phi^{n})\|_0^2.
	\end{align}
        By incorporating \eqref{eq_a} into \eqref{sec3_eq6} and
        summing up equations \eqref{sec3_eq6} and \eqref{sec3_eq7}, we get
	\begin{align}\label{sec3_eq8}
	&\|\nabla\phi^{n+1}\|_0^2-\|\nabla\phi^{n}\|_0^2+\|\nabla\phi^{n+1}-\nabla\phi^{n}\|_0^2+\frac{\lambda}{2}(\|\phi^{n+1}\|_0^2-\|\phi^{n}\|_0^2+\|\phi^{n+1}-\phi^{n}\|_0^2)
	\nonumber\\
	&\qquad +\frac{\Delta t}{2}\|\nabla\mu^{n+1}\|_0^2+\frac{\Delta t}{2}\|\nabla\Delta\phi^{n+1}\|_0^2+\lambda\Delta t\|\Delta\phi^{n+1}\|_0^2\leq |\xi_{_{1A}}^{n+1}|^4\Delta t\|\nabla h(\phi^{n})\|_0^2.
	\end{align}
        
	To deal with the term on the right hand side of \eqref{sec3_eq8},
        we use an idea from \cite{Temam1997Infinite,Shen2018Convergence}.
        Noting that $\nabla h(\phi^{n}) = h'(\phi^{n})\nabla\phi^{n}$,
        equation \eqref{def_h} and the relation \eqref{Assume_H2}, we have
	\begin{equation}\label{sec3_eq9}
	\|\nabla h(\phi^{n})\|_0^2\leq \|h'(\phi^{n})\|_{0,\infty}^2\|\nabla\phi^{n}\|_0^2\leq C\left(\|\nabla\phi^{n}\|_0^2+\|\nabla\phi^{n}\|_0^2\|\phi^{n}\|_{0,\infty}^{4}\right).
	\end{equation}
	Let $\breve{\phi}^n= \frac{1}{|\Omega|}\int_{\Omega}\phi^nd\bm{x}$.
        Lemma \ref{sec3_Lemma1} implies that
        \begin{equation}
        |\breve{\phi}^n|^2=
        |\breve{\phi}^0|^2\leq \frac{1}{|\Omega|}\|\phi^0\|_0^2\leq C. 
        \end{equation}
	Using Sobolev embedding theorems $H^1(\Omega) \hookrightarrow L^{\infty}(\Omega)$ ($d=1$), $H^{1+2\sigma}(\Omega) \hookrightarrow L^{\infty}(\Omega)$ for any $\sigma>0$ ($d=2$), the Agmon's inequality ($d=3$, see \cite{Temam1997Infinite,Shen2018Convergence} for more details) and the interpolation inequality about the spaces $H^s(\Omega)$, we deduce that
	\begin{eqnarray}\label{sec3_eq10}
	\|\phi^n-\breve{\phi}^n\|_{0,\infty}\leq\left\{\begin {array}{lll}
	C\|\nabla\phi^n\|_0        & \mbox{for}\ d =1,\\
	C\|\nabla\phi^n\|_0^{1-\sigma}\|\nabla\Delta\phi^n\|_0^{\sigma}      & \mbox{for}\ d=2,\\
	C\|\nabla\phi^n\|_0^{3/4}\|\nabla\Delta\phi^n\|_0^{1/4}        & \mbox{for}\ d=3.	\end{array}\right.
	\end{eqnarray}
    By the triangle inequality, we have
    \[
	\|\phi^n\|_{0,\infty}^{4}\leq 8\left(\|\phi^n-\breve{\phi}^n\|_{0,\infty}^4+\|\breve{\phi}^n\|_{0,\infty}^4\right).
	\]
	By setting $\sigma=1/4$ and using $\|\breve{\phi}^n\|_{0,\infty}^4\leq C$,
        we arrive at the relation
	\begin{eqnarray*}
    \|\nabla\phi^n\|_0^2\|\phi^n\|_{0,\infty}^4\leq\left\{\begin {array}{lll}
	C\|\nabla\phi^n\|_0^2+C\|\nabla\phi^n\|_0^6       & \mbox{for}\ d =1,\\
	C\|\nabla\phi^n\|_0^2+C\|\nabla\phi^n\|_0^{5}\|\nabla\Delta\phi^n\|_0     & \mbox{for}\ d=2,\\
	C\|\nabla\phi^n\|_0^2+C\|\nabla\phi^n\|_0^{5}\|\nabla\Delta\phi^n\|_0        & \mbox{for}\ d=3.	\end{array}\right.
	\end{eqnarray*}
        
	Applying  the following inequality
    \begin{equation}\label{sec3_xi}
	\xi_{_{1A}}^{n+1}=\frac{R^n}{\sqrt{\bm{E}[\phi^{n}]}+\frac{\Delta t}{2\sqrt{\bm{E}[\phi^{n}]}}\int_{\Omega}|\nabla \mu^n|_0^2d\bm{x}}\leq\frac{C(M)}{\|\phi^n\|_1}
	\end{equation}
	and \eqref{sec3_eq5}, i.e. $\xi_{_{1A}}^{n+1}\leq C(M)$, we obtain
	\begin{eqnarray*}
	|\xi_{_{1A}}^{n+1}|^4\|\nabla\phi^{n}\|_0^2\|\phi^{n}\|_{0,\infty}^{4}\leq\left\{\begin {array}{lll}
	C(M)\|\nabla\phi^n\|_0^2+\frac{C(M)}{\|\phi^n\|_1^4}\|\nabla\phi^n\|_0^6     & \mbox{for}\ d =1,\\
	C(M)\|\nabla\phi^n\|_0^2 +\frac{C(M)}{\|\phi^n\|_1^4}\|\nabla\phi^n\|_0^5\|\nabla\Delta\phi^n\|_0     & \mbox{for}\ d=2,\\
	C(M)\|\nabla\phi^n\|_0^2 +\frac{C(M)}{\|\phi^n\|_1^4}\|\nabla\phi^n\|_0^5\|\nabla\Delta\phi^n\|_0      & \mbox{for}\ d=3.	\end{array}\right.
	\end{eqnarray*}
	By using the Cauchy Schwarz inequality, for any $\epsilon_i>0$ ($i$=$1$, $2$), there exist constants $C(\epsilon_i, M)$ depending on $\epsilon_i$ and $M$, such that
	\begin{eqnarray}\label{sec3_eq11}
	|\xi_{_{1A}}^{n+1}|^4\|\nabla h(\phi^{n})\|_0^2\leq\left\{\begin {array}{lll}
	C(M)\|\nabla\phi^n\|_0^2      & \mbox{for}\ d =1,\\
	C(M)\|\nabla\phi^n\|_0^2+C(\epsilon_1, M)\|\nabla\phi^n\|_0^2+\epsilon_1\|\nabla\Delta\phi^n\|_0^2         & \mbox{for}\ d=2,\\
	C(M)\|\nabla\phi^n\|_0^2+C(\epsilon_2, M)\|\nabla\phi^n\|_0^2+\epsilon_2\|\nabla\Delta\phi^n\|_0^2        & \mbox{for}\ d=3.	\end{array}\right.
	\end{eqnarray}
	We set $\epsilon_1=\epsilon_2=\frac{1}{2}$ and combine
        the above inequalities with \eqref{sec3_eq8}, and then
	\begin{align}\label{sec3_eq12}
	&\|\nabla\phi^{n+1}\|_0^2-\|\nabla\phi^{n}\|_0^2+\frac{\lambda}{2}(\|\phi^{n+1}\|_0^2-\|\phi^{n}\|_0^2)+\frac{\Delta t}{2}\|\nabla\mu^{n+1}\|_0^2
	\nonumber\\
	&\qquad +\frac{\Delta t}{2}(\|\nabla\Delta\phi^{n+1}\|_0^2-\|\nabla\Delta\phi^{n}\|_0^2)+\lambda\Delta t\|\Delta\phi^{n+1}\|_0^2\leq C(M)\Delta t\|\nabla\phi^{n}\|_0^2.
	\end{align}
        We conclude the proof by 
	summing up the above relation for indices from 
        $0$ to $n$ and by using the discrete Gronwall lemma \ref{sec2_Lemma4}.
\end{proof}

\begin{myTheorem}\label{sec3_Lemma4}
  Suppose $\phi_{in}\in H^4(\Omega)$, and that
  the conditions for Lemmas \ref{sec2_Lemma1} and \ref{sec2_Lemma2} hold.
  The following inequality holds for
  all $n$ with the scheme~\eqref{CH1_weak},
	\begin{align*}
	\|\Delta\phi^{n+1}\|_0^2+\frac{\Delta t}{2}\|\Delta^2\phi^{n+1}\|_0^2+\frac{\Delta t}{2}\sum_{k=0}^{n}\|\Delta^2\phi^{k+1}\|_0^2+2\lambda\Delta t\sum_{k=0}^{n}\|\nabla\Delta\phi^{k+1}\|_0^2\leq\widehat{C}_2,
	\end{align*}
	where $\widehat{C}_2=\|\Delta\phi^{0}\|_0^2+\frac{\Delta t}{2}\|\Delta^2\phi^{0}\|_0^2+C(M)T$.
\end{myTheorem}
\begin{proof} 
  The proof is similar to that for the SAV scheme in \cite{Shen2018Convergence},
  by using the Lemmas \ref{sec2_Lemma1} and \ref{sec2_Lemma2}.
\end{proof}

\paragraph{Error Estimate}
%
We next examine the errors of the solution to the Cahn-Hilliard
equation with the scheme \eqref{CH1_weak}.
Let 
\begin{equation} \label{def_error}
e_{\phi}^n=\phi^n-\phi(t^n),\quad  e_{\mu}^n=\mu^n-\mu(t^n)\quad  {\rm and} \quad e_{R}^n=R^n-R(t^n).
\end{equation}
Assume that $\phi_{in} \in H^4(\Omega)$
and the solution ($\phi,\mu$) of equations \eqref{CH}-\eqref{boundary} satisfies
\begin{align}\label{regularity_1}
&\phi \in L^{\infty}(0,T;W^{3,\infty}(\Omega)), \qquad  \phi_t \in L^{4}(0,T;H^{1}(\Omega)) \cap L^{2}(0,T;H^1(\Omega)),
\nonumber\\
&\qquad\qquad \phi_{tt} \in L^{2}(0,T;H^{-1}(\Omega)), \qquad \mu \in L^{\infty}(0,T;H^{1}(\Omega)).
\end{align}
In light of Lemma \ref{sec2_Lemma1} and Theorems \ref{sec3_Lemma3} and \ref{sec3_Lemma4},
we conclude that
\begin{equation}\label{sec4_eq1}
\|\phi(t)\|_{2} \leq C, \quad \|\phi^n\|_{2}\leq C,
\end{equation}
where the constant $C$ is dependent on $T$, $\phi_{in}$, and $\Omega$.
Since $H^2(\Omega) \hookrightarrow L^{\infty}(\Omega)$, we conclude that 
\begin{equation}\label{sec4_eq2}
|h(\phi)|, |h'(\phi)|, |h''(\phi)|, |h(\phi^n)|, |h'(\phi^n)|, |h''(\phi^n)|\leq C.
\end{equation}
Based on the relation $R(t)=\sqrt{\bm{E}[\phi]}$
and the equation \eqref{energy_0}, we have
\begin{align}\label{sec4_eq3}
\frac{d^2R}{dt^2} &= -\frac{1}{4\sqrt{\bm{E}[\phi]^{3}}}\left(\int_{\Omega}\big{(}\nabla \phi\cdot\nabla\phi_t+\lambda\phi\phi_{t}+h(\phi)\phi_t\big{)}d\bm{x}\right)^2
\nonumber\\
&\ \ +\frac{1}{2\sqrt{\bm{E}[\phi]}}\int_{\Omega}\left(|\nabla\phi_t|^2
+\nabla\phi\cdot\nabla\phi_{tt}+\lambda|\phi_{t}|^2+\lambda\phi\phi_{tt}+h'(\phi)|\phi_{t}|^2+h(\phi)\phi_{tt}\right)d\bm{x}.
\end{align}
Combining \eqref{boundary}, \eqref{regularity_1}, \eqref{sec4_eq1}
and \eqref{sec4_eq2} with \eqref{sec4_eq3}, we deduce that
\begin{equation}\label{sec4_eq4}
\int_{0}^T\left|\frac{d^2R}{dt^2}\right|^2dt\leq C\int_{0}^T\left(\|\phi_{t}\|_{1}^4+\|\phi_{t}\|_{1}^2+\|\phi_{tt}\|_{-1}^2\right)dt,
\end{equation}
where $\phi_{tt}$ denotes the second time derivative of $\phi$.

The truncation errors $T_{\phi_{1A}}^{n+1}$ and $T_{R_{1A}}^{n+1}$ are defined by
\begin{subequations}\label{truncation1A}
	\begin{align}
	\label{truncation1A_eq1}
	&\frac{\phi(t^{n+1})-\phi(t^{n})}{\Delta t} = \Delta\mu(t^{n+1})+\frac{1}{\Delta t}T_{\phi_{1A}}^{n+1},\\
	\label{truncation1A_eq2}
	&\mu(t^{n+1})=-\Delta\phi(t^{n+1})+ \lambda\phi(t^{n+1}) + \frac{R(t^{n+1})^2}{\bm{E}[\phi(t^{n+1})]}h(\phi(t^{n+1})),\\
	\label{truncation1A_eq3}
	&\frac{R(t^{n+1})-R(t^{n})}{\Delta t} = -\frac{R(t^{n})}{2\bm{E}[\phi(t^{n})]} \int_{\Omega}|\nabla\mu(t^{n})|^2d\bm{x}+\frac{1}{\Delta t}T_{R_{1A}}^{n+1},
	\end{align}
\end{subequations}
where
\begin{equation}\label{sec4_eq5}
T_{\phi_{1A}}^{n+1} = \int_{t^n}^{t^{n+1}}(t^n-t)\phi_{tt}(t)dt
\quad {\rm and} \quad T_{R_{1A}}^{n+1} = \int_{t^n}^{t^{n+1}}(t^{n+1}-t)\frac{d^2R(t)}{dt^2}dt.
\end{equation}

With the above definitions and relations,
the errors of the scheme \eqref{CH1_weak} is summarized by the following result.
\begin{myTheorem}\label{sec4_Lemma4}  
  Suppose the condition \eqref{regularity_1},
  and the conditions for Theorems \ref{sec3_Lemma3} and \ref{sec3_Lemma4} hold. 
  The following result holds with sufficiently small $\Delta t$,
  \begin{equation}\label{sec4_eq6}
	\frac{1}{2}\|\nabla e_{\phi}^{n+1}\|_0^2+\frac{\lambda}{2}\|e_{\phi}^{n+1}\|_0^2+\frac{\Delta t}{2}\|\nabla e_{\mu}^{n+1}\|_0^2+|e_{R}^{n+1}|^2\leq\widehat{C}_3\Delta t^2,
	\end{equation}
	where $\widehat{C}_3=C\exp(\Delta t\sum_{k=0}^{n}\frac{r^{k}}{1-r^{k}\Delta t})\int_{0}^{t^{n+1}}\left(\|\phi_t(s)\|_{1}^4+\|\phi_t(s)\|_{1}^2+\|\phi_{tt}(s)\|_{-1}^2\right)ds$, $r^k=1+\|\nabla\mu^k\|_0^2$ and the constant $C$ is dependent on $T$, $\phi_{in}$, $\Omega$, $\|\phi\|_{L^{\infty}(0,T;W^{3,\infty}(\Omega))}$ and $\|\mu\|_{L^{\infty}(0,T;H^{1}(\Omega))}$. 
\end{myTheorem}

\begin{proof} By subtracting \eqref{truncation1A} from \eqref{CH1_weak}, we have
	\begin{subequations}\label{error1A}
		\begin{align}
		\label{error1A_eq1}
		&\frac{e_{\phi}^{n+1}-e_{\phi}^{n}}{\Delta t}= \Delta e_{\mu}^{n+1}-\frac{1}{\Delta t}T_{\phi_{1A}}^{n+1},\\
		\label{error1A_eq2}
		&e_{\mu}^{n+1}=-\Delta e_{\phi}^{n+1}+ \lambda e_{\phi}^{n+1} + A_1^{n+1},\\
		\label{error1A_eq3}
		&\frac{e_{R}^{n+1}-e_{R}^{n}}{\Delta t}= -\frac{1}{2}A_2^{n+1}-\frac{1}{\Delta t}T_{R_{1A}}^{n+1},
		\end{align}
	\end{subequations}
	where 
	\begin{align*}
	A_1^{n+1}&=|\xi_{_{1A}}^{n+1}|^2h(\phi^{n})-\frac{R(t^{n+1})^2}{\bm{E}[\phi(t^{n+1})]}h(\phi(t^{n+1}))\\
	&=\frac{e_R^{n+1}(R^{n+1}+R(t^{n+1}))}{\bm{E}[\phi^{n}]}h(\phi^{n})
	+R(t^{n+1})^2\left(\frac{h(\phi^{n})}{\bm{E}[\phi^{n}]}-\frac{h(\phi(t^{n}))}{\bm{E}[\phi(t^{n})]}\right)\\
	&\quad+R(t^{n+1})^2\left(\frac{h(\phi(t^{n}))}{\bm{E}[\phi(t^{n})]}-\frac{h(\phi(t^{n+1}))}{\bm{E}[\phi(t^{n+1})]}\right),\\
	A_2^{n+1}&=\frac{\xi_{_{1A}}^{n+1}}{\sqrt{\bm{E}[\phi^{n}]}}\int_{\Omega}|\nabla\mu^{n}|^2d\bm{x}-\frac{R(t^{n})}{\bm{E}[\phi(t^{n})]} \int_{\Omega}|\nabla\mu(t^{n})|^2d\bm{x}\\
	&=\frac{e_{R}^{n+1}}{\bm{E}[\phi^{n}]}\int_{\Omega}|\nabla\mu^{n}|^2d\bm{x}+
	\frac{R(t^{n+1})}{\bm{E}[\phi^{n}]}\int_{\Omega}(|\nabla\mu^{n}|^2-|\nabla\mu(t^{n})|^2)d\bm{x}\\
	& \quad + R(t^{n+1})\left(\frac{1}{\bm{E}[\phi^{n}]}-\frac{1}{\bm{E}[\phi(t^{n})]}\right)\int_{\Omega}|\nabla\mu(t^{n})|^2d\bm{x}+\frac{R(t^{n+1})-R(t^{n})}{\bm{E}[\phi(t^{n})]}\int_{\Omega}|\nabla\mu(t^{n})|^2d\bm{x}.
	\end{align*}
        
	Taking the inner product between \eqref{error1A_eq1} and $\Delta te_{\mu}^{n+1}$
        and between \eqref{error1A_eq2} and $e_{\phi}^{n+1}-e_{\phi}^{n}$,
        multiplying \eqref{error1A_eq3} by $2\Delta te_{R}^{n+1}$, and
        combining the resultant equations, we get
	\begin{subequations}\label{sec4_eq7}
		\begin{align}
		\label{sec4_7_eq1}
		&\frac{1}{2}(\|\nabla e_{\phi}^{n+1}\|_0^2-\|\nabla e_{\phi}^{n}\|_0^2+\|\nabla e_{\phi}^{n+1}-\nabla e_{\phi}^{n}\|_0^2)+\frac{\lambda}{2}(\|e_{\phi}^{n+1}\|_0^2-\|e_{\phi}^{n}\|_0^2
		\nonumber\\
		&\qquad+\|e_{\phi}^{n+1}-e_{\phi}^{n}\|_0^2)+\Delta t\|\nabla e_{\mu}^{n+1}\|_0^2=-(T_{\phi_{1A}}^{n+1},e_{\mu}^{n+1})- (A_1^{n+1},e_{\phi}^{n+1}-e_{\phi}^{n}),\\
		\label{sec4_7_eq2}
		&|e_{R}^{n+1}|^2-|e_{R}^{n}|^2+|e_{R}^{n+1}-e_{R}^{n}|^2= -\Delta tA_2^{n+1}e_{R}^{n+1}-2e_{R}^{n+1}T_{R_{1A}}^{n+1}.
		\end{align}
	\end{subequations}
	By the Taylor expansion theorem, we deal with the truncation errors as follows,
	\begin{align*}
	-(T_{\phi_{1A}}^{n+1},e_{\mu}^{n+1})&\leq\frac{\Delta t}{8}\|\nabla e_{\mu}^{n+1}\|_0^2+\frac{2}{\Delta t}\|(-\Delta)^{-1/2}T_{\phi_{1A}}^{n+1}\|_0^2\\
	&\leq\frac{\Delta t}{8}\|\nabla e_{\mu}^{n+1}\|_0^2+C\Delta t^2\int_{t^{n}}^{t^{n+1}}\|\phi_{tt}(s)\|_{-1}^2ds,\\
	-2e_{R}^{n+1}T_{R_{1A}}^{n+1}&\leq\Delta t|e_{R}^{n+1}|^2+\frac{1}{\Delta t}|T_{R_{1A}}^{n+1}|^2\leq \Delta t|e_R^{n+1}|^2+
	C\Delta t^2\int_{t^{n}}^{t^{n+1}}\left|\frac{d^2R(s)}{dt^2}\right|^2ds,
	\end{align*}
        where $(-\Delta)^{-1/2}$ denotes the power of
        $-\Delta$ by the spectral theory of self-adjoint operators.
	We treat the $A_1^{n+1}$ term on the right-hand side
        of \eqref{sec4_7_eq1} as follows. 
	\begin{equation}\label{sec4_eq8}
          \begin{split}
	  &-e_R^{n+1}(R^{n+1}+R(t^{n+1}))
          \left(\frac{h(\phi^{n})}{\bm{E}[\phi^{n}]},e_{\phi}^{n+1}-e_{\phi}^{n}\right) \\
          &=e_R^{n+1}(R^{n+1}+R(t^{n+1}))\left(\frac{h(\phi^{n})}{\bm{E}[\phi^{n}]},\Delta t\Delta e_{\mu}^{n+1}-T_{\phi_{1A}}^{n+1}\right) \\
	&\leq Ce_{R}^{n+1}\left(\Delta t\|\nabla e_{\mu}^{n+1}\|_0+\|(-\Delta)^{-1/2}T_{\phi_{1A}}^{n+1}\|_0\right)\bigg{\|}\frac{\nabla h(\phi^{n})}{\bm{E}[\phi^{n}]}\bigg{\|}_0 \\
	&= Ce_{R}^{n+1}\left(\Delta t\|\nabla e_{\mu}^{n+1}\|_0+\|(-\Delta)^{-1/2}T_{\phi_{1A}}^{n+1}\|_0\right)\bigg{\|}\frac{ h'(\phi^{n})\nabla\phi^{n}}{\bm{E}[\phi^{n}]}\bigg{\|}_0 \\
	&\leq\frac{\Delta t}{8}\|\nabla e_{\mu}^{n+1}\|_0^2+C\Delta t|e_{R}^{n+1}|^2+
	  C\Delta t^2\int_{t^{n}}^{t^{n+1}}\|\phi_{tt}(s)\|_{-1}^2ds.
          \end{split}
	\end{equation}
	Additionally, 
	\begin{equation*}
          \begin{split}
	  &-R(t^{n+1})^2\left(\frac{h(\phi^{n})}{\bm{E}[\phi^{n}]}-\frac{h(\phi(t^{n}))}{\bm{E}[\phi(t^{n})]},e_{\phi}^{n+1}-e_{\phi}^{n}\right) \\
	&\leq\frac{\Delta t}{8}\|\nabla e_{\mu}^{n+1}\|_0^2+C\Delta t\Big{\|}\frac{\nabla h(\phi^{n})}{\bm{E}[\phi^{n}]}-\frac{\nabla h(\phi(t^{n}))}{\bm{E}[\phi(t^{n})]}\Big{\|}_0^2
	    +C\Delta t^2\int_{t^{n}}^{t^{n+1}}\|\phi_{tt}(s)\|_{-1}^2ds,
            \end{split}
	\end{equation*}
        \begin{equation*}
          \begin{split}
	&-R(t^{n+1})^2\left(\frac{h(\phi(t^{n}))}{\bm{E}[\phi(t^{n})]}-\frac{h(\phi(t^{n+1}))}{\bm{E}[\phi(t^{n+1})]},e_{\phi}^{n+1}-e_{\phi}^{n}\right) \\
	&\leq\frac{\Delta t}{8}\|\nabla e_{\mu}^{n+1}\|_0^2+C\Delta t\Big{\|}\frac{\nabla h(\phi(t^{n}))}{ \bm{E}[\phi(t^{n})]}-\frac{\nabla h(\phi(t^{n+1}))}{\bm{E}[\phi(t^{n+1})]}\Big{\|}_0^2
	+C\Delta t^2\int_{t^{n}}^{t^{n+1}}\|\phi_{tt}(s)\|_{-1}^2ds.
        \end{split}
	\end{equation*}

        Note that $\int_{\Omega}\phi(t)d\bm{x}$ is a constant
        and $\int_{\Omega}T_{\phi_{1A}}^{n+1}d\bm{x}=0$.
        Noting the definition of $\bm{E}[\phi]$ and
        that $H(s) \in C^3(\mathbb{R})$, we have
	\begin{equation}\label{sec4_eq9}
          \begin{split}
	    &\bm{E}[\phi^{n}]-\bm{E}[\phi(t^{n})] \\
            &= \frac{1}{2}\int_{\Omega}(\nabla\phi^n+\nabla \phi(t^n))\nabla e_\phi^nd\bm{x}+\frac{\lambda}{2}\int_{\Omega}(\phi^n+\phi(t^n))e_\phi^nd\bm{x}+\int_{\Omega}\big{(}H(\phi^n)-H(\phi(t^n)\big{)}d\bm{x} \\
	&\leq C \|\nabla e_\phi^n\|_0+C\|e_\phi^n\|_0+\int_{\Omega}H'\big{(}\theta\phi^n+(1-\theta)\phi(t^{n})\big{)}\big{(}\phi^n-\phi(t^{n})\big{)}d\bm{x} \\
	&\leq C \|\nabla e_\phi^n\|_0+C\|e_\phi^n\|_0.
        \end{split}
	\end{equation}
	We rewrite the term
        $\frac{\nabla h(\phi^{n})}{\bm{E}[\phi^{n}]}-\frac{\nabla h(\phi(t^{  n}))}{\bm{E}[\phi(t^{n})]}$  into
	\begin{equation*}
	\frac{\nabla h(\phi^{n})}{\bm{E}[\phi^{n}]}-\frac{\nabla h(\phi(t^{  n}))}{\bm{E}[\phi(t^{n})]}=\frac{\nabla h(\phi^{n})-\nabla h(\phi(t^{n}))}{\bm{E}[\phi^{n}]}
	+\frac{\nabla h(\phi(t^{n}))\big{(}\bm{E}[\phi(t^{n})]-\bm{E}[\phi^{n}]\big{)}}{  \bm{E}[\phi^{n}]\bm{E}[\phi(t^{n})]}.
	\end{equation*}
	It follows from the H${\rm\ddot{o}}$lder's inequality and Sobolev embedding theorem
        that,
        \begin{align*}
	\|\nabla h(\phi^{n})-\nabla h(\phi(t^{n}))\|_0\leq&\|(h'(\phi^{n})-h'(\phi(t^{n})))\nabla\phi(t^{n})\|_0+\|h'(\phi^{n})\nabla e_{\phi}^n\|_0
	\nonumber\\
	\leq&C\|\nabla\phi(t^{n})e_{\phi}^n\|_{0}+C\|\nabla e_{\phi}^n\|_0
	\leq C(\|\nabla\phi(t^{n})\|_{0,3}\|e_{\phi}^n\|_{0,6}+\|\nabla e_{\phi}^n\|_0)
	\nonumber\\
	\leq&C\|\phi(t^{n})\|_{2}\|e_{\phi}^n\|_{1}+C\|\nabla e_{\phi}^n\|_0\leq C(\|\nabla e_{\phi}^n\|_0+\|e_{\phi}^n\|_0).
	\end{align*}
	Then, 
	\begin{align}\label{sec4_eq11}
	&\Big{\|}\frac{\nabla h(\phi^{n})}{\bm{E}[\phi^{n}]}-\frac{\nabla h(\phi(t^{  n}))}{\bm{E}[\phi(t^{n})]} \Big{\|}_0^2
	\nonumber\\
	=&\Bigg{\|}\frac{\nabla h(\phi^{n})-\nabla h(\phi(t^{n}))}{\bm{E}[\phi^{n}]}
	+\frac{\nabla h(\phi(t^{n}))\big{(}\bm{E}[\phi(t^{n})]-\bm{E}[\phi^{n}]\big{)}}{  \bm{E}[\phi^{n}]\bm{E}[\phi(t^{n})]}\Bigg{\|}_0^2
	\nonumber\\
	\leq&C\|\nabla h(\phi^{n})-\nabla h(\phi(t^{n}))\|_0^2+C\|\nabla h(\phi(t^{n}))\|_0^2\big{|}\bm{E}[\phi^{n}]-\bm{E}[\phi(t^{n})]\big{|}^2
	\nonumber\\
	\leq&C(\|\nabla e_{\phi}^n\|_0^2+\|e_{\phi}^n\|_0^2).
	\end{align}
	Similarly, 
	\begin{align}\label{sec4_eq12}
	&\Big{\|}\frac{\nabla h(\phi(t^{n}))}{\bm{E}[\phi(t^{n})]}-\frac{\nabla h(\phi(t^{n+1}))}{\bm{E}[\phi(t^{n+1})]}\Big{\|}_0^2
	\nonumber\\
	\leq&C\|\nabla h(\phi(t^{n}))-\nabla h(\phi(t^{n+1}))\|_0^2+C\|\nabla h(\phi(t^{n+1}))\|_0^2\big{|}\bm{E}[\phi(t^{n+1})]-\bm{E}[\phi(t^{n})]\big{|}^2
	\nonumber\\
	\leq&C(\|\nabla\phi(t^{n+1})-\nabla\phi(t^{n})\|_0^2+\|\phi(t^{n+1})-\phi(t^{n})\|_0^2)\nonumber\\
	\leq&C\Delta t\int_{t^n}^{t^{n+1}}\|\phi_t(s)\|_{1}^2ds.
	\end{align}

	Next, we treat the right-hand side of \eqref{sec4_7_eq2} as follows: 
	\begin{equation}\label{sec4_eq13}
          \left\{
          \begin{split}
	    &-\frac{\Delta t|e_{R}^{n+1}|^2}{\bm{E}[\phi^{n}]}\int_{\Omega}|\nabla\mu^{n}|^2d\bm{x}
            \leq C\Delta t\|\nabla\mu^{n}\|_0^2|e_{R}^{n+1}|^2,\\
	    &-\frac{\Delta te_{R}^{n+1}R(t^{n+1})}{\bm{E}[\phi^{n}]}\int_{\Omega}(|\nabla\mu^{n}|^2-|\nabla\mu(t^{n})|^2)d\bm{x} \\
            &\qquad\qquad
            \leq\frac{\Delta te_{R}^{n+1}R(t^{n+1})}{\bm{E}[\phi^{n}]}\|\nabla e_\mu^{n}\|_0\|\nabla\mu^{n}+\nabla\mu(t^{n})\|_0 \\
	    &\qquad\qquad
            \leq C\Delta t(\|\nabla\mu^{n}\|_0^2+1)|e_{R}^{n+1}|^2+\frac{\Delta t}{2}\|\nabla e_\mu^{n}\|_0^2, \\
	    &-\Delta tR(t^{n+1})e_{R}^{n+1}\left(\frac{1}{\bm{E}[\phi^{n}]}-\frac{1}{\bm{E}[\phi(t^{n})]}\right)\int_{\Omega}|\nabla\mu(t^{n})|^2d\bm{x} \\
            &\qquad\qquad
            \leq C\Delta t\|\nabla\mu(t^{n})\|_0^2\big{(}|e_{R}^{n+1}|^2+\left|\bm{E}[\phi^{n}]-\bm{E}[\phi(t^{n})]\right|^2\big{)} \\
	    &\qquad\qquad
            \leq C\Delta t\left(|e_{R}^{n+1}|^2+\|\nabla e_{\phi}^n\|_0^2+\|e_{\phi}^n\|_0^2\right), \\
	    &-\Delta te_{R}^{n+1}\frac{R(t^{n+1})-R(t^{n})}{\bm{E}[\phi(t^{n})]}\int_{\Omega}|\nabla\mu(t^{n})|^2d\bm{x} \\
            &\qquad\qquad
            \leq C\Delta t\|\nabla\mu(t^{n})\|_0^2\left(|e_{R}^{n+1}|^2+\big{|}R(t^{n+1})-R(t^{n})\big{|}^2\right) \\
	    &\qquad\qquad
            \leq C\Delta t|e_{R}^{n+1}|^2+C\Delta t^2\int_{t^n}^{t^{n+1}}\left|\frac{dR(s)}{dt}\right|^2ds.
          \end{split}
          \right.
	\end{equation}
	By combining the above inequalities with \eqref{sec4_7_eq1} and \eqref{sec4_7_eq2},  we have
	\begin{align}\label{sec4_10}
	&\frac{1}{2}(\|\nabla e_{\phi}^{n+1}\|_0^2-\|\nabla e_{\phi}^{n}\|_0^2)+\frac{\lambda}{2}(\|e_{\phi}^{n+1}\|_0^2-\|e_{\phi}^{n}\|_0^2) +|e_{R}^{n+1}|^2-|e_{R}^{n}|^2+\frac{\Delta t}{2}(\|\nabla e_{\mu}^{n+1}\|_0^2-\|\nabla e_\mu^{n}\|_0^2)
	\nonumber\\
	&\qquad+\frac{1}{2}\|\nabla e_{\phi}^{n+1}-\nabla e_{\phi}^{n}\|_0^2+\frac{\lambda}{2}\|e_{\phi}^{n+1}-e_{\phi}^{n}\|_0^2+|e_{R}^{n+1}-e_{R}^{n}|^2
	\nonumber\\
	&\leq C\Delta t(1+\|\nabla\mu^n\|_0^2)|e_{R}^{n+1}|^2+C\Delta t\big{(}\|\nabla e_{\phi}^n\|_0^2+\|e_{\phi}^n\|_0^2\big{)} +C\Delta t^2\int_{t^n}^{t^{n+1}}\|\phi_t(s)\|_1^2ds
	\nonumber\\
	&\ \ +C\Delta t^2\int_{t^{n}}^{t^{n+1}}\|\phi_{tt}(s)\|_{-1}^2ds +C\Delta t^2\int_{t^{n}}^{t^{n+1}}\left|\frac{d^2R(s)}{dt^2}\right|^2ds+C\Delta t^2\int_{t^n}^{t^{n+1}}\left|\frac{dR(s)}{dt}\right|^2ds.
	\end{align}
	We sum up the above inequality for the indices
        from $0$ to $n$ and use the discrete Gronwall lemma \ref{sec2_Lemma3}
        to finish the proof.
\end{proof}

\subsubsection{Scheme 1B}\label{gPAV_1B}

An alternative algorithm, in some sense reciprocal to
the scheme presented in the previous section,
is as follows.
Let $\phi^0$, $\mu^0$ and $R^0$ be defined by \eqref{CH1_initial}.
Given ($\phi^n$, $R^n$),
we compute ($\phi^{n+1}$, $\mu^{n+1}$,$R^{n+1}$) by the following procedure,
\begin{subequations}\label{CH1_weak_1B}
	\begin{align}
	\label{CH1_eq1_1B}
	&\frac{\phi^{n+1}-\phi^{n}}{\Delta t} = \Delta\mu^{n+1},\\
	\label{CH1_eq2_1B}
	&\mu^{n+1} =-\Delta\phi^{n+1}+ \lambda\phi^{n+1} + |\xi_{_{1B}}^{n}|^2h(\phi^{n}),\\
	\label{CH1_eq3_1B}
	&\frac{R^{n+1}-R^{n}}{\Delta t} = -\frac{\xi_{_{1B}}^{n+1}}{2\sqrt{\bm{E}[\phi^{n+1}]}}\int_{\Omega}|\nabla\mu^{n+1}|^2d\bm{x},
	\end{align}
\end{subequations}
with the boundary conditions
\begin{equation}\label{boundary_1B}
\nabla\phi^{n+1}\cdot\bm{n}=\nabla\mu^{n+1}\cdot\bm{n}=0,
\quad \text{on}\ \partial\Omega,
\end{equation}
 where 
\begin{equation}\label{xi_0_1B}
\xi_{_{1B}}^{n+1} = \frac{R^{n+1}}{\sqrt{\bm{E}[\phi^{n+1}]}}.
\end{equation}
Note that $\xi_{_{1B}}^{n+1}$ is again an approximation of the constant
$\frac{R(t)}{\sqrt{E(t)}}=1$.

\begin{myRemark}
  In this scheme the equations \eqref{CH1_eq1_1B}--\eqref{CH1_eq2_1B}
  are not coupled with the equations \eqref{CH1_eq3_1B} and \eqref{xi_0_1B},
  because of the explicit treatments of $h(\phi^n)$ and $\xi_{_{1B}}^{n}$
  in \eqref{CH1_eq2_1B}. Therefore, the computations for
  ($\phi^{n+1},\mu^{n+1}$) and for $R^{n+1}$ are de-coupled 
  with this scheme.
  
\end{myRemark}

  Substituting the $\xi^{n+1}_{_{1B}}$ expression in \eqref{xi_0_1B}
  into equation \eqref{CH1_eq3_1B} leads to
	\begin{equation}\label{xi_1B}
	\xi^{n+1}_{_{1B}}=\frac{R^n}{\sqrt{\bm{E}[\phi^{n+1}]}+\frac{\Delta t}{2\sqrt{\bm{E}[\phi^{n+1}]}}\int_{\Omega}|\nabla \mu^{n+1}|^2d\bm{x}}.
	\end{equation}
	Since $R^0>0$, we conclude by induction that $\xi^{n}_{_{1B}}>0$ for all $n$.

Given $\phi^{n}$, $\xi_{_{1B}}^{n}$ and $R^{n}$, we first compute
$\phi^{n+1}$ and $\mu^{n+1}$ by solving  equations
\eqref{CH1_eq1_1B}--\eqref{CH1_eq2_1B}, together with
the boundary conditions \eqref{boundary_1B}.
Then, we compute $\bm{E}[\phi^{n+1}]$ and $\xi_{_{1B}}^{n+1}$
 by equations \eqref{energy_0} and \eqref{xi_1B}, respectively.
$R^{n+1}$ can then be computed based on equation \eqref{xi_0_1B} as follows,
\begin{equation}\label{Rnp1}
R^{n+1} = \xi_{_{1B}}^{n+1}\sqrt{E[\phi^{n+1}]}.
\end{equation}
We therefore conclude that $R^{n+1}>0$ for all $n$ with this scheme.

Similar to Scheme 1A from Section \ref{gPAV}, this scheme requires
the solution of the Cahn-Hilliard field equation only once per
time step. Its operation count per time step is comparable to
that of Scheme 1A, and is approximately a half of those of
the original gPAV scheme~\cite{YangD2020} and
the SAV scheme~\cite{ShenXY2018,YangLD2019}.
Note that in Scheme 1A $R^{n+1}$ is computed first, followed by
the fields $(\phi^{n+1},\mu^{n+1})$. In contrast,
in the current scheme the fields $(\phi^{n+1},\mu^{n+1})$ are
computed first, followed by the variables $(\xi_{_{1B}}^{n+1},R^{n+1})$.

\paragraph{Stability Properties}


The scheme given by equations \eqref{CH1_weak_1B}--\eqref{xi_0_1B}
is unconditionally stable. Its stability properties are summarized by
the following results.

\begin{myLemma}\label{sec3_Lemma1_1B} The scheme \eqref{CH1_weak_1B} is mass conserving
  in the sense that $(\phi^{n+1},1)=(\phi^n,1)$.
\end{myLemma}
\begin{proof} Integrating equation \eqref{CH1_eq1_1B} over $\Omega$ and using
  the boundary condition \eqref{boundary_1B} lead to the result.
\end{proof}

\begin{myLemma}\label{sec3_Lemma2_1B}
  With the scheme \eqref{CH1_weak_1B},
  \begin{align}
    &
    0<R^{n+1} \leq R^{n}\leq M, \label{sec3_eq1_1B} \\
    &
    0<\xi_{_{1B}}^{n+1}
    \leq \frac{M}{\sqrt{C_0}},
    \label{sec3_eq5_1B}
  \end{align}
for some constant $C_0>0$, and a constant $M$
that depends only on $\Omega$, $\phi_{in}$ and $c_0$.
\end{myLemma}

\begin{myTheorem}\label{sec3_Lemma3_1B}
  Suppose $\phi_{in}\in H^3(\Omega)$ and the condition \eqref{Assume_H2} holds. The following inequality holds with the scheme \eqref{CH1_weak_1B},
        \begin{align*}
        \|\nabla\phi^{n+1}\|_0^2+\frac{\lambda}{2}\|\phi^{n+1}\|_0^2
        +\frac{\Delta t}{2}\|\nabla\Delta\phi^{n+1}\|_0^2+\lambda\Delta t\sum_{k=0}^{n}\|\Delta\phi^{k+1}\|_0^2+\frac{\Delta t}{2}\sum_{k=0}^{n}\|\nabla\mu^{k+1}\|_0^2\leq \widehat{C}_1,
        \end{align*}
        where $\widehat{C}_1$ is the constant as given
         in Theorem \ref{sec3_Lemma3}. 
\end{myTheorem}
\noindent The proof of this theorem is provided in the Appendix A.

\begin{myTheorem}\label{sec3_Lemma4_1B}
  Suppose $\phi_{in}\in H^4(\Omega)$,
  and the conditions for Lemmas \ref{sec2_Lemma1} and \ref{sec2_Lemma2} hold. The following inequality holds,
	\begin{align*}
	\|\Delta\phi^{n+1}\|_0^2+\frac{\Delta t}{2}\|\Delta^2\phi^{n+1}\|_0^2+\frac{\Delta t}{2}\sum_{k=0}^{n}\|\Delta^2\phi^{k+1}\|_0^2+\frac{\lambda\Delta t}{2}\sum_{k=0}^{n}\|\nabla\Delta\phi^{k+1}\|_0^2\leq\widehat{C}_2,
	\end{align*}
	where the constant $\widehat{C}_2$
        is given in Theorem \ref{sec3_Lemma4}.
\end{myTheorem}
\begin{proof}
  The proof is similar to that for the SAV scheme in \cite{Shen2018Convergence},
  by using the Lemmas \ref{sec2_Lemma1} and \ref{sec2_Lemma2}.
\end{proof}

\paragraph{Error Estimate}
We define the errors of the variables by \eqref{def_error}.
%
%
Suppose that the solution $(\phi,\mu)$ of equations \eqref{CH}-\eqref{boundary} satisfies \eqref{regularity_1}. 
Based on Lemma \ref{sec2_Lemma2} and Theorems \ref{sec3_Lemma3_1B} and
\ref{sec3_Lemma4_1B}, we have the same results
expressed by the
inequalities \eqref{sec4_eq1}, \eqref{sec4_eq2} and \eqref{sec4_eq3}, i.e.
\begin{align*}
&\|\phi(t^n)\|_{2}\leq C, \quad \|\phi^n\|_{2}\leq C,\\
&|h(\phi)|, |h'(\phi)|, |h''(\phi)|, |h(\phi^n)|, |h'(\phi^n)|, |h''(\phi^n)|\leq C,\\
&\int_{0}^T\left|\frac{d^2R}{dt^2}\right|^2dt\leq C\int_{0}^T\left(\|\phi_{t}\|_{1}^4+\|\phi_{t}\|_{1}^2+\|\phi_{tt}\|_{-1}^2\right)dt.
\end{align*}

The truncation errors $T_{\phi_{1B}}^{n+1}$ and $T_{R_{1B}}^{n+1}$ are given by
\begin{subequations}\label{truncation1B}
	\begin{align}
	\label{truncation1B_eq1}
	&\frac{\phi(t^{n+1})-\phi(t^{n})}{\Delta t} = \Delta\mu(t^{n+1})+\frac{1}{\Delta t}T_{\phi_{1B}}^{n+1},\\
	\label{truncation1B_eq2}
	&\mu(t^{n+1})=-\Delta\phi(t^{n+1})+ \lambda\phi(t^{n+1}) + \frac{R(t^{n+1})^2}{\bm{E}[\phi(t^{n+1})]}h(\phi(t^{n+1})),\\
	\label{truncation1B_eq3}
	&\frac{R(t^{n+1})-R(t^{n})}{\Delta t} = -\frac{R(t^{n+1})}{2\bm{E}[\phi(t^{n+1})]} \int_{\Omega}|\nabla\mu(t^{n+1})|^2d\bm{x}+\frac{1}{\Delta t}T_{R_{1B}}^{n+1},
	\end{align}
\end{subequations}
where
\begin{equation}\label{sec4_eq5_1B}
T_{\phi_{1B}}^{n+1} = \int_{t^n}^{t^{n+1}}(t^n-t)\phi_{tt}(t)dt\quad {\rm and} \quad T_{R_{1B}}^{n+1} = \int_{t^n}^{t^{n+1}}(t^{n}-t)\frac{d^2R(t)}{dt^2}dt.
\end{equation}

\begin{myTheorem}\label{sec4_Lemma4_1B}
   Suppose the condition \eqref{regularity_1} and
  the conditions for Theorems \ref{sec3_Lemma3_1B} and \ref{sec3_Lemma4_1B} hold.
  We have the following result with sufficiently small $\Delta t$, 
  \begin{equation*}
  \frac{1}{2}\|\nabla e_{\phi}^{n+1}\|_0^2+ \frac{\lambda}{2}\|e_{\phi}^{n+1}\|_0^2+|e_{R}^{n+1}|^2+\frac{\Delta t}{2}\sum_{k=0}^n\|\nabla e_{\mu}^{k+1}\|_0^2\leq\widehat{C}_4\Delta t^2,
  \end{equation*}
  where $\widehat{C}_4=C\exp(\Delta t\sum_{k=0}^{n}\frac{r^{k+1}}{1-r^{k+1}\Delta t})\int_{0}^{t^{n+1}}\left(\|\phi_t(s)\|_{1}^4+\|\phi_t(s)\|_{1}^2+\|\phi_{tt}(s)\|_{-1}^2\right)ds$, $r^k=1+\|\nabla\mu^k\|_0^2$ and the constant $C$ depends on $T$, $\phi_{in}$, $\Omega$, $\|\phi\|_{L^{\infty}(0,T;W^{3,\infty}(\Omega))}$ and $\|\mu\|_{L^{\infty}(0,T;H^{1}(\Omega))}$. 
\end{myTheorem}
\noindent The proof of this theorem is provided in the Appendix A.

\subsection{Second-Order Schemes}
\label{BDF2}

We next present two second-order schemes for
solving the reformulated system of equations, both of
which are unconditionally energy stable.
Similar to their first-order counterparts
from Section \ref{gPAV_scheme}, these schemes solve
the Cahn-Hilliard field equation only once per time
step. A prominent feature of these schemes lies
in that the Cahn-Hilliard field equation and the dynamic equation
for the auxiliary variable are discretized in time
by different methods, the former by the backward differentiation
formula (BDF2) and the latter by the Crank-Nicolson scheme (CN2).
This allows the computation of the auxiliary variable, and
ensures the positivity of its computed values, in
a very straightforward way.
We provide stability analyses for both schemes,
as well as the error estimate for the second scheme.
Due to a technical difficulty caused by
its multi-step nature, the error estimate for the first
scheme (Scheme 2A) is not available at this time.

\subsubsection{Scheme 2A}\label{BDF2A_scheme}

Suppose ($\phi^0$, $\mu^0$, $R^0$) is given by \eqref{CH1_initial}.
Define
\begin{equation*}
\left.\phi^{n-1}\right|_{n=0}=\phi^0, \quad
\mu^0 = -\Delta\phi^0 + \lambda\phi^0 + h(\phi^0), \quad
\left.\mu^{n-1}\right|_{n=0}=\mu^0. 
\end{equation*}
Given $\phi^n$, $R^n$, $\phi^{n-1}$ and $\mu^{n-1}$  for $n\geq 0$, we compute
$\phi^{n+1}$, $\mu^{n+1}$ and $R^{n+1}$ as follows,
\begin{subequations}\label{CH2}
	\begin{align}
	\label{CH2_eq1}
	&\frac{3\phi^{n+1}-4\phi^{n}+\phi^{n-1}}{2\Delta t} = \Delta\mu^{n+1},\\
	\label{CH2_eq2}
	&\mu^{n+1} =-\Delta\phi^{n+1}+ \lambda\phi^{n+1} + |\xi^{n+1}_{_{2A}}|^2h(\overline{\phi}^{n}),\\
	\label{CH2_eq3}
	&\frac{R^{n+1}-R^{n}}{\Delta t} =
        -\frac{\xi^{n+1}_{_{2A}}}{2\sqrt{\bm{E}[\widetilde{\phi}^{n+1/2}]}}\int_{\Omega}\left|\nabla\widetilde{\mu}^{n+1/2}\right|^2d\bm{x}, \\
        &
        \bm{n}\cdot\nabla\phi^{n+1} = \bm{n}\cdot\nabla\mu^{n+1} = 0,
        \quad \text{on} \ \partial\Omega,
        \label{CH2_eq4}
	\end{align}
\end{subequations}
where 
\begin{equation}\label{xi_2}
\xi^{n+1}_{_{2A}} = \frac{R^{n+1}}{\sqrt{\bm{E}[\overline{\phi}^{n}]}}.
\end{equation}
The symbols in the above equations are defined by
\[
\overline{\phi}^{n}=2\phi^{n}-\phi^{n-1}, \quad \quad
\widetilde{\phi}^{n+1/2}=\frac{3}{2}\phi^{n}-\frac{1}{2}\phi^{n-1}, \quad\quad
\widetilde{\mu}^{n+1/2} = \frac32\mu^n - \frac12\mu^{n-1}.
\]
Obviously, $\overline{\phi}^n$ is a second-order explicit
approximation of $\phi^{n+1}$, and
$\widetilde{\mu}^{n+1/2}$ is a second-order explicit
approximation of $\mu$ at step $(n+1/2)$,
both by extrapolations.
It follows that $\xi^{n+1}_{_{2A}}$ in \eqref{xi_2} is a second-order
approximation of the constant $\frac{R(t)}{\sqrt{E(t)}}=1$
at step $(n+1)$.

Notice that we have used BDF2 to approximate $\frac{\partial\phi}{\partial t}$
in \eqref{CH2_eq1} and enforced this equation at step $(n+1)$.
On the other hand, we approximate $\frac{dR}{dt}$
by the Crank-Nicolson form, and enforce all terms in
the equation \eqref{CH2_eq3}, except for the variable $\xi^{n+1}_{_{2A}}$,
at the time step $(n+1/2)$.
Note that the $\xi^{n+1}_{_{2A}}$ variable in \eqref{CH2_eq3} is approximated
at the time step $(n+1)$ according to \eqref{xi_2}. This is
a crucial point, and it allows $R^{n+1}$ to be computed
by a well-defined formula and ensures the positivity of its
values. It should be appreciated that this approximation of $\xi^{n+1}_{_{2A}}$
does not spoil the second-order accuracy of the overall scheme,
because $\xi^{n+1}_{_{2A}}$ is an approximation of the constant $\frac{R(t)}{\sqrt{E(t)}}=1$.

Thanks to the explicit nature of $\widetilde{\mu}^{n+1/2}$ in \eqref{CH2_eq3},
the computation for $R^{n+1}$ from \eqref{CH2_eq3} and \eqref{xi_2}
is un-coupled with the computations for $\phi^{n+1}$ and $\mu^{n+1}$
from \eqref{CH2_eq1}--\eqref{CH2_eq2}.
%
Substituting the $\xi^{n+1}_{_{2A}}$ expression in \eqref{xi_2}
into equation \eqref{CH2_eq3}, one finds
\begin{equation}\label{xi_3}
  \xi^{n+1}_{_{2A}}=\frac{R^n}{\sqrt{\bm{E}[\overline{\phi}^{n}]}+\frac{\Delta t}{2\sqrt{\bm{E}[\widetilde{\phi}^{n+1/2}]}}\int_{\Omega}|\nabla\widetilde{\mu}^{n+1/2}|^2d\bm{x}}.
\end{equation}
Since $R^0>0$, we conclude by induction that $\xi^{n}_{_{2A}}>0$ for all $n\geq 0$.
It then follows that $R^{n}>0$ for all $n\geq 0$.
%

To implement this scheme,
we first compute $\xi^{n+1}_{_{2A}}$ from equation \eqref{xi_3} and
$R^{n+1}$ from \eqref{xi_2}.
Then we solve equations \eqref{CH2_eq1}--\eqref{CH2_eq2}
for $\phi^{n+1}$ and $\mu^{n+1}$.
It can be noted that the Cahn-Hilliard field equation
\eqref{CH2_eq1}--\eqref{CH2_eq2} is solved only once per time step
with this scheme.

\paragraph{Stability Properties}


The scheme given by equations \eqref{CH2_eq1}--\eqref{xi_2}
is unconditionally stable. The following lemmas and theorems
summarize its stability properties.

\begin{myLemma}\label{sec5_Lemma0} The scheme \eqref{CH2} is mass conserving
  in the sense that $(\phi^{n+1},1)=(\phi^0,1)$.
\end{myLemma}
\begin{proof}
  Taking the $L^2$ inner product between equation \eqref{CH2_eq1} and
  the constant $1$ leads to
  $(\phi^{n+1},1) = \frac43(\phi^n,1) - \frac13(\phi^{n-1},1)$ for all $n\geq 0$.
  Since $\phi^{n-1}|_{n=0}=\phi^0$ by definition,
  we conclude by induction that $(\phi^{n},1)=(\phi^0,1)$
  for all $n\geq 0$.
\end{proof}

\begin{myLemma}\label{sec5_Lemma1}
  With the scheme \eqref{CH2},
  \begin{align}
    &
	  0<R^{n+1} \leq R^{n}\leq M, \label{sec5_eq1} \\
          &
          0< \xi^{n+1}_{_{2A}}\leq \frac{M}{\sqrt{C_0}}, 
	\end{align}
        for some constant $C_0>0$, and a constant $M$ that
        depends on $\Omega$, $\phi_{in}$ and $c_0$.
\end{myLemma}


\begin{myTheorem}\label{sec5_Lemma2}
  Suppose $\phi_{in}\in H^3(\Omega)$ and the condition \eqref{Assume_H2} holds. The following inequality
  holds for all $n$ with the scheme \eqref{CH2},
	\begin{align*}
	&\|\nabla\phi^{n+1}\|_0^2+\|\nabla(2\phi^{n+1}-\phi^{n})\|_0^2+\frac{1}{2}\|\nabla(\phi^{n+1}-\phi^{n})\|_0^2+\frac{\lambda}{2}(\|\phi^{n+1}\|_0^2+\|2\phi^{n+1}-\phi^{n}\|_0^2)
	\\
	&\quad +\lambda\Delta t\|\Delta\phi^{n+1}\|_0^2+\frac{3\Delta t}{2}\|\nabla\Delta\phi^{n+1}\|_0^2+\frac{\Delta t}{2}\|\nabla\Delta(\phi^{n+1}-\phi^n)\|_0^2+\Delta t\sum_{k=0}^{n}\|\nabla\mu^{k+1}\|_0^2\leq\widehat{C}_5,
	\end{align*}
	where $\widehat{C}_5=(2\|\nabla\phi^{0}\|_0^2+\lambda\|\phi^{0}\|_0^2+\lambda\Delta t\|\Delta\phi^{0}\|_0^2+\frac{3\Delta t}{2}\|\nabla\Delta\phi^{0}\|_0^2)\exp\left(C(M)T\right)$.
\end{myTheorem}
\noindent The proof of this theorem is provided in the Appendix A.

\begin{myTheorem}\label{sec5_Lemma3}
  Suppose $\phi_{in}\in H^4(\Omega)$, and that
  the conditions for Lemmas \ref{sec2_Lemma1} and \ref{sec2_Lemma2} hold. The following inequality holds for all $n$ with the scheme \eqref{CH2}, 
	\begin{align*}
	&\frac{1}{2}\|\Delta\phi^{n+1}\|_0^2+\frac{1}{2}\|\Delta(2\phi^{n+1}-\phi^{n})\|_0^2+\frac{1}{2}\|\Delta(\phi^{n+1}-\phi^{n})\|_0^2+\lambda\Delta t\|\nabla\Delta\phi^{n+1}\|_0^2\\
	&\qquad +\frac{3\Delta t}{2}\|\Delta^2\phi^{n+1}\|_0^2+\frac{\Delta t}{2}\|\Delta^2(\phi^{n+1}-\phi^{n})\|_0^2+\lambda\Delta t\sum_{k=0}^{n}\|\nabla\Delta\phi^{k+1}\|_0^2\leq\widehat{C}_6,
	\end{align*}
	where $\widehat{C}_6=\|\Delta\phi^{0}\|_0^2+\frac{3\Delta t}{2}\|\Delta^2\phi^{0}\|_0^2+\lambda\Delta t\|\nabla\Delta\phi^0\|_0^2+C(M)T$.
\end{myTheorem}
\noindent The proof of this theorem is provided in the Appendix A.

\subsubsection{Scheme 2B}
\label{BDF2_scheme}

Suppose ($\phi^0$, $\mu^0$, $R^0$) are given by \eqref{CH1_initial}, and
let $\phi^{n-1}|_{n=0}=\phi^0$ and $R^{n-1}|_{n=0}=R^0$.
Given ($\phi^n$, $\phi^{n-1}$, $R^n$, $R^{n-1}$) for $n\geq 0$, we compute
($\phi^{n+1}$, $\mu^{n+1}$, $R^{n+1}$) as follows,
\begin{subequations}\label{CH2_2B}
	\begin{align}
	\label{CH2_eq1_2B}
	&\frac{3\phi^{n+1}-4\phi^{n}+\phi^{n-1}}{2\Delta t} = \Delta\mu^{n+1},\\
	\label{CH2_eq2_2B}
	&\mu^{n+1} =-\Delta\phi^{n+1}+ \lambda\phi^{n+1} + |\widehat{\xi}^{n}_{_{2B}}|^2h(\overline{\phi}^{n}),\\
	\label{CH2_eq3_2B}
	&\frac{R^{n+1}-R^{n}}{\Delta t} = -\frac{\xi^{n+1}_{_{2B}}}{2\sqrt{E[\widetilde{\phi}^{n+1/2}]}}\int_{\Omega}\left|\nabla\mu^{n+1/2}\right|^2d\bm{x},\\
        &
        \bm{n}\cdot\nabla\phi^{n+1} = \bm{n}\cdot\nabla\mu^{n+1} = 0,
        \quad \text{on} \ \partial\Omega,
        \label{CH2_eq4_2B}
	\end{align}
\end{subequations}
where 
\begin{equation}\label{xi_2_2B}
\xi^{n+1}_{_{2B}} = \frac{R^{n+1}}{\sqrt{E[\phi^{n+1}]}}.
\end{equation}
The symbols in the above equations are defined by
\begin{align*}
  &
  \overline{R}^n = 2R^n - R^{n-1}, \quad
  \overline{\phi}^n = 2\phi^n - \phi^{n-1}, \quad  
  \quad \widehat{\xi}^{n}_{_{2B}} = \frac{\overline{R}^{n}}{\sqrt{E[\overline{\phi}^{n}]}}, \\
  &
  \widetilde{\phi}^{n+1/2} = \frac32\phi^n - \frac12\phi^{n-1}, \quad
  \mu^{n+1/2} = \frac12(\mu^{n+1} + \mu^n).
\end{align*}
It can be noted that $\overline{R}^n$ and $\overline{\phi}^n$
are second-order explicit approximations of $R^{n+1}$ and $\phi^{n+1}$, respectively.
So $\widehat{\xi}^{n}_{_{2B}}$ is a second-order explicit approximation
of the constant $\frac{R(t)}{\sqrt{E(t)}}=1$ at the time step $(n+1)$.

It should be noted that equations \eqref{CH2_eq1_2B}--\eqref{CH2_eq2_2B}
are enforced at the time step $(n+1)$, while equation \eqref{CH2_eq3_2B} is
enforced at the step $(n+1/2)$, except for the term $\xi^{n+1}_{_{2B}}$,
which is approximated at the time step $(n+1)$.
Similar to Scheme 2A in the previous subsection, this treatment
of $\xi^{n+1}_{_{2B}}$ here does not spoil the second-order accuracy of the scheme
and ensures the positivity of the computed values for $R^{n+1}$
and $\xi_{2B}^{n+1}$.

Substitution of the $\xi_{_{2B}}^{n+1}$ expression in \eqref{xi_2_2B}
into equation \eqref{CH2_eq3_2B} leads to
	\begin{equation}\label{xi_3_2B}
	\xi^{n+1}_{_{2B}}=\frac{R^n}{\sqrt{\bm{E}[\phi^{n+1}]}+\frac{\Delta t}{2\sqrt{\bm{E}[\widetilde{\phi}^{n+1/2}]}}\int_{\Omega}|\nabla\mu^{n+1/2}|^2d\bm{x}}.
	\end{equation}
	Since $R^0>0$,
        we conclude by induction that $\xi_{_{2B}}^{n+1}>0$ and $R^{n+1}>0$
        for all $n\geq 0$ with the current scheme.

Thanks to the explicit nature of $\widehat{\xi}^n_{_{2B}}$, the field equations
\eqref{CH2_eq1_2B}--\eqref{CH2_eq2_2B} are de-coupled from the equation
\eqref{CH2_eq3_2B}. To compute $(\phi^{n+1},\mu^{n+1},R^{n+1})$, we can first solve
\eqref{CH2_eq1_2B}--\eqref{CH2_eq2_2B} and \eqref{CH2_eq4_2B}
for $\phi^{n+1}$ and $\mu^{n+1}$. Then we compute $\xi^{n+1}_{_{2B}}$ by \eqref{xi_3_2B},
and compute $R^{n+1}$ by equation \eqref{xi_2_2B}.


\paragraph{Stability Properties}


This scheme is also unconditionally energy stable. Its
stability properties are summarized by the following results.

\begin{myLemma}\label{sec5_Lemma1_2B} The scheme \eqref{CH2_2B} is mass conserving
  in the sense that $(\phi^{n+1},1) = (\phi^0,1)$.
\end{myLemma}

\begin{myLemma}\label{sec5_Lemma2_2B} The scheme \eqref{CH2_2B} satisfies,
  for all $n$,
  \begin{align}
    &
    0<R^{n+1} \leq R^{n}\leq M, \label{sec5_eq1_2B} \\
    &
    0< \xi^{n+1}_{_{2B}}\leq \frac{M}{\sqrt{C_0}}, \\
    &
    \left|\widehat{\xi}^{n}_{_{2B}}\right|\leq \frac{3M}{\sqrt{C_0}},
  \end{align}
  for some constant $C_0>0$, and a constant $M$ that depends on
  $\Omega$, $\phi_{in}$ and $c_0$.
\end{myLemma}


\begin{myTheorem}\label{sec5_Lemma3_2B}  Suppose $\phi_{in}\in H^3(\Omega)$ and the condition \eqref{Assume_H2} holds.
  The following inequality holds for all n with the scheme \eqref{CH2_2B},
	\begin{align*}
	&\|\nabla\phi^{n+1}\|_0^2+\|\nabla(2\phi^{n+1}-\phi^{n})\|_0^2+\frac{1}{2}\|\nabla(\phi^{n+1}-\phi^{n})\|_0^2+\frac{\lambda}{2}(\|\phi^{n+1}\|_0^2+\|2\phi^{n+1}-\phi^{n}\|_0^2)\\
	&\quad +\lambda\Delta t\|\Delta\phi^{n+1}\|_0^2+\frac{3\Delta t}{2}\|\nabla\Delta\phi^{n+1}\|_0^2+\frac{\Delta t}{2}\|\nabla\Delta(\phi^{n+1}-\phi^{n})\|_0^2+\Delta t\sum_{k=0}^{n}\|\nabla\mu^{k+1}\|_0^2\leq\widehat{C}_5,
	\end{align*}
	where $\widehat{C}_5$
       is given in Theorem \ref{sec5_Lemma2}.
\end{myTheorem}
\noindent The proof of this theorem is provided in the Appendix A.

\begin{myTheorem}\label{sec5_Lemma4_2B} Suppose $\phi^0\in H^4(\Omega)$,
  and that the conditions for
  Lemmas \ref{sec2_Lemma1} and \ref{sec2_Lemma2} hold.
  Then the following inequality holds with the scheme \eqref{CH2_2B},
	\begin{align*}
	&\frac{1}{2}\|\Delta\phi^{n+1}\|_0^2+\frac{1}{2}\|\Delta(2\phi^{n+1}-\phi^{n})\|_0^2+\frac{1}{2}\|\Delta(\phi^{n+1}-\phi^{n})\|_0^2+\lambda\Delta t\|\nabla\Delta\phi^{n+1}\|_0^2\\
	&\qquad +\frac{3\Delta t}{2}\|\Delta^2\phi^{n+1}\|_0^2+\frac{\Delta t}{2}\|\Delta^2(\phi^{n+1}-\phi^{n})\|_0^2+\lambda\Delta t\sum_{k=0}^{n}\|\nabla\Delta\phi^{k+1}\|_0^2\leq\widehat{C}_6,
	\end{align*}
	where $\widehat{C}_6$
        is given in Theorem \ref{sec5_Lemma3}.
\end{myTheorem}
\noindent The proof of this theorem is provided in the Appendix A.

\paragraph{Error Estimate}
Assume that 
\begin{align}\label{regularity_2}
&\phi \in L^{\infty}(0,T;W^{3,\infty}(\Omega)), \quad  \phi_t \in L^{\infty}(0,T;L^{2}(\Omega))\cap L^{4}(0,T;H^{1}(\Omega))\cap L^{2}(0,T;H^{1}(\Omega)),\nonumber\\
&\phi_{tt} \in L^{2}(0,T;H^{1}(\Omega)), \quad \phi_{ttt} \in L^{\infty}(0,T;H^{-1}(\Omega)), \quad \mu \in L^{\infty}(0,T;H^{1}(\Omega)).
\end{align}
By Lemma \ref{sec2_Lemma2} and Theorems \ref{sec5_Lemma3_2B}
and \ref{sec5_Lemma4_2B}, we can also
arrive at the boundedness properties in
\eqref{sec4_eq1} and \eqref{sec4_eq2}.
In light of the relation $R(t)=\sqrt{E[\phi]}$, we have
\begin{equation}\label{sec5_eq10_2B}
\frac{d^3R}{dt^3} =\frac{3}{8\sqrt{\bm{E}[\phi]^{5}}}\left(\frac{dE}{dt}\right)^3-\frac{3}{4\sqrt{\bm{E}[\phi]^{3}}}\frac{dE}{dt}\frac{d^2E}{dt^2}+\frac{1}{2\sqrt{\bm{E}[\phi]}}\frac{d^3E}{dt^3},
\end{equation}
where
\begin{align*}
\frac{d^2E}{dt^2}&=\int_{\Omega}\left(|\nabla\phi_t|^2
+\nabla\phi\cdot\nabla\phi_{tt}+\lambda|\phi_{t}|^2+\lambda\phi\phi_{tt}+h'(\phi)|\phi_{t}|^2+h(\phi)\phi_{tt}\right)d\bm{x},\\
\frac{d^3E}{dt^3}&=\int_{\Omega}\left(3\nabla \phi_t\cdot\nabla\phi_{tt}+\nabla \phi\cdot\nabla\phi_{ttt}+3\lambda\phi_t\phi_{tt}+\lambda\phi\phi_{ttt}+h''(\phi)\phi_t^3
+3h'(\phi)\phi_t\phi_{tt}+h(\phi)\phi_{ttt}\right)d\bm{x}.
\end{align*}
It follows that
\begin{equation}\label{sec5_eq11}
\int_{0}^T\left|\frac{d^3R}{dt^3}\right|^2dt\leq C\int_{0}^T\left(\|\phi_{t}\|_{1}^4+\|\phi_{t}\|_{1}^2+\|\phi_{tt}\|_{1}^2+\|\phi_{ttt}\|_{-1}^2\right)dt.
\end{equation}

Based on the Taylor expansion theorem, we arrive at
\begin{subequations}\label{truncation2_2B}
	\begin{align}
	\label{truncation2_eq1_2B}
	&\frac{3\phi(t^{n+1})-4\phi(t^{n})+\phi(t^{n-1})}{2\Delta t} = \Delta\mu(t^{n+1})+\frac{1}{\Delta t}T_{\phi_{2B}}^{n+1},\\
	\label{truncation2_eq2_2B}
	&\mu(t^{n+1})=-\Delta\phi(t^{n+1})+ \lambda\phi(t^{n+1}) + \frac{R(t^{n+1})^2}{\bm{E}[\phi(t^{n+1})]}h(\phi(t^{n+1})),\\
	\label{truncation2_eq3_2B}
	&\frac{R(t^{n+1})-R(t^{n})}{\Delta t} = -\frac{R(t^{n+1})}{\sqrt{\bm{E}[\phi(t^{n+1})]}}\frac{1}{2\sqrt{\bm{E}[\phi(t^{n+1/2})]}} \int_{\Omega}\left|\nabla\mu(t^{n+1/2})\right|^2d\bm{x}+\frac{1}{\Delta t}T_{R_{2B}}^{n+1},
	\end{align}
\end{subequations}
where 
\begin{equation}\label{sec5_eq13_2B}
  \left\{
  \begin{split}
   &T_{\phi_{2B}}^{n+1}  =\int_{t^{n}}^{t^{n+1}}(t-t^{n})^2\phi_{ttt}(t)dt-\frac{1}{4}\int_{t^{n-1}}^{t^{n+1}}(t-t^{n-1})^2\phi_{ttt}(t)dt,\\
   & T_{R_{2B}}^{n+1} = \frac{1}{2}\int_{t^{n+1/2}}^{t^{n+1}}(t^{n+1}-t)^2\frac{d^3R}{dt^3}(t)dt-\frac{1}{2}\int_{t^{n}}^{t^{n+1/2}}(t^n-t)^2\frac{d^3R}{dt^3}(t)dt.
  \end{split}
  \right.
\end{equation}

\begin{myTheorem}\label{sec5_Lemma5_2B}
  Suppose the condition \eqref{regularity_2}, 
  and the conditions for Theorems
  \ref{sec5_Lemma3_2B} and \ref{sec5_Lemma4_2B} hold.
  The following inequality holds for sufficiently small $\Delta t$,
	\begin{equation*}
	  \frac{1}{2}\left(\|\nabla e_{\phi}^{n+1}\|_0^2+\|\nabla(2e_{\phi}^{n+1}-e_{\phi}^{n})\|_0^2\right)+\frac{\lambda}{2}\left(\|e_{\phi}^{n+1}\|_0^2+\|2e_{\phi}^{n+1}-e_{\phi}^{n}\|_0^2\right)+\frac{\Delta t}{2}\|\nabla e_{\mu}^{n+1}\|_0^2
          +\left|e_{R}^{n+1}\right|^2\leq\widehat{C}_7\Delta t^4,
	\end{equation*}
	where $\widehat{C}_7=C\exp(\Delta t\sum_{k=0}^{n+1}\frac{r^{k+1/2}}{1-r^{k+1/2}\Delta t})\int_{0}^{t^{n+1}}\left(\|\phi_{t}(s)\|_{1}^4+\|\phi_{t}(s)\|_{1}^2+\|\phi_{tt}(s)\|_{1}^2+\|\phi_{ttt}(s)\|_{-1}^2\right)ds$,
        $r^{k+1/2}=1+\|\nabla\mu^{k+1/2}\|_0^2$,
        and the constant $C$  depends on
        $T$, $\phi_{in}$, $\Omega$, $\|\phi\|_{L^{\infty}(0,T;W^{3,\infty}(\Omega))}$, $\|\phi_t\|_{L^{\infty}(0,T;L^{2}(\Omega))}$ and $\|\mu\|_{L^{\infty}(0,T;H^{1}(\Omega))}$. 
\end{myTheorem}
\noindent The proof of this theorem is provided
in the Appendix A.

%

\vspace{0.1in}

\begin{myRemark}
  \label{rem:rem_final}

  The four schemes presented in this section share several
  common characteristics:
  (i) They are all unconditionally energy stable.
  (ii) The Cahn-Hilliard field equation only needs to be
  computed once per time step, by solving linear algebraic systems
  with a constant coefficient matrix.
  (iii) The auxiliary variable is given by a well-defined
  explicit form, and its computed values are guaranteed
  to be positive.

\end{myRemark}

\begin{myRemark}
  \label{rem:rem_A}
  In the analysis of these schemes (Schemes 1A/1B and 2A/2B)
  we have focused on the
  boundary conditions given by \eqref{boundary}.
  We would like to point out that the stability properties
  about these schemes proved in this section equally hold
  with periodic boundary conditions for the domain.
  
\end{myRemark}


\section{Numerical Examples}
\label{sec:tests}

In this section we provide numerical results to verify
the stability and error analysis of the proposed numerical schemes
from the previous section.
The convergence rates of these schemes are first demonstrated
using a manufactured analytic solution. We then look into
the coalescence of an array of drops and show
that the proposed schemes produce stable and accurate numerical
results.

 In the forthcoming numerical experiments,
 we add the mobility coefficient and the interfacial thickness parameter
 to the Cahn-Hilliard equation as follows so that it resembles the applications
 (e.g.~in two-phase flows) more closely, 
\begin{equation}\label{eq:CHEstd}
\phi_t=m_0\Delta \mu+f(\bs x,t),\quad \mu=-\beta \Delta \phi +h(\phi),
\end{equation}
where $h(\phi)=H'(\phi)$ with $H(\phi)=\frac{\beta}{4\eta^2}(\phi^2-1)^2$.
Here, the constant $m_0$ ($m_0>0$)
is the mobility of the interface, $\eta$ is a characteristic
scale of the interfacial thickness,
$\beta$ is the mixing energy density coefficient and is related to the surface tension by $\beta=\frac{3}{2\sqrt{2}}\sigma \eta$, where the constant $\sigma$ is the interface surface tension.
$f(\bs x,t)$ is a prescribed source term for testing
the convergence rate only,
and will be set to $f=0$ in practical simulations.
For simplicity and efficiency,
we will consider periodic boundary conditions 
in the following tests.
These algorithms are employed to numerically integrate the governing equation \eqref{eq:CHEstd} in time from $t=t_0$ to $t=t_f$ ($t_0$ and $t_f$ to be specified below).


\subsection{Convergence Rates}

\begin{figure}[tbp]
\begin{center}
  \subfigure[$1$st-order scheme vs errors ]{ \includegraphics[scale=.35]{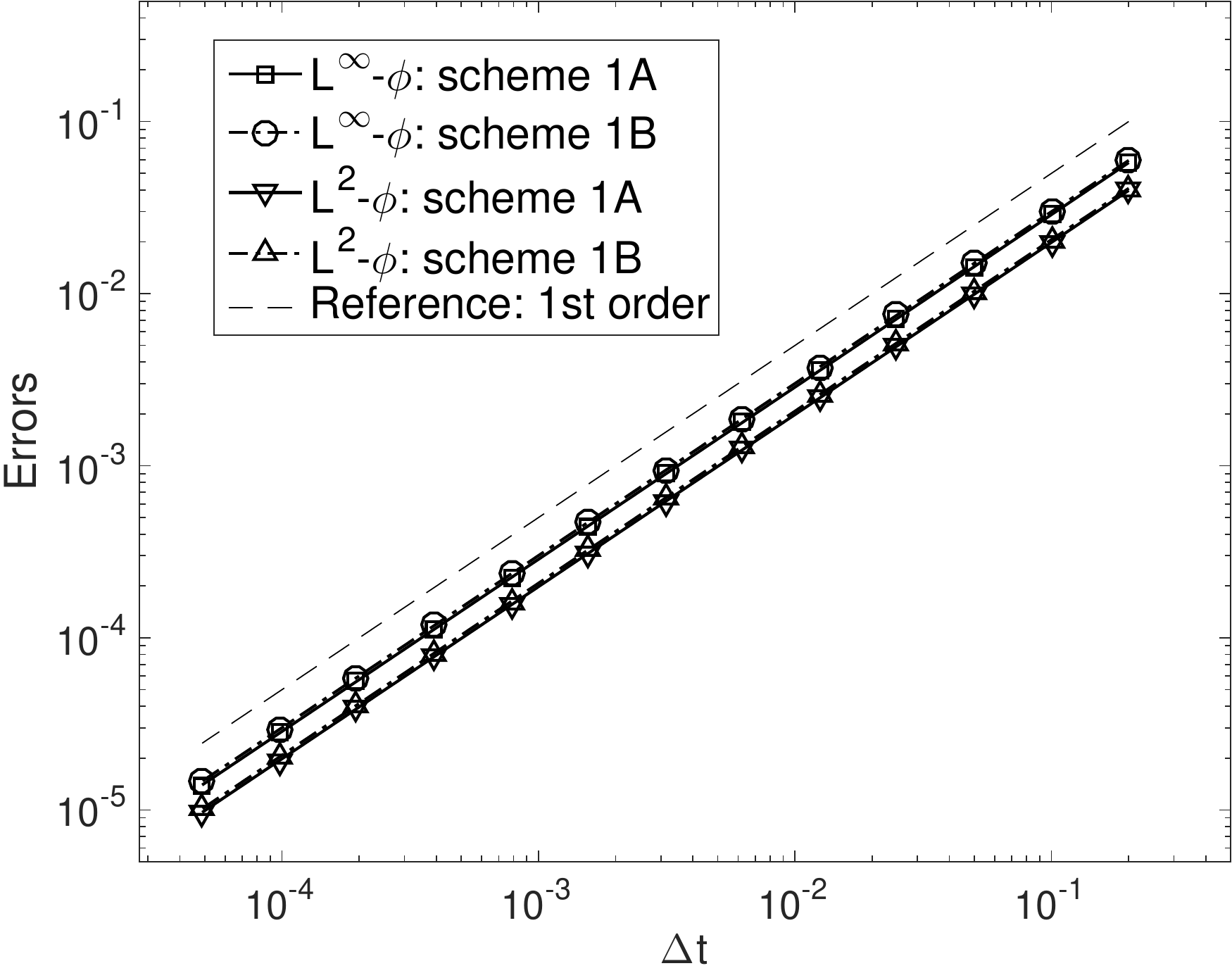}}\qquad
  \subfigure[2nd-order scheme vs errors]{ \includegraphics[scale=.35]{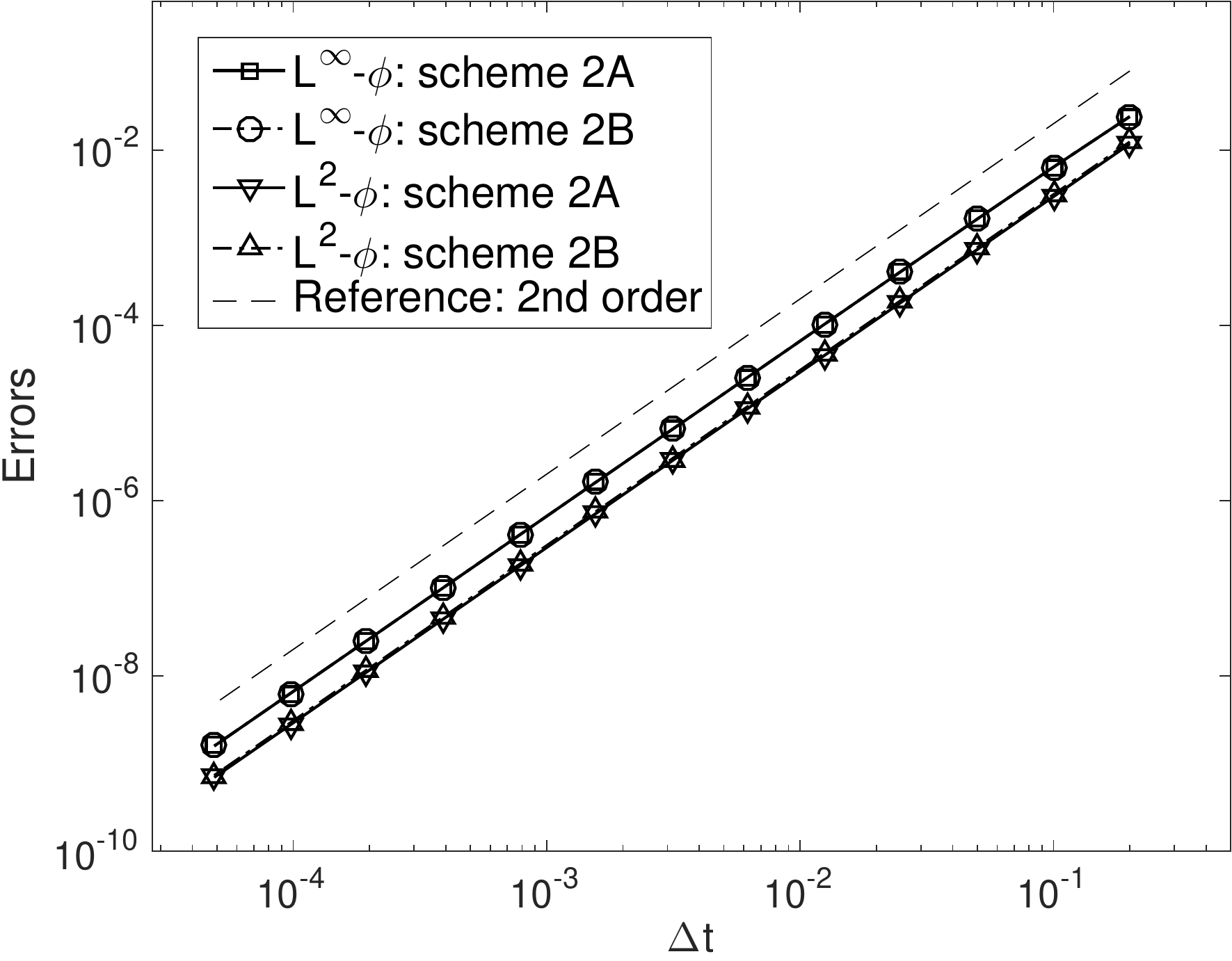}}\\
  \caption{\small Temporal convergence rates in
    $L^{\infty}$ and $L^2$ norms  for the proposed schemes: (a) first-order schemes;
    (b) second-order schemes,
    with $t_0=0.1,$ $t_f=1.1$, $(N_x, N_y)=(20,20)$ and $c_0=1$.
  }
   \label{fig: CHtemporaltest}
\end{center}
\end{figure}

We first test the convergence rates of the proposed
methods using a manufactured analytic solution. Consider equation \eqref{eq:CHEstd} in the domain $\Omega=[0,2]\times [0,2]$ with a manufactured solution
\begin{equation}\label{eq:consoln}
\phi(\bs x,t)=\cos(\pi x)\cos(\pi y)\sin(t).
\end{equation}
The external source term $f(\bs x,t)$ in \eqref{eq:CHEstd} is chosen such that this equation is satisfied by the analytic expression
given in \eqref{eq:consoln}. Periodic conditions are assumed
for the boundaries in the $x$ and $y$ directions.
We employ the
Fourier spectral method for spatial discretization
throughout this section. Let $N_x$ and $N_y$ denote the number
of Fourier collocation points along $x$ and $y$ directions,
respectively. In the simulations, we set  $(N_x,N_y)=(20,20)$, with which the spatial discretization error is negligible compared with the  temporal
discretization error. Other parameters are $t_0=0.1$, $t_f=1.1$, $m_0=0.01$, $\beta=0.01$, $\eta=0.1$ and $c_0=1.$ The $L^{\infty}$ and $L^2$ errors at $t=1.1$ are
plotted respectively for the Schemes 1A/1B and 2A/2B
in Figure \ref{fig: CHtemporaltest}.
We can observe the expected convergence rate for all  cases.
The error curves corresponding to Schemes 1A and 1B, and also for Schemes 2A and 2B,
essentially overlap with each other,
indicating a negligible difference in the convergence rates.

\subsection{Coalescence of an Array of Drops}

\begin{figure}[tbp]
\begin{center}
 \subfigure[$t=0$ ]{ \includegraphics[scale=.23]{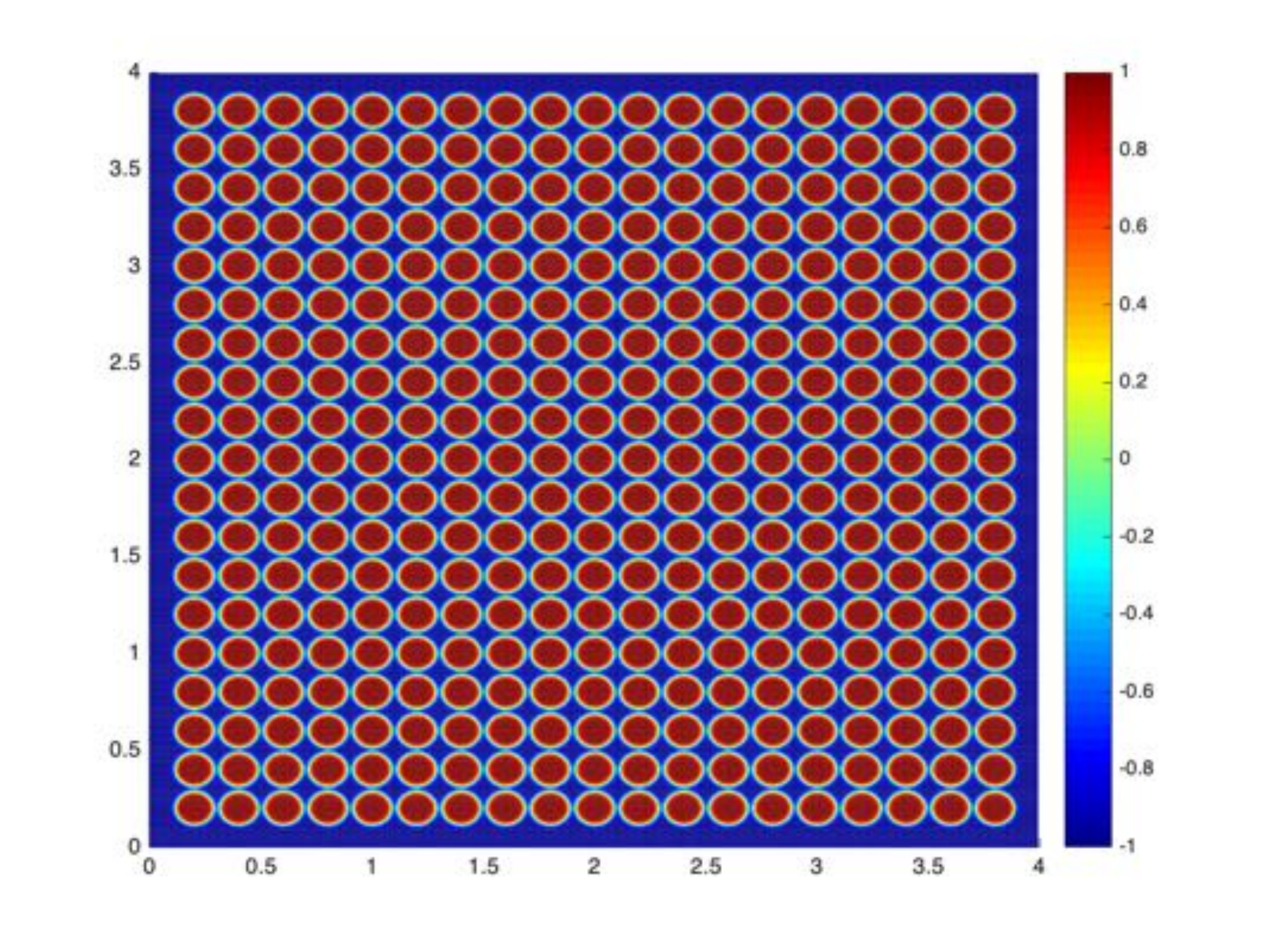}}
  \subfigure[$t=10$ ]{ \includegraphics[scale=.23]{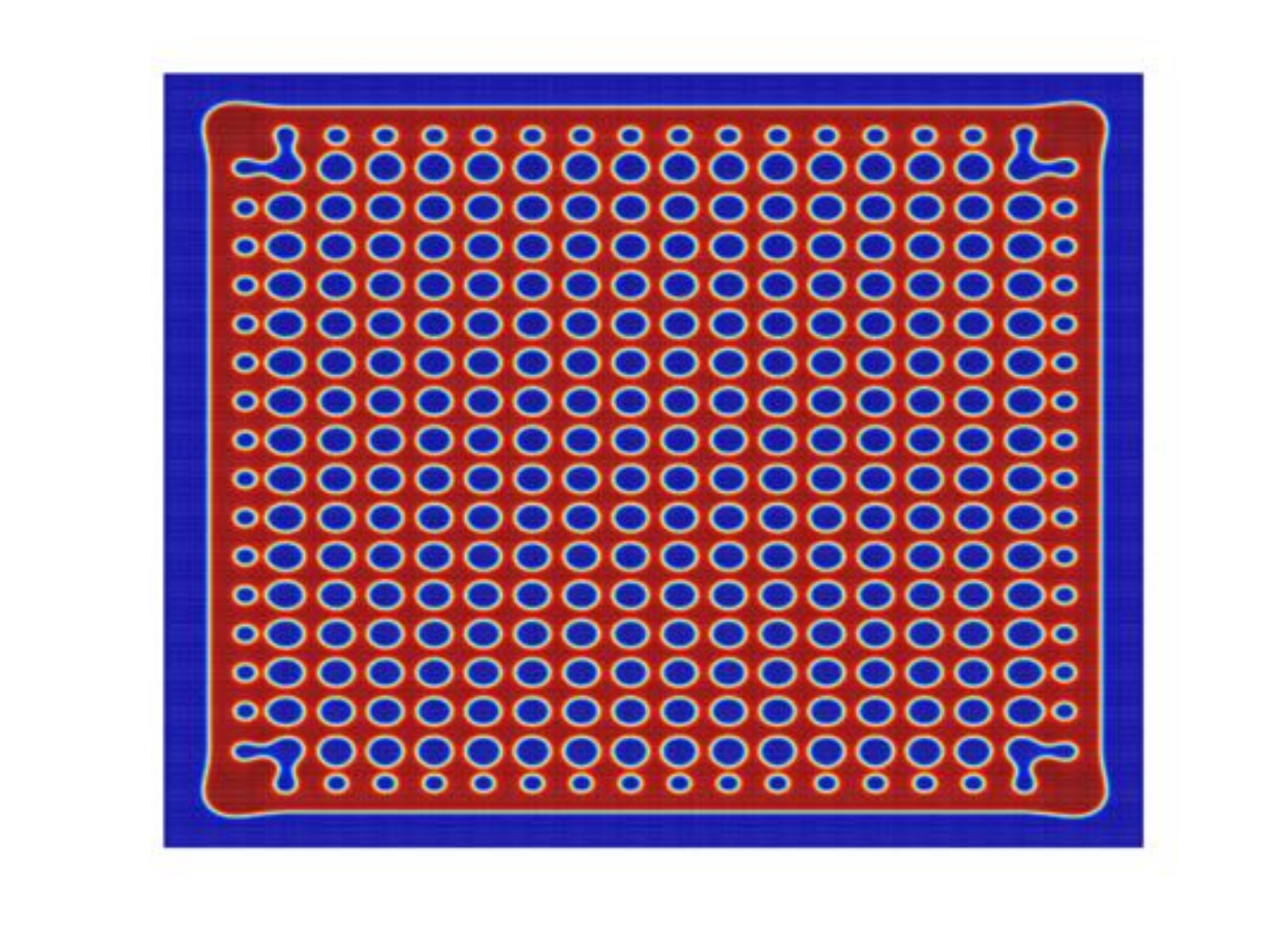}}
  \subfigure[$t=20$ ]{ \includegraphics[scale=.23]{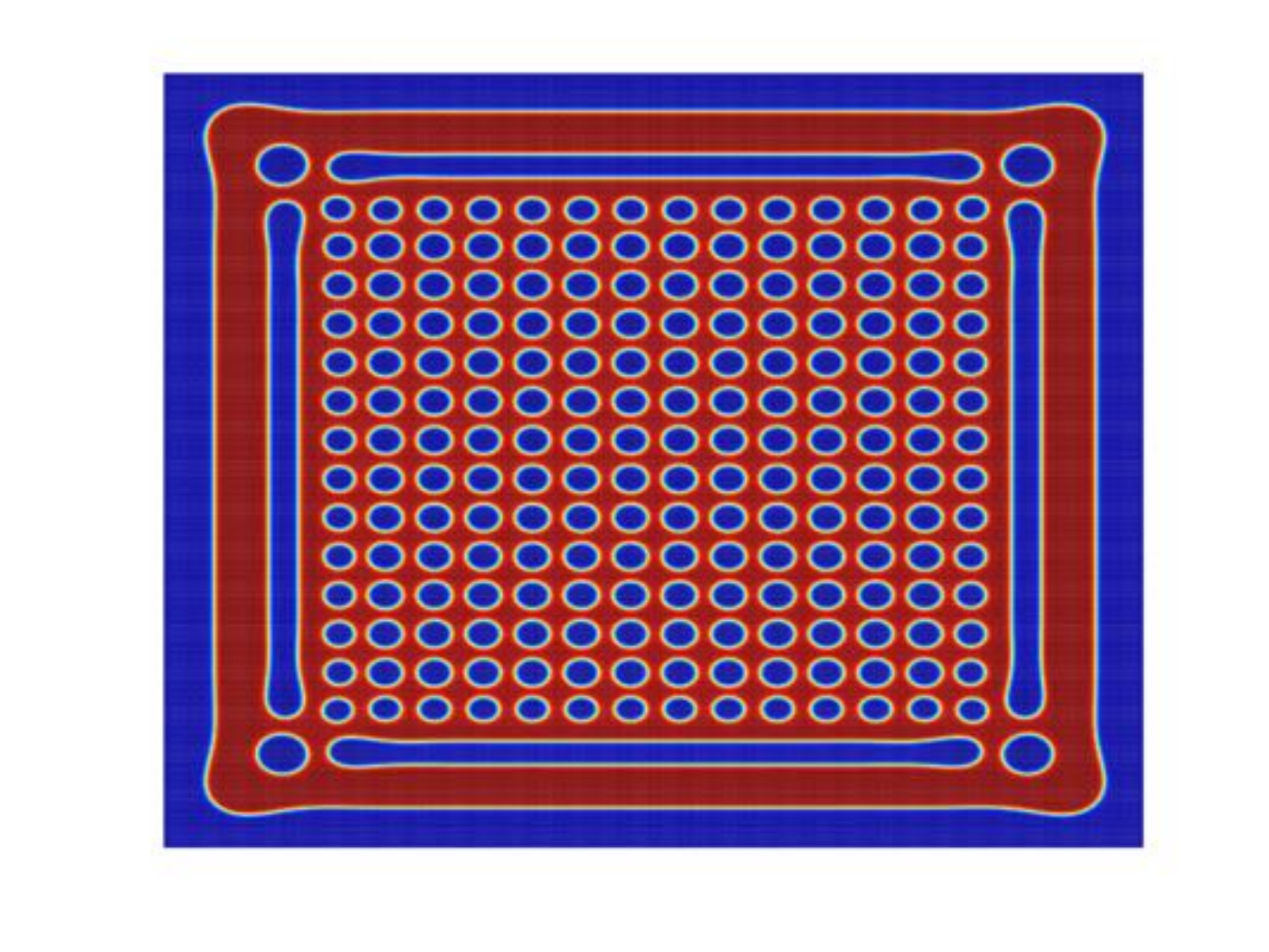}}\\
  \subfigure[$t=30$]{ \includegraphics[scale=.23]{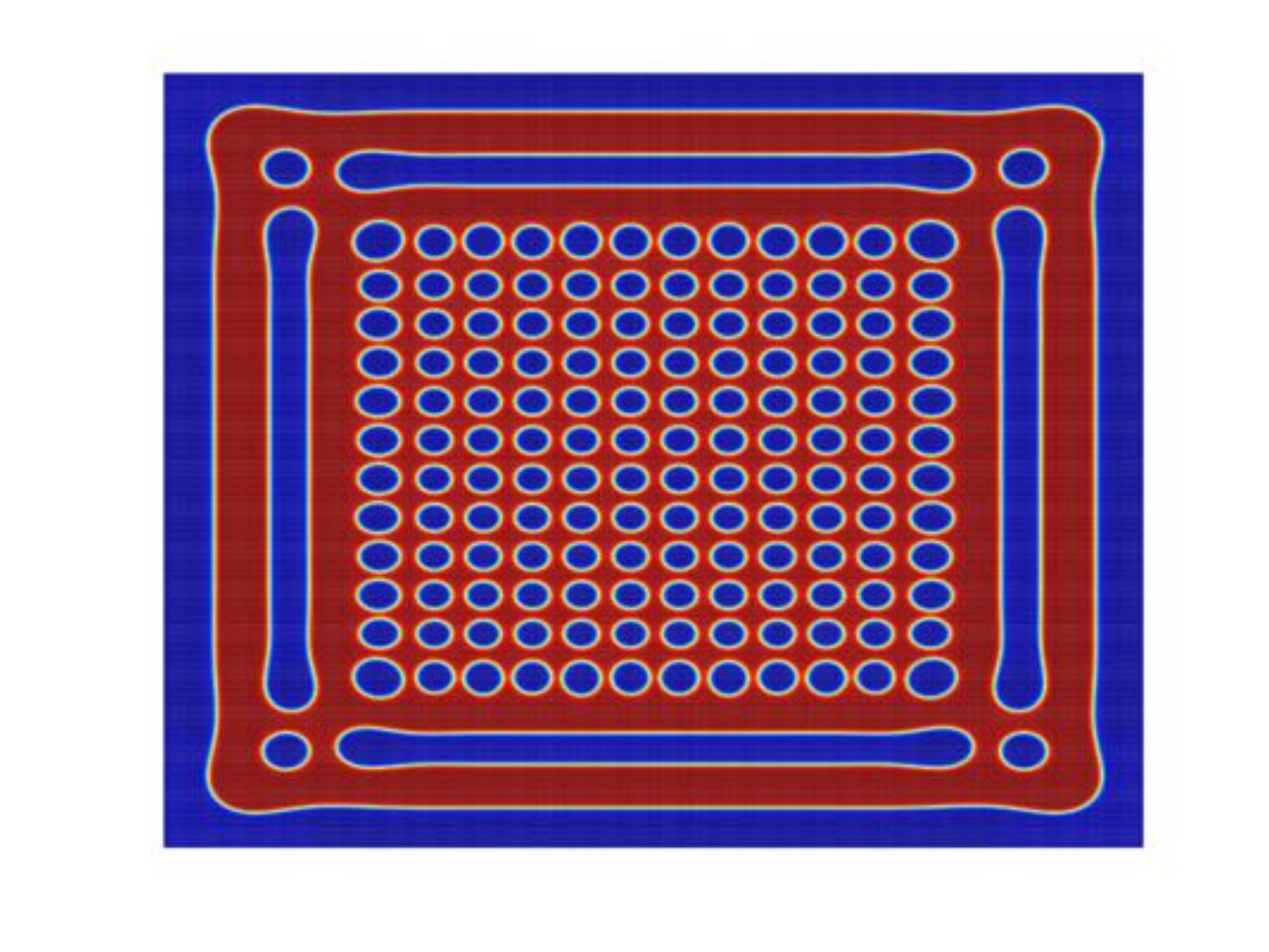}}
    \subfigure[$t=40$]{ \includegraphics[scale=.23]{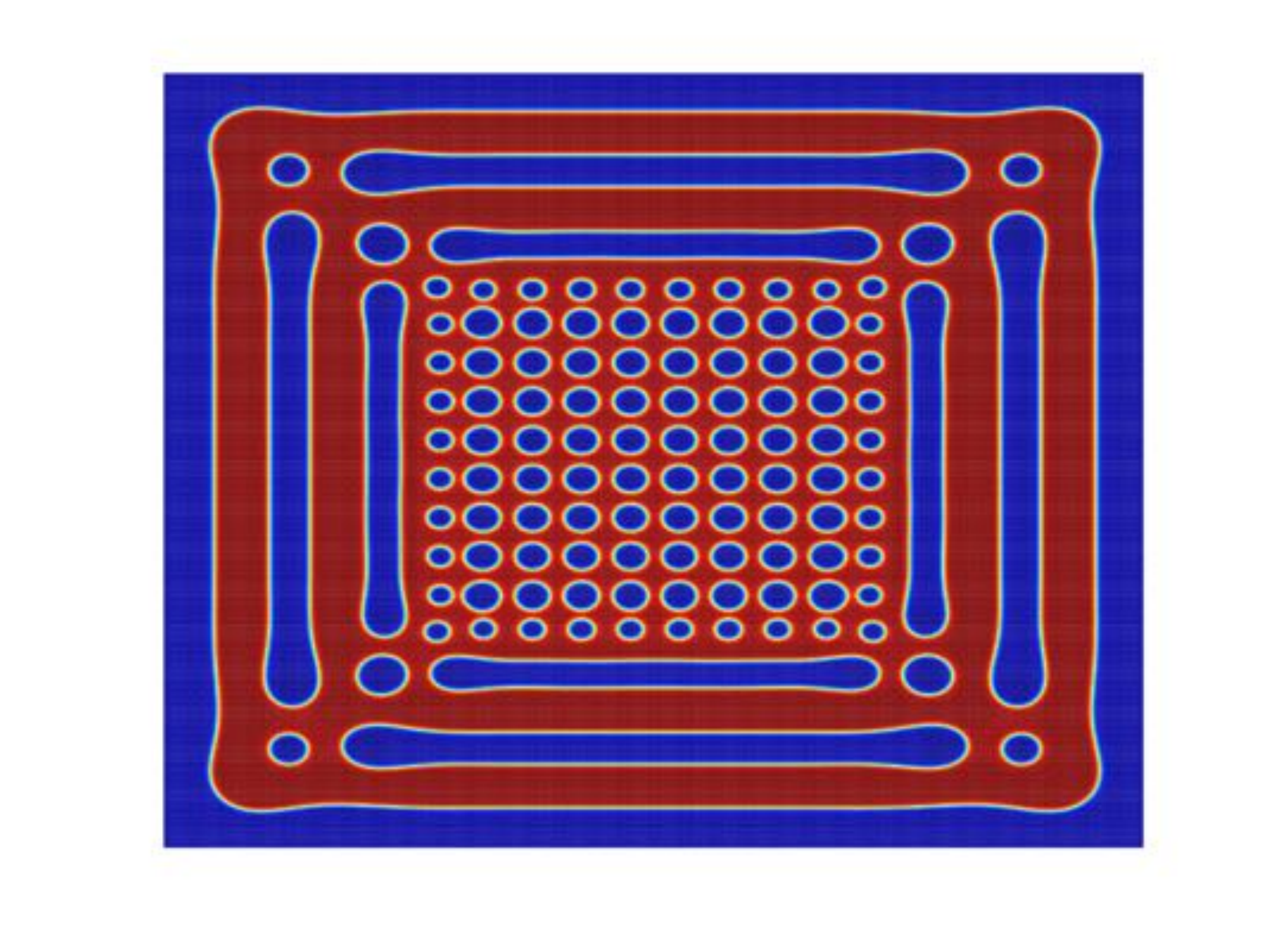}}
      \subfigure[$t=50$ ]{ \includegraphics[scale=.23]{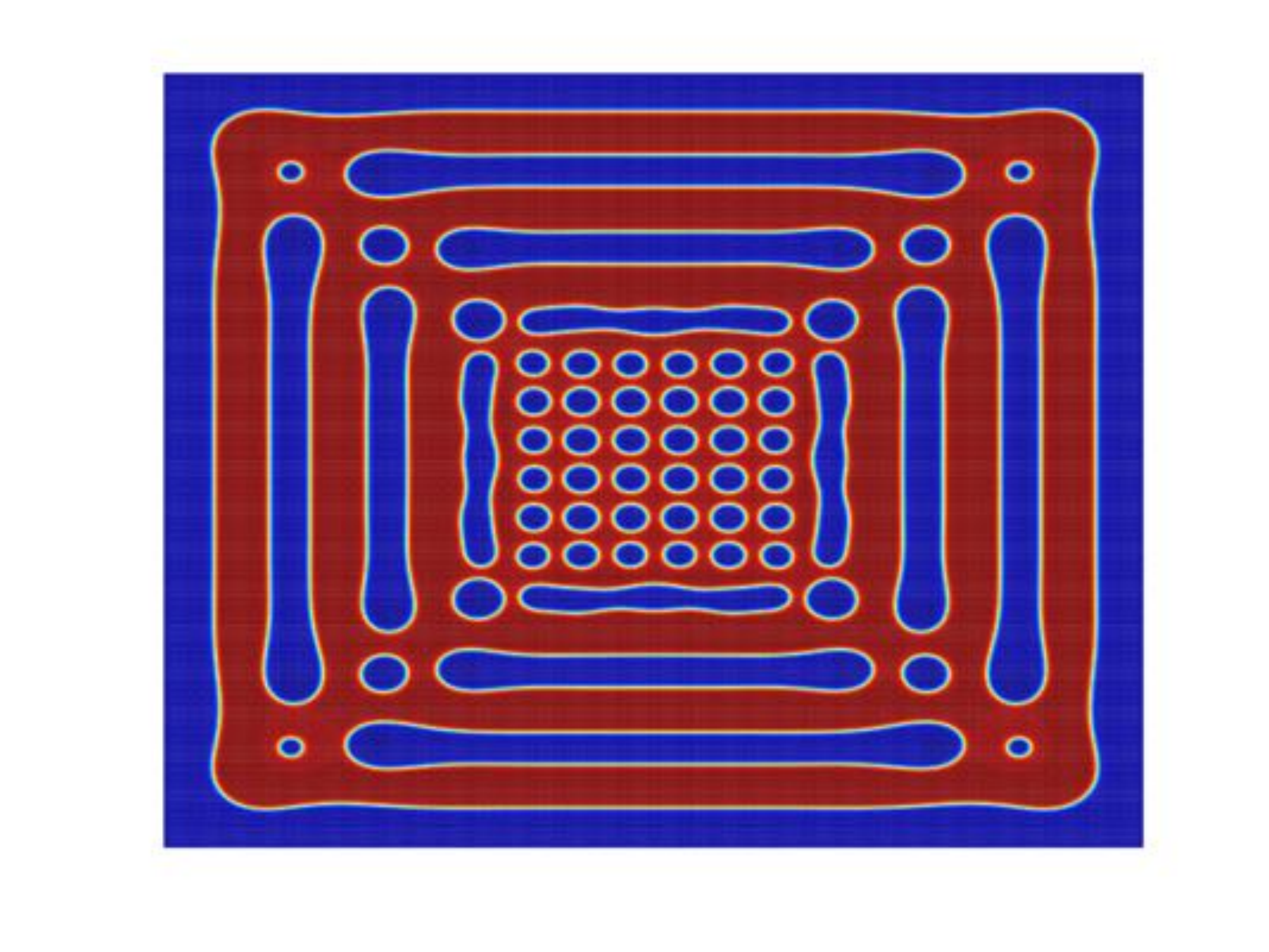}}\\
  \subfigure[$t=60$]{ \includegraphics[scale=.23]{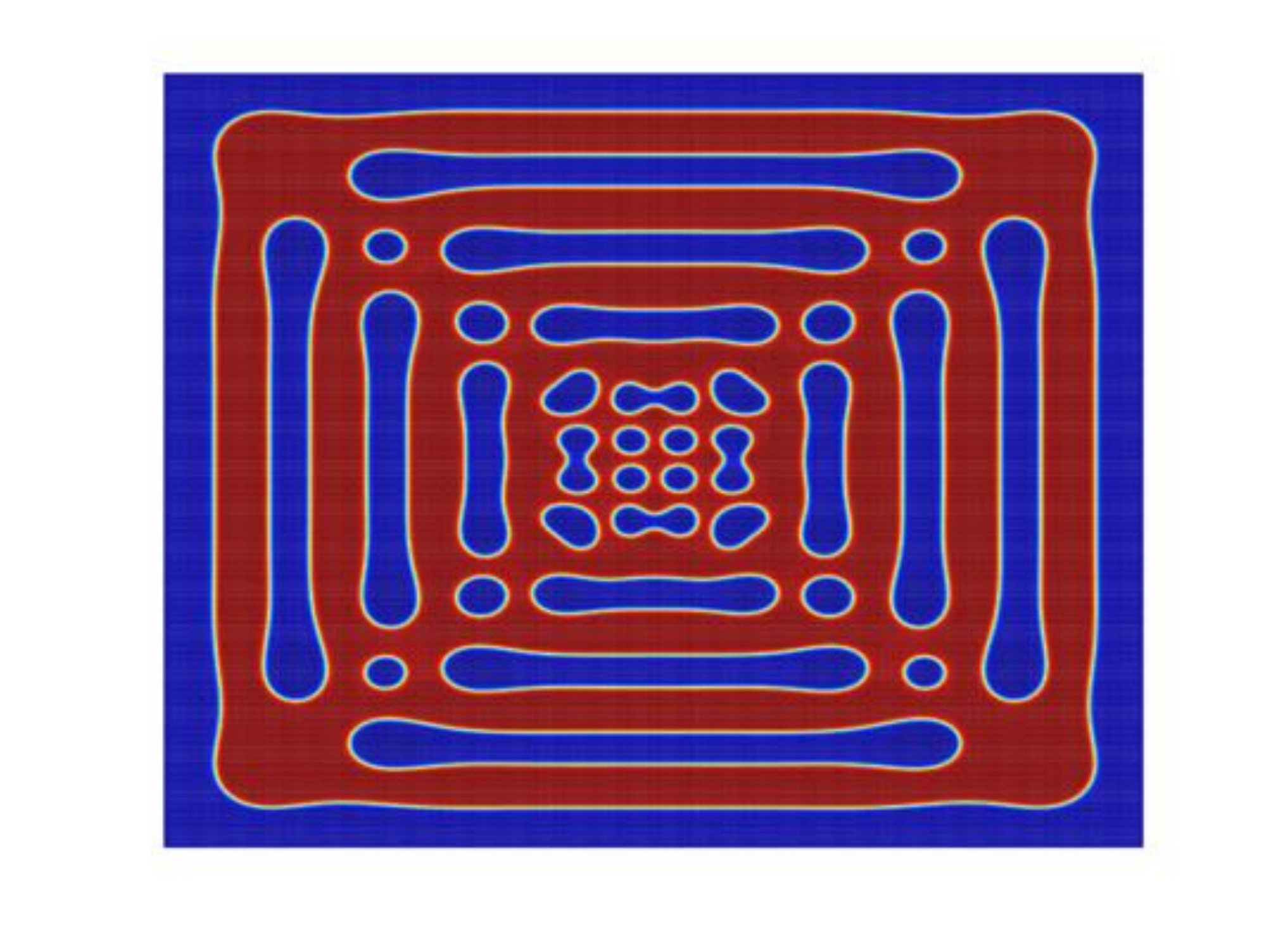}}
   \subfigure[$t=70$]{ \includegraphics[scale=.23]{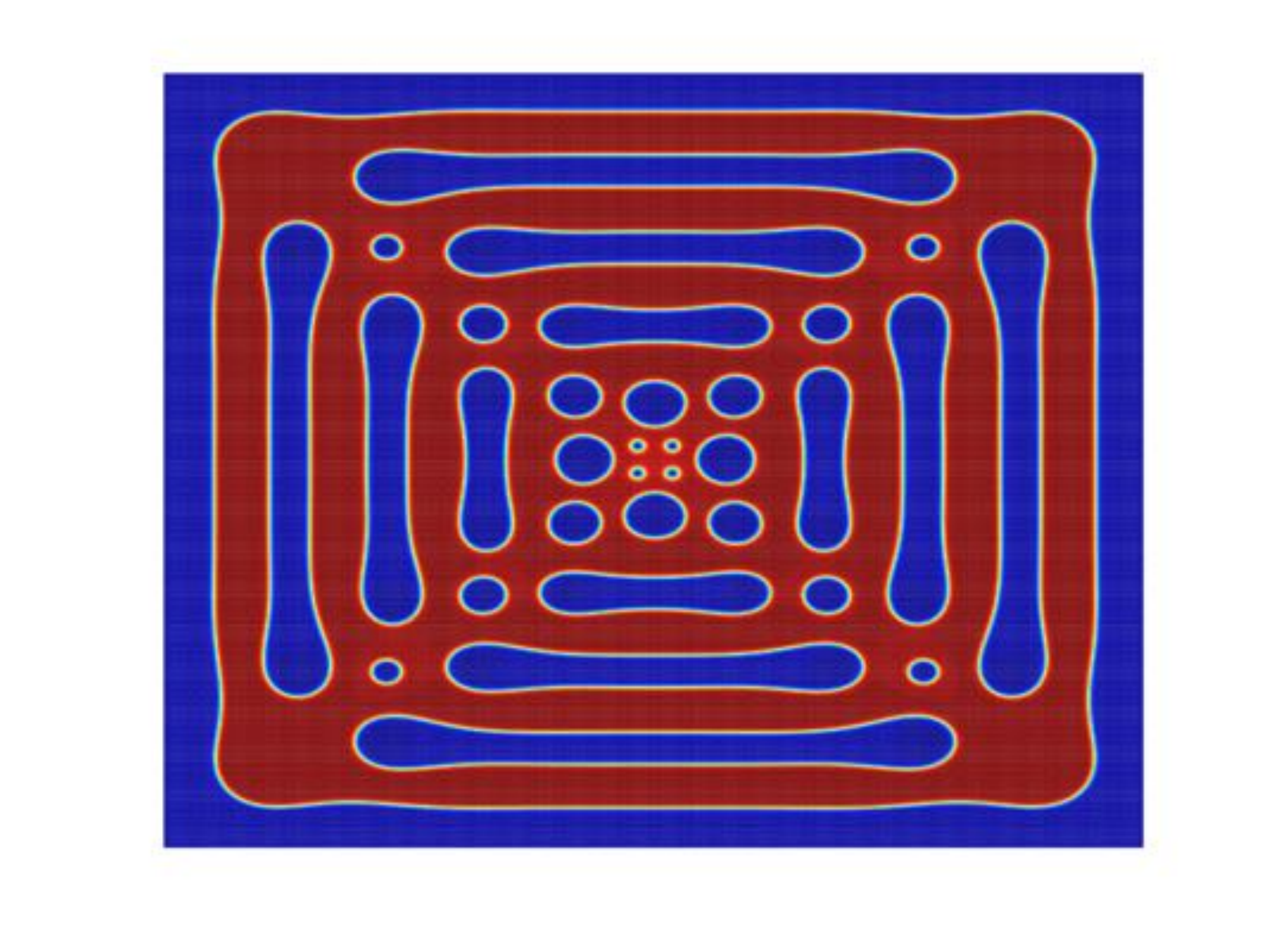}}
    \subfigure[$t=80$]{ \includegraphics[scale=.23]{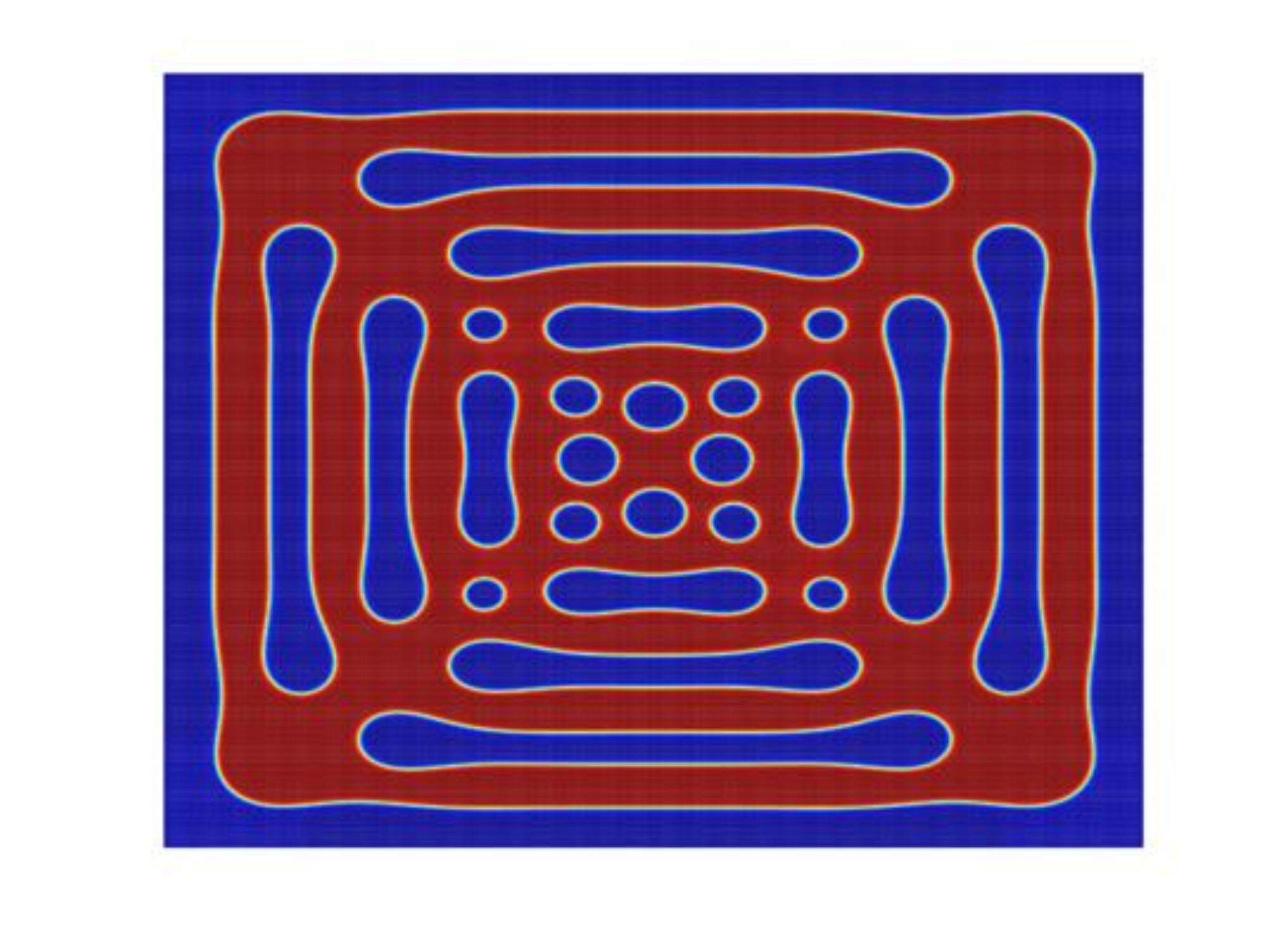}}
    \caption{\small Temporal sequence of snapshots showing the evolution
      of an array of $361$ drops governed by the Cahn-Hilliard equation.
      Simulation results are obtained using Scheme 2A 
      with $\Delta t=10^{-3}.$
    }
   \label{fig:CHEfield}
\end{center}
\end{figure}

Another test problem we would like to consider is the evolution and interaction of $361$
circular drops of one material,
with their centers arranged on a $19\times 19$ grid (see Figure \ref{fig:CHEfield}(a)),
which are immersed in another material.
We assume that the evolution of the material regions is described by the Cahn-Hilliard equation.
The computational domain is taken to be $[0,4]\times[0,4]$,
and the initial phase field distribution is given by
\begin{equation}\label{eq: CHEini}
\phi_0(\bs x,t=0)=360-\sum_{i=1}^{19}\sum_{j=1}^{19} \frac{\tanh \big( \sqrt{(x-x_i)^2+(y-y_j)^2}-R_0 \big)}{\sqrt{2}\eta},
\end{equation}
where $R_0$ is the initial drop radius with $R_0=0.085,$
and $x_i=0.2\times i$ and $y_j=0.2\times j$ for $ i,\,j=1,2,...,19$.
We employ $512$ grid points in both $x$ and $y$ directions in
the Fourier spectral discretization. The other simulation parameters are
$m_0=10^{-6},$ $\sigma=151.15,$ $\eta=0.01$, $\beta=\frac{3}{2\sqrt{2}}\sigma \eta$ and $c_0=1$.
We set $f(\bs x,t)=0$ in \eqref{eq:CHEstd} and periodic boundary conditions
are prescribed on the domain boundaries in both directions.

\begin{figure}
  \begin{center}
    \subfigure[$t=0.15$]{\includegraphics[scale=0.23]{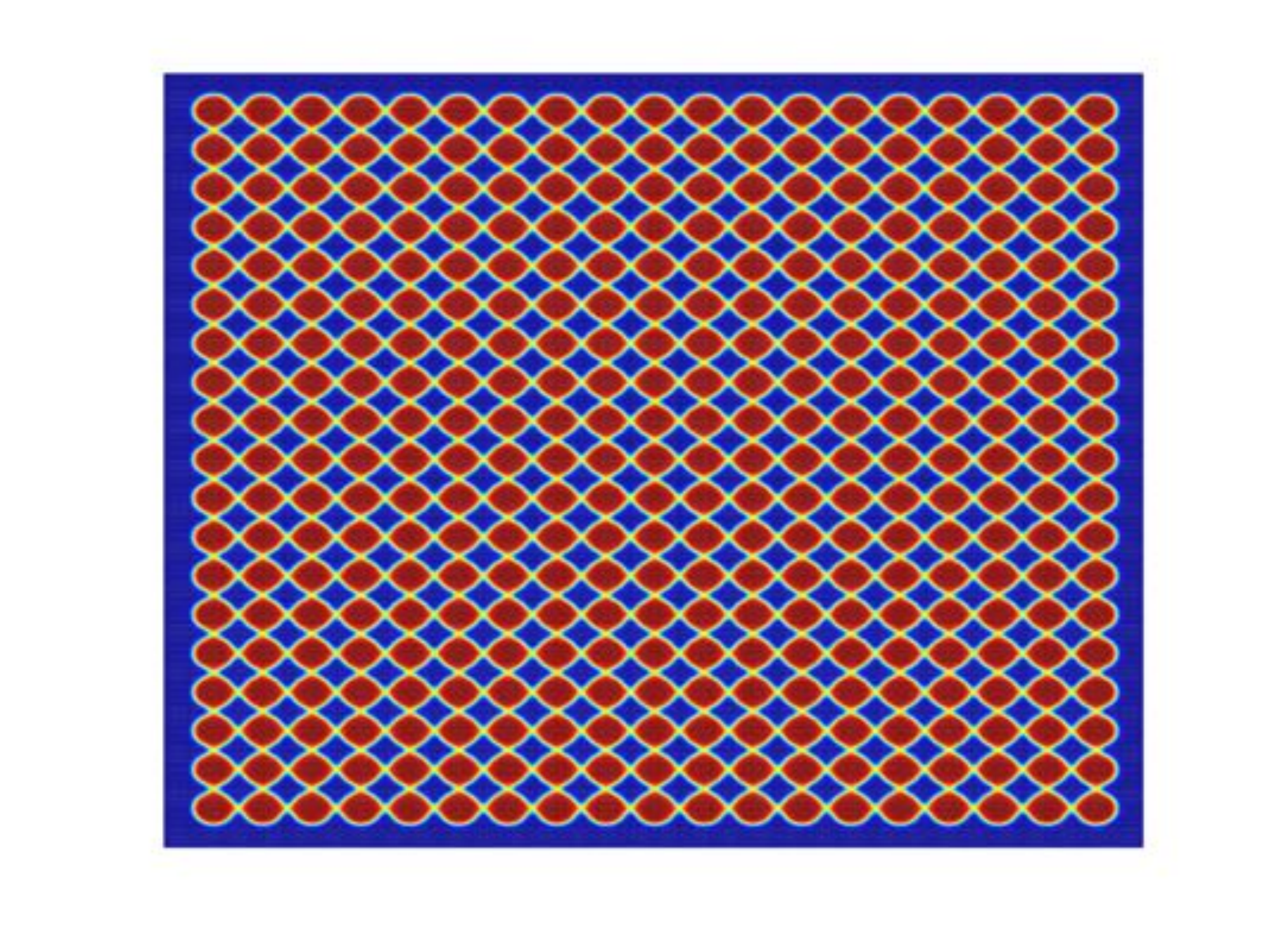}}
    \subfigure[$t=0.2$]{\includegraphics[scale=0.23]{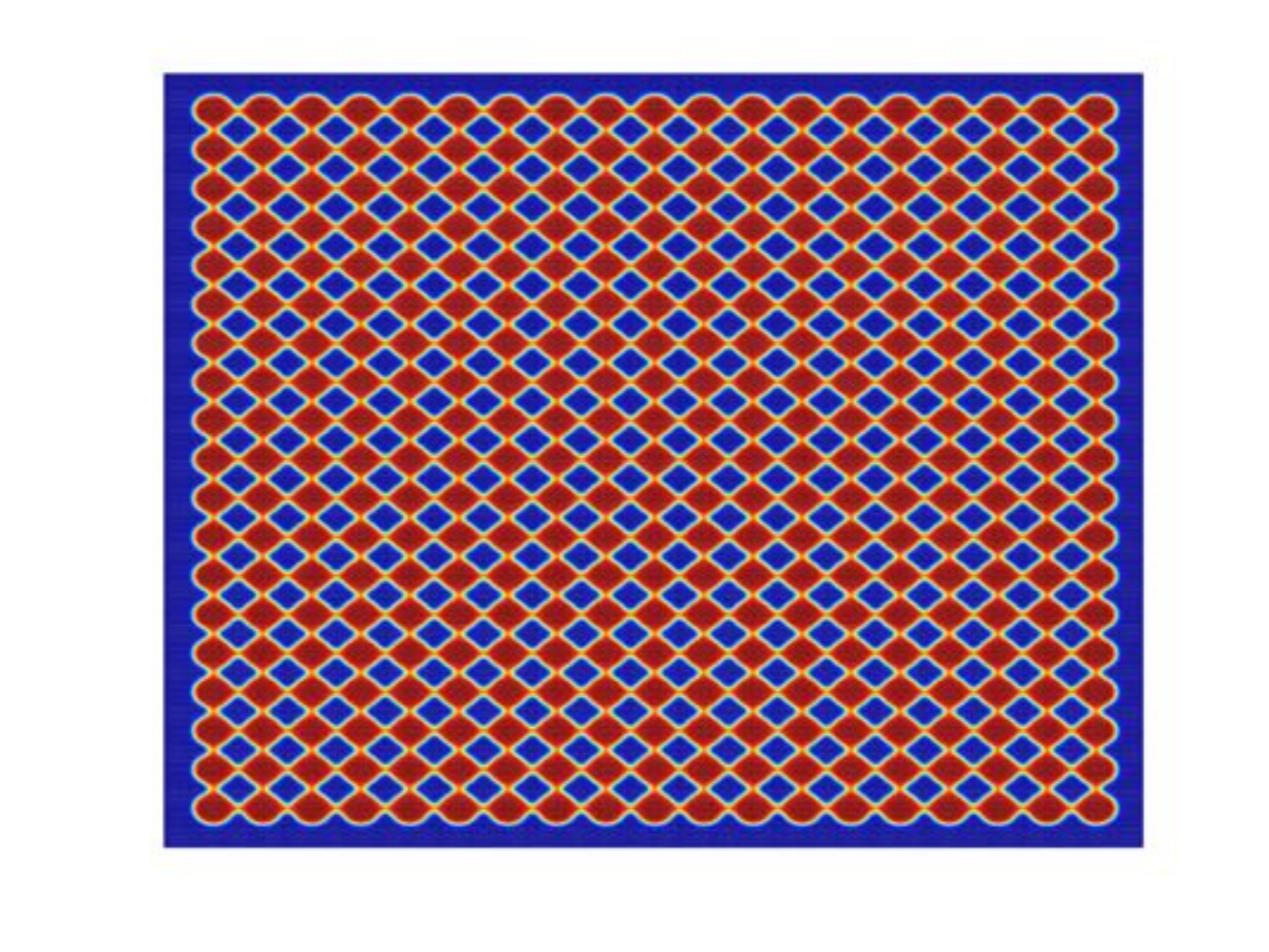}} \\
    \subfigure[$t=0.3$]{\includegraphics[scale=0.23]{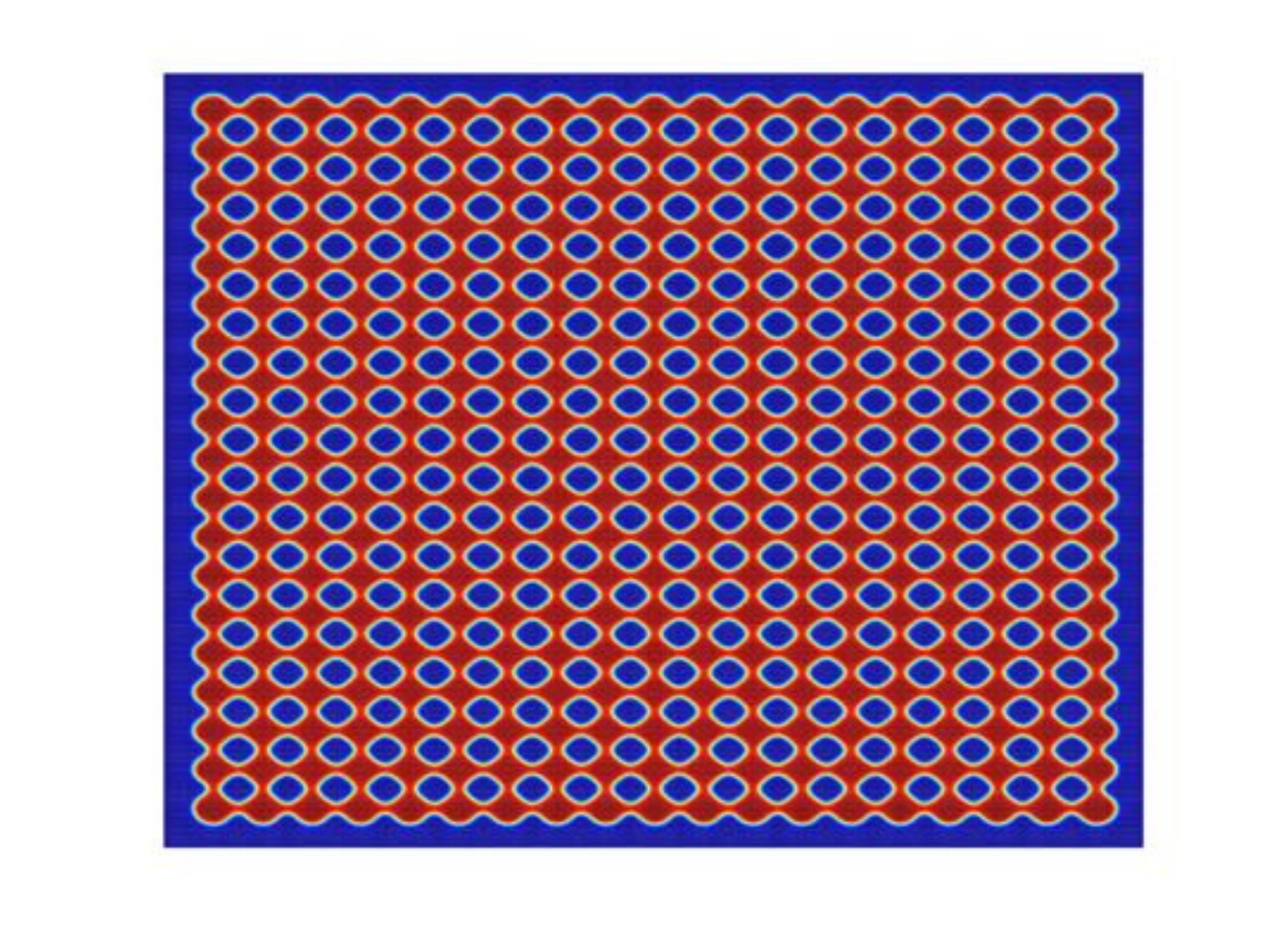}}
    \subfigure[$t=1.0$]{\includegraphics[scale=0.23]{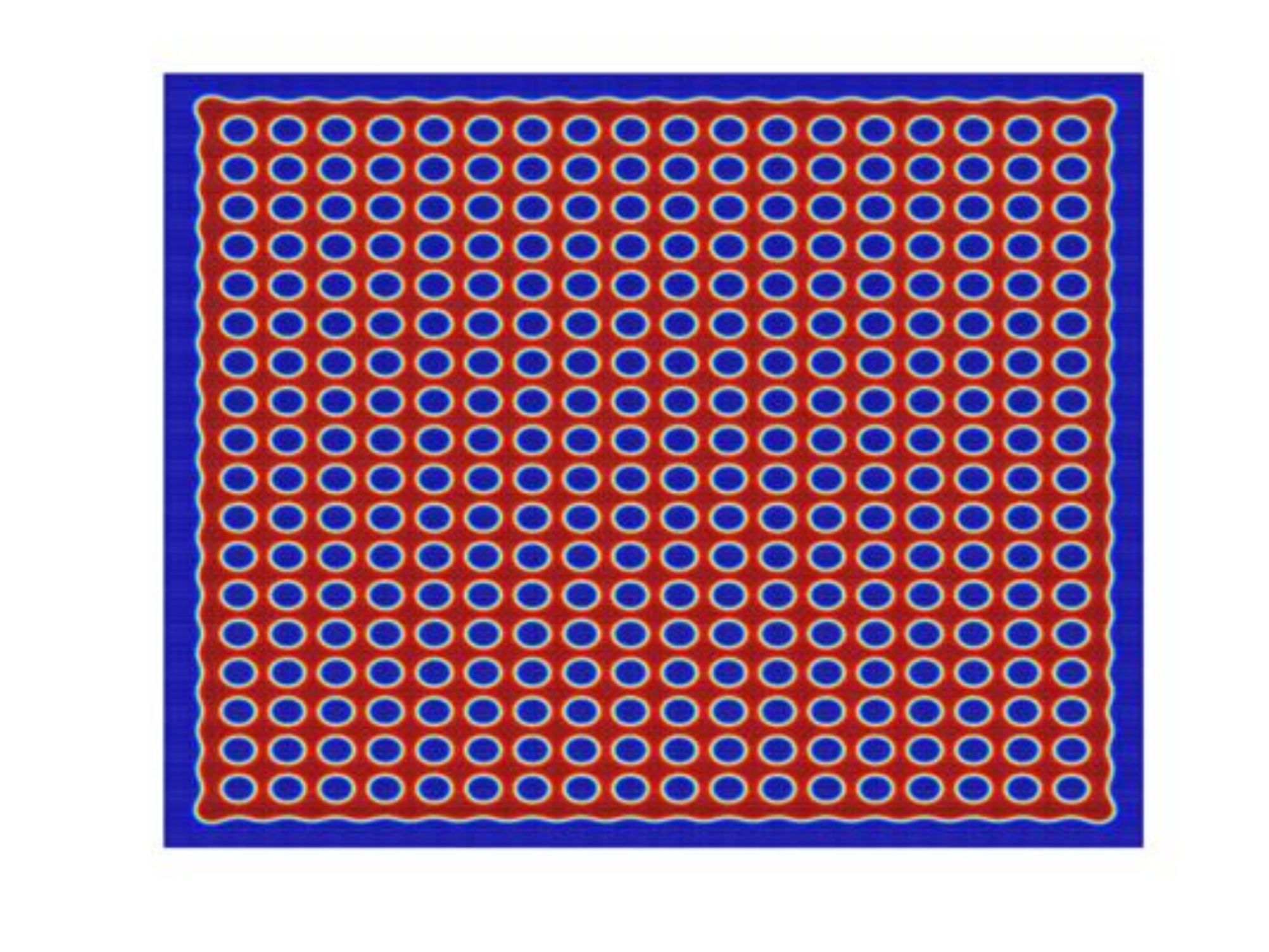}}
  \end{center}
  \caption{
    Early evolution of the two material regions: (a) $t=0.15$,
    (b) $t=0.2$, (c) $t=0.3$, and (d) $t=1.0$. The initial drops of
    material one (red) merge to form the new background material.
    The space between the initial drops filled with material two (blue)
    form a new $18\times 18$ array  of blue drops
    in the new background material (red).
  }
  \label{fig:field_early}
\end{figure}

The regions for the two materials (circular drops, and the background)
are observed to evolve and coalesce to form coarser regions.
This process is visualized in Figure \ref{fig:CHEfield}
with a long temporal sequence of snapshots of the phase field distributions
obtained using Scheme 2A with $\Delta{t}=10^{-3}.$ The first material is marked
by red and the other material is marked by blue. Increasingly coarser
regions can be observed to form over time.
Comparison between Figures \ref{fig:CHEfield}(a) and (b)
indicates that the roles (foreground/background) of the two materials
seem to have reversed early in the evolution. The first material (initial
red drops) evolves into a new background material, while
the second material (initial blue background) form blue drops in
the red background; see Figure \ref{fig:CHEfield}(b).
This process is illustrated in Figure \ref{fig:field_early}
with four snapshots at the early stage of the evolution.
We can observe that the initial $19\times 19$ array of red drops coalesce to form
a new background material, while the second material
in the spaces between the red drops
evolves into a new $18\times 18$ array of blue drops immersed in
the red background material.

\begin{figure}[tb]
\begin{center}
  \subfigure[$\Delta{t}=10^{-4}$ ]{ \includegraphics[scale=.23]{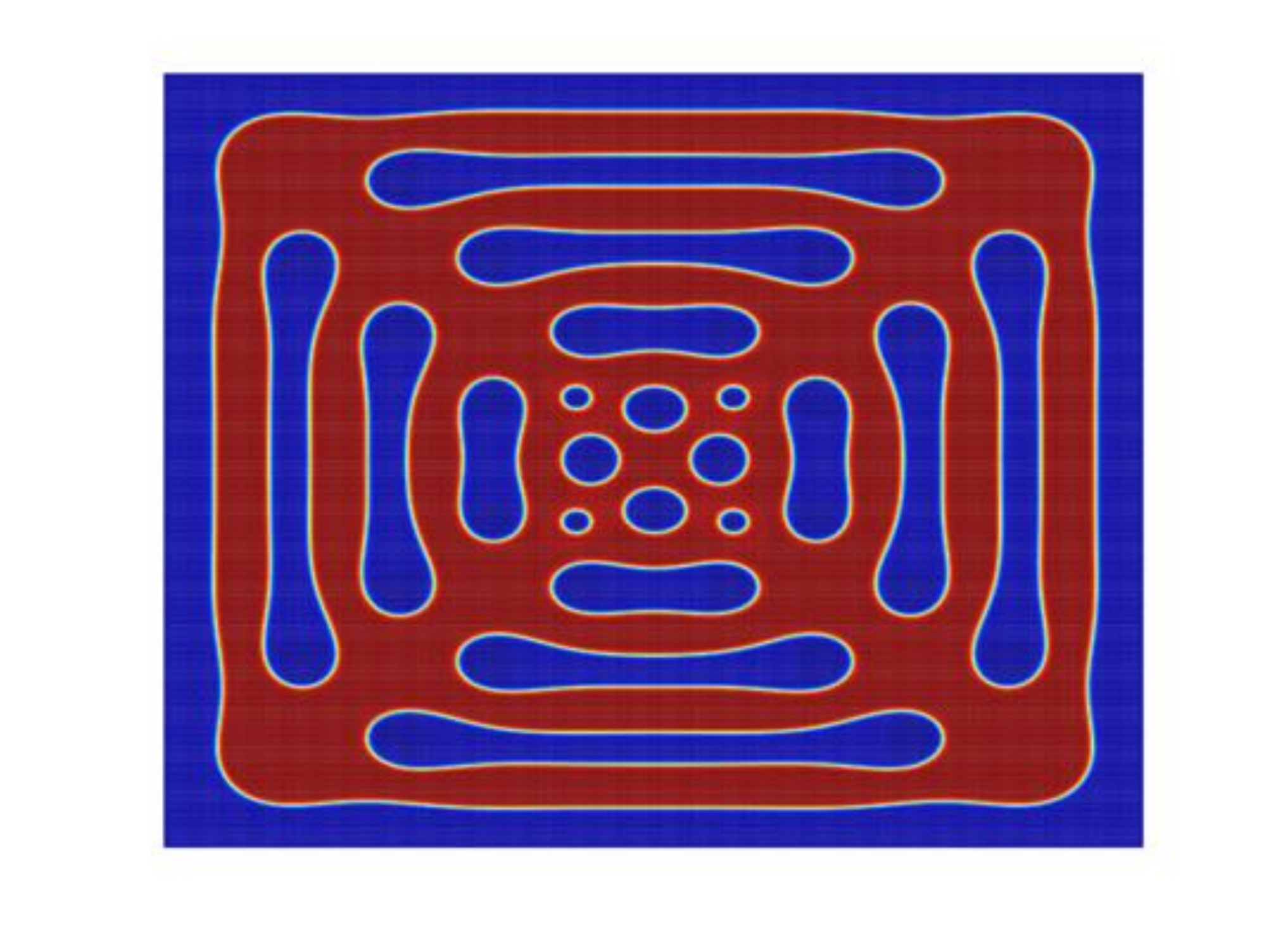}}
  \subfigure[$\Delta{t}=10^{-3}$ ]{ \includegraphics[scale=.23]{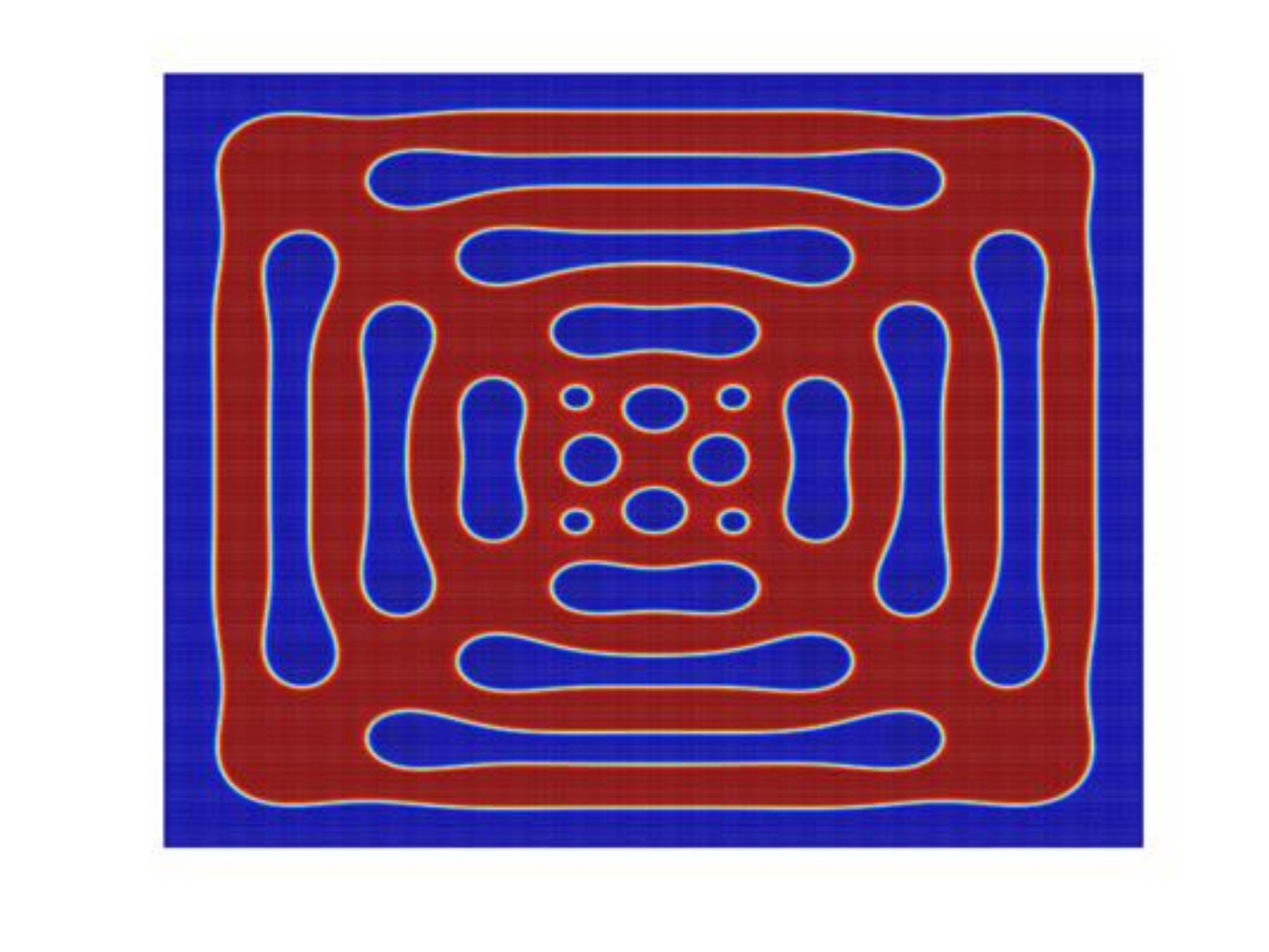}} 
  \subfigure[$\Delta t=10^{-2}$ ]{ \includegraphics[scale=.23]{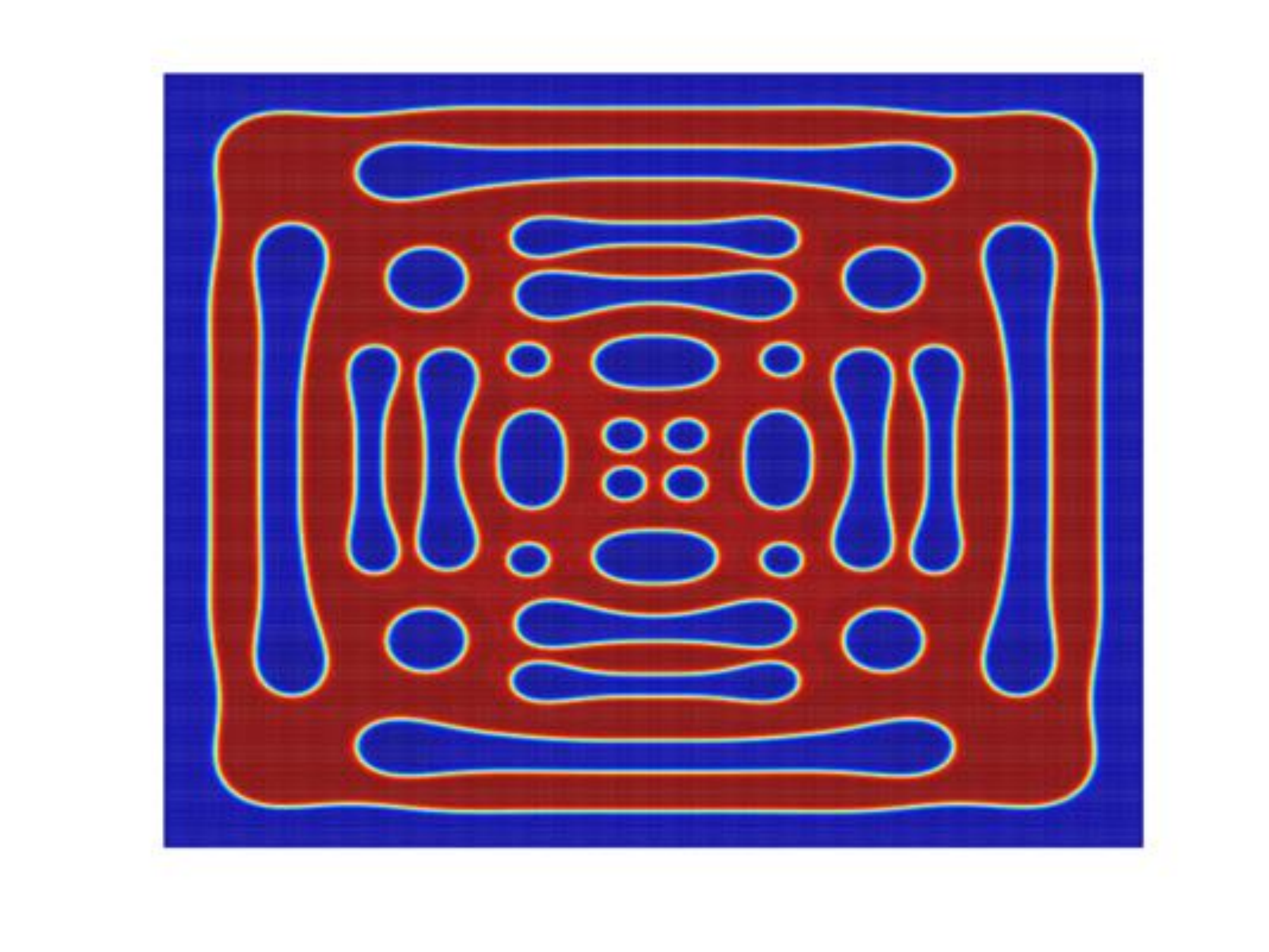}}
  \caption{\small Coalescence of arrays of 361 circles: snapshots of the phase field function at $t=100$ computed using Scheme 2A with (a) $\Delta{t}=10^{-4}$, (b) $\Delta{t}=10^{-3}$, (c) $\Delta{t}=10^{-2}.$}
   \label{fig: CHEfieldcompare}
\end{center}
\end{figure}

The distribution of the material
interface at $t=100$ obtained with several time step sizes,
ranging from $\Delta{t}=10^{-4}$ to $\Delta{t}=10^{-2}$, computed using Scheme 2A
are shown in Figure \ref{fig: CHEfieldcompare}.
It is observed that
the results obtained with $\Delta{t}=10^{-4}$ and $\Delta{t}=10^{-3}$ are essentially the same.
With the larger time step size $\Delta{t}=10^{-2},$ we can observe some
differences in the material distribution from those obtained using
smaller $\Delta{t}$ values, indicating that the simulation starts to lose accuracy
with this step size.

\begin{figure}[tbp]
\begin{center}
  \subfigure[Scheme 2B: $\xi$ vs $\Delta{t} $ ]{ \includegraphics[scale=.32]{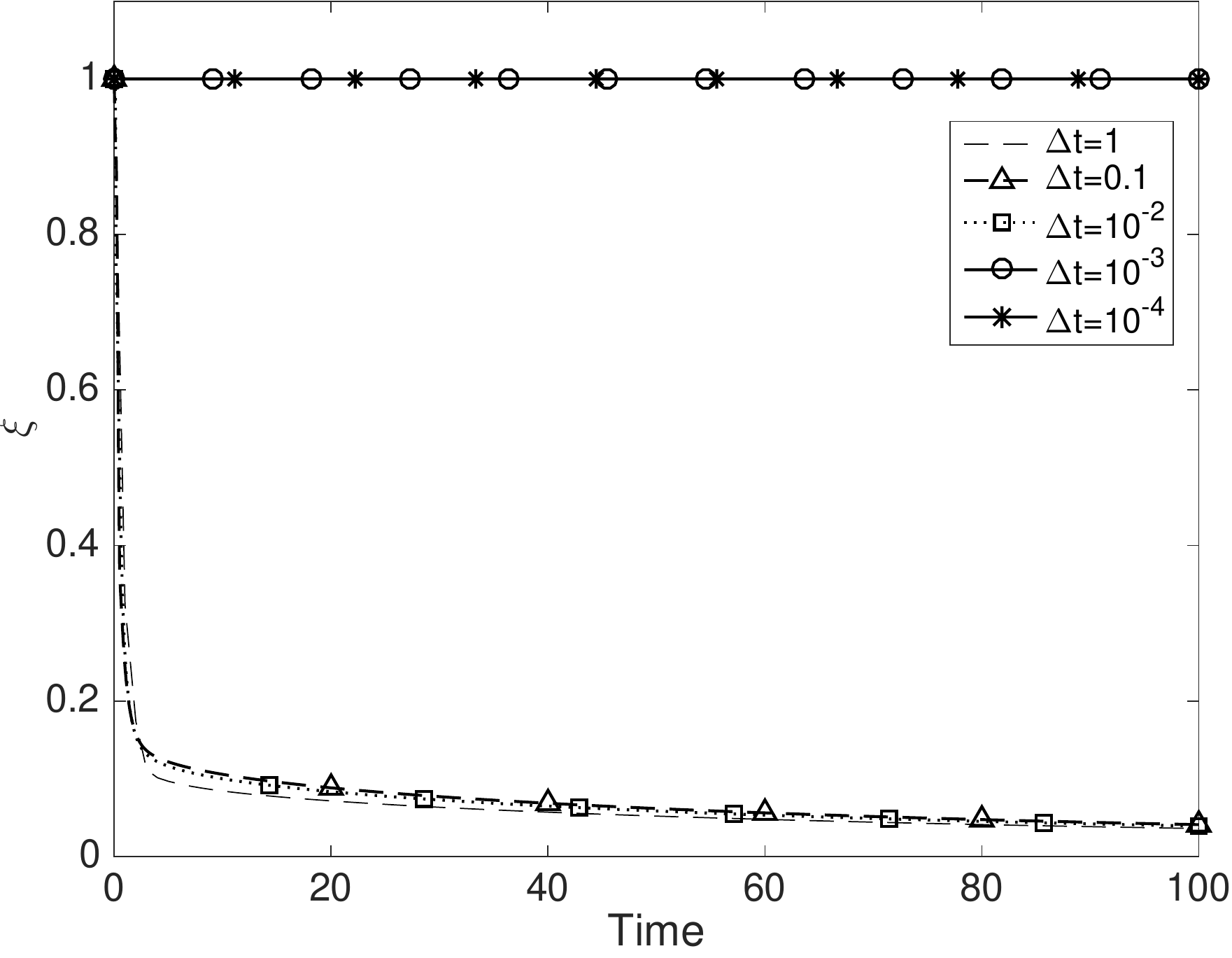}}\qquad
  \subfigure[$\xi$ with $\Delta{t}=10^{-3}$]{ \includegraphics[scale=.32]{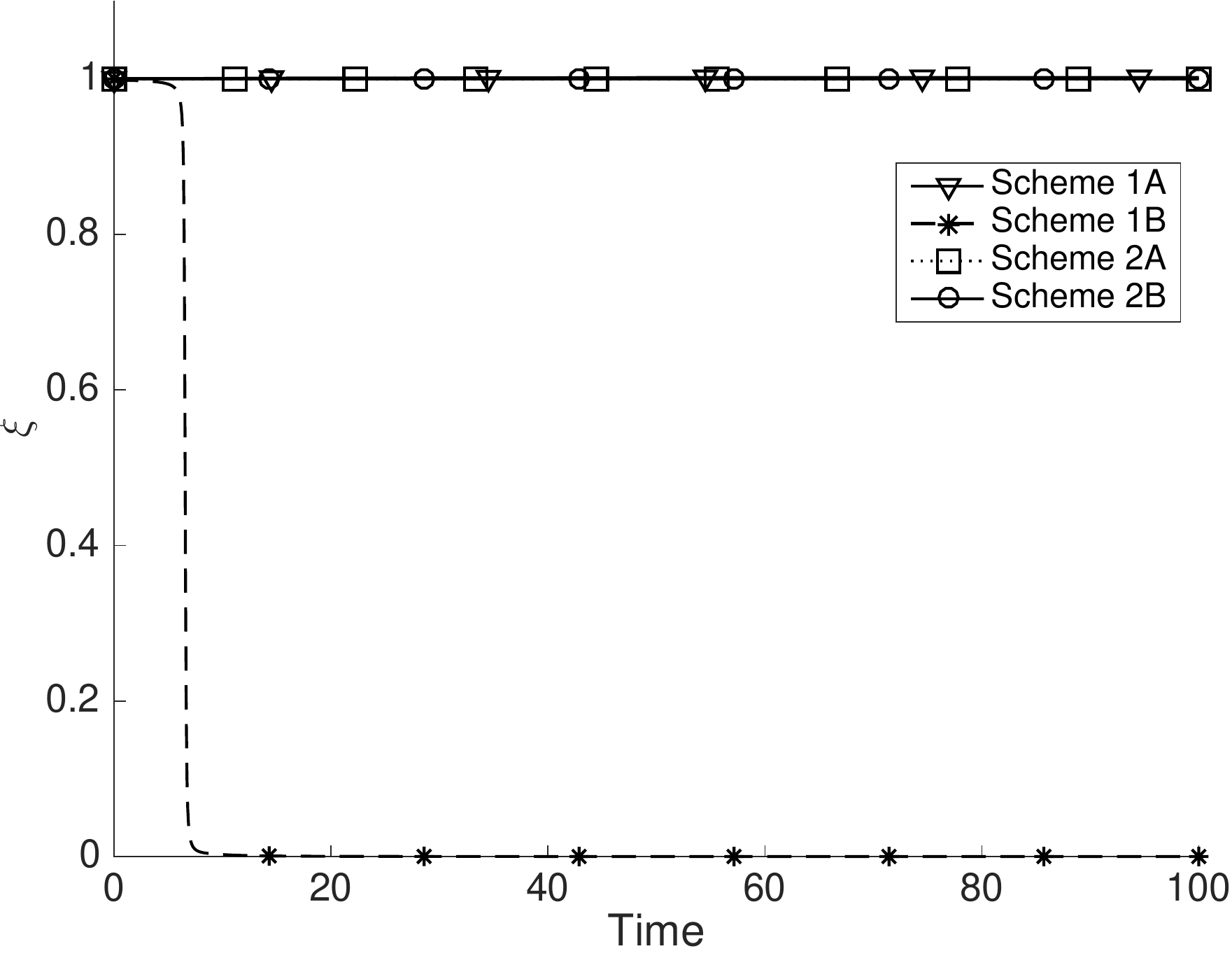}}
  \caption{\small  Coalescence of an array of drops: time histories  of $\xi(t)=R(t)/\sqrt{E(t)}$ corresponding to (a) a range of time step sizes
    $\Delta{t}=1,10^{-1}, 10^{-2}, 10^{-3}, 10^{-4}$, computed using
    Scheme 2B,
    and (b) computed using different schemes (Schemes 1A/1B, 2A/2B)
    with a fixed $\Delta{t}=10^{-3}.$  
   \label{fig: xi}}
\end{center}
\end{figure}

Note that the quantity $\xi=\frac{R(t)}{\sqrt{E(t)}}$ is an approximation of the unit value. This $\xi$ can serve as an indicator of the accuracy of the simulations. If the deviation of $\xi$ from the unit value is small, then the simulation tends to be more accurate. In Figure \ref{fig: xi}(a), we depict the time histories of $\xi$,
computed using Scheme 2B  with various time
step sizes ranging from $\Delta{t}=1$  to $10^{-4}$. It can be observed that $\xi$ remains close to 1 for small time steps $10^{-3} \sim 10^{-4}.$ While for relatively larger time step sizes $1 \sim 10^{-2},$ $\xi$ exhibits an obvious deviation from 1, suggesting
that the simulation results are no longer accurate.  In Figure \ref{fig: xi}(b), we compare the time histories of $\xi$,
obtained using the four schemes (1A/1B and 2A/2B) with $\Delta{t}=10^{-3}$.
The schemes 1A, 2A and 2B all produce quite accurate simulation results with
this time step size, with the computed $\xi$ taking essentially the unit value.
On the other hand, the $\xi$ computed by Scheme 1B has the unit value
initially, and at about $t=10$ it decreases sharply to a small positive value
(on the order $10^{-6}$) and remains at that level for the rest of the
simulation. These results indicate that among the four schemes
developed here the Scheme 1B might be somewhat inferior in terms of accuracy
to the other schemes under the same conditions.

\begin{figure}[tbp]
\begin{center}
  \subfigure[Scheme 2A]{ \includegraphics[scale=.3]{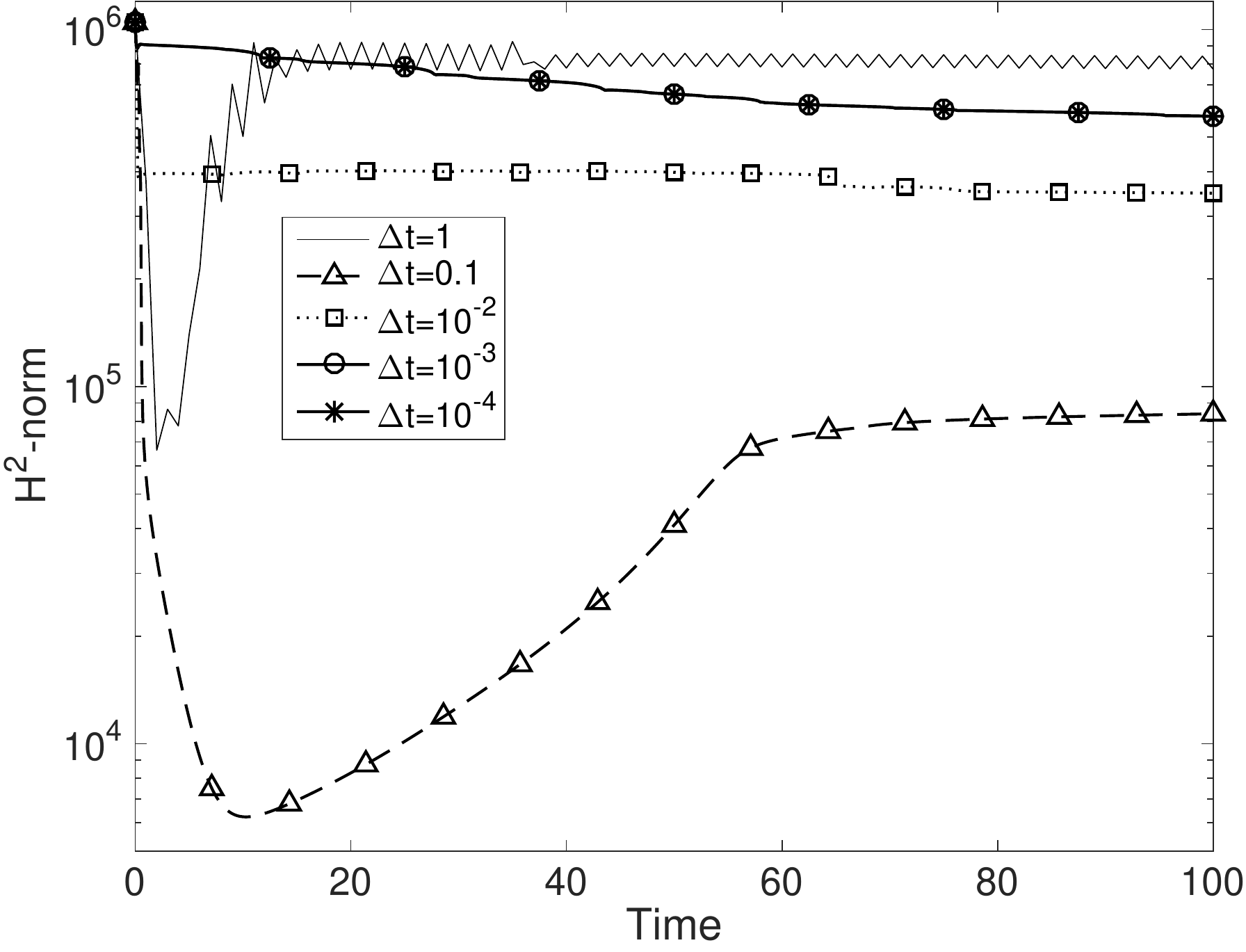}}
 \subfigure[different schemes with $\Delta{t}=1$]{ \includegraphics[scale=.3]{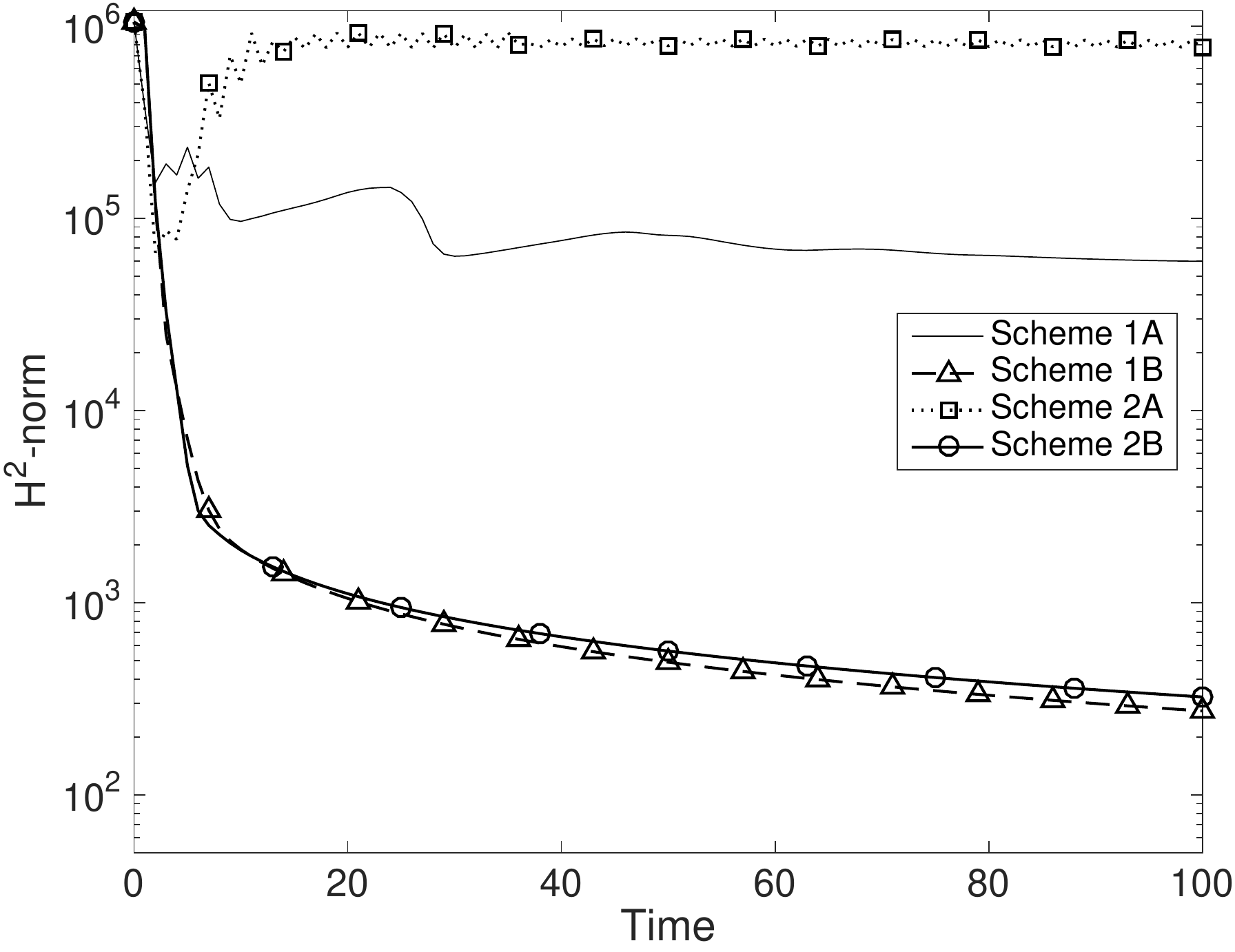}}\qquad
 \caption{\small  Coalescence of an array of drops:
   time histories  of $\| \phi \|_{H^2(\Omega)} $ corresponding to
   (a) a range of time step sizes $\Delta{t}=1,10^{-1}, 10^{-2}, 10^{-3}, 10^{-4}$,
   computed using Scheme 2A,
   and (b) computed using different schemes
   with a fixed $\Delta{t}=1.$
   \label{fig: CHEH2}}
\end{center}
\end{figure}

To validate the stability analysis of the schemes
in the previous section, we look into 
the time histories of the $H^2$ norm of the phase field function $\phi$
in Figure \ref{fig: CHEH2}. 
Figure \ref{fig: CHEH2}(a) shows 
the time histories of the $H^2$-norm of the phase field function $\phi$
corresponding to a number of time step sizes,
ranging from $\Delta{t}=10^{-4}$ to $\Delta{t}=1$,
obtained using Scheme 2A.
It is observed that with smaller $\Delta t$ values
the $H^2$ norm decreases over time, and for larger $\Delta t$ values
it remains approximately
at some constant level over time (except for an initial dip at
the early stage of the simulation). These characteristics
signify the stability of the computations.
In Figure \ref{fig: CHEH2}(b),
we fix $\Delta{t}=1$ and depict the time histories
of the $H^2$-norm of $\phi$ obtained using the different schemes
developed herein.
Since the $\Delta{t}$ is quite large,
we do not expect these simulations to be accurate.
Nonetheless, it can be observed that the $H^2$-norms
are all bounded, indicating the stability
of the proposed schemes.

\begin{figure}[tb]
\begin{center}
  \subfigure[energy history]{ \includegraphics[scale=.35]{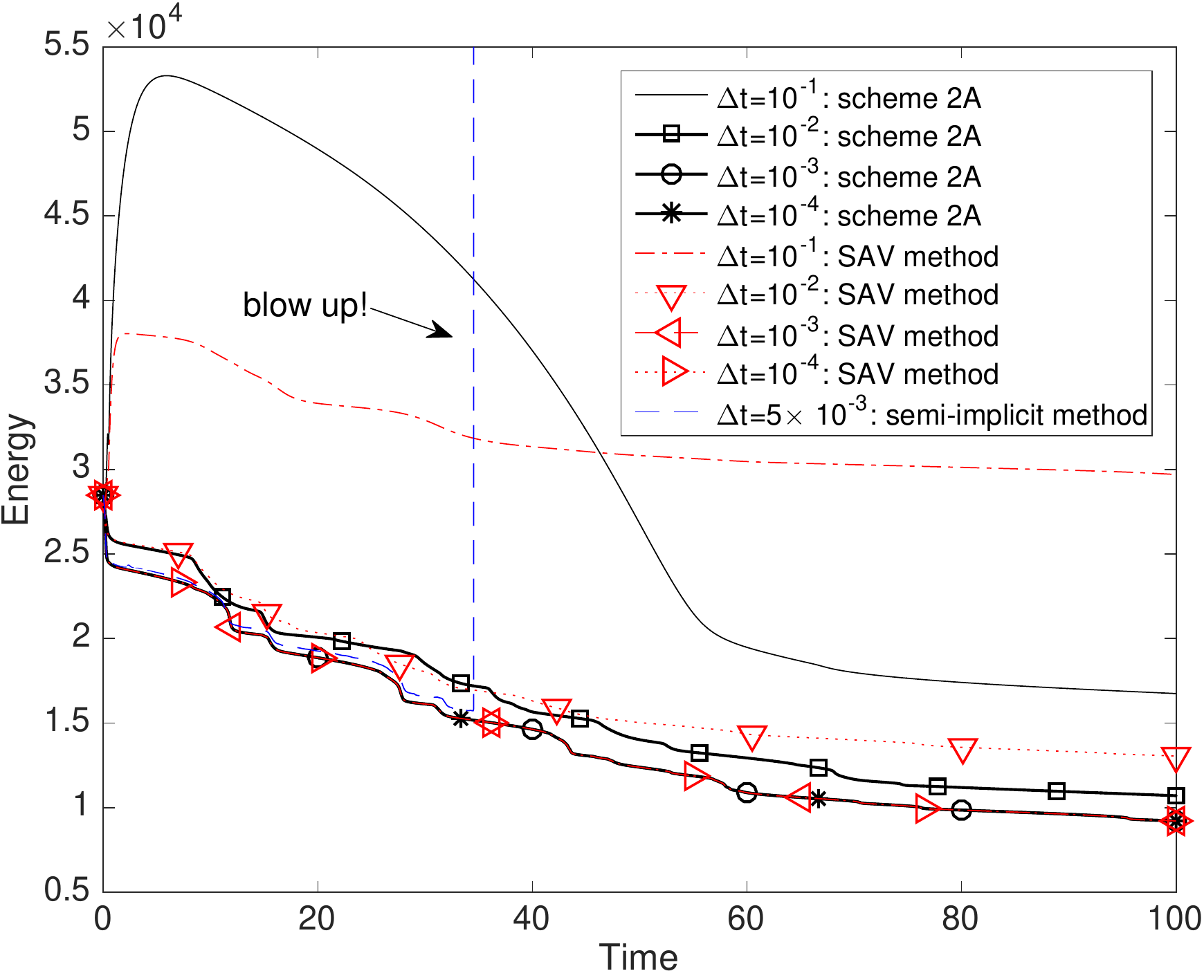}}\\
    \subfigure[history of $R(t)$ with current Scheme 2A]{ \includegraphics[scale=.25]{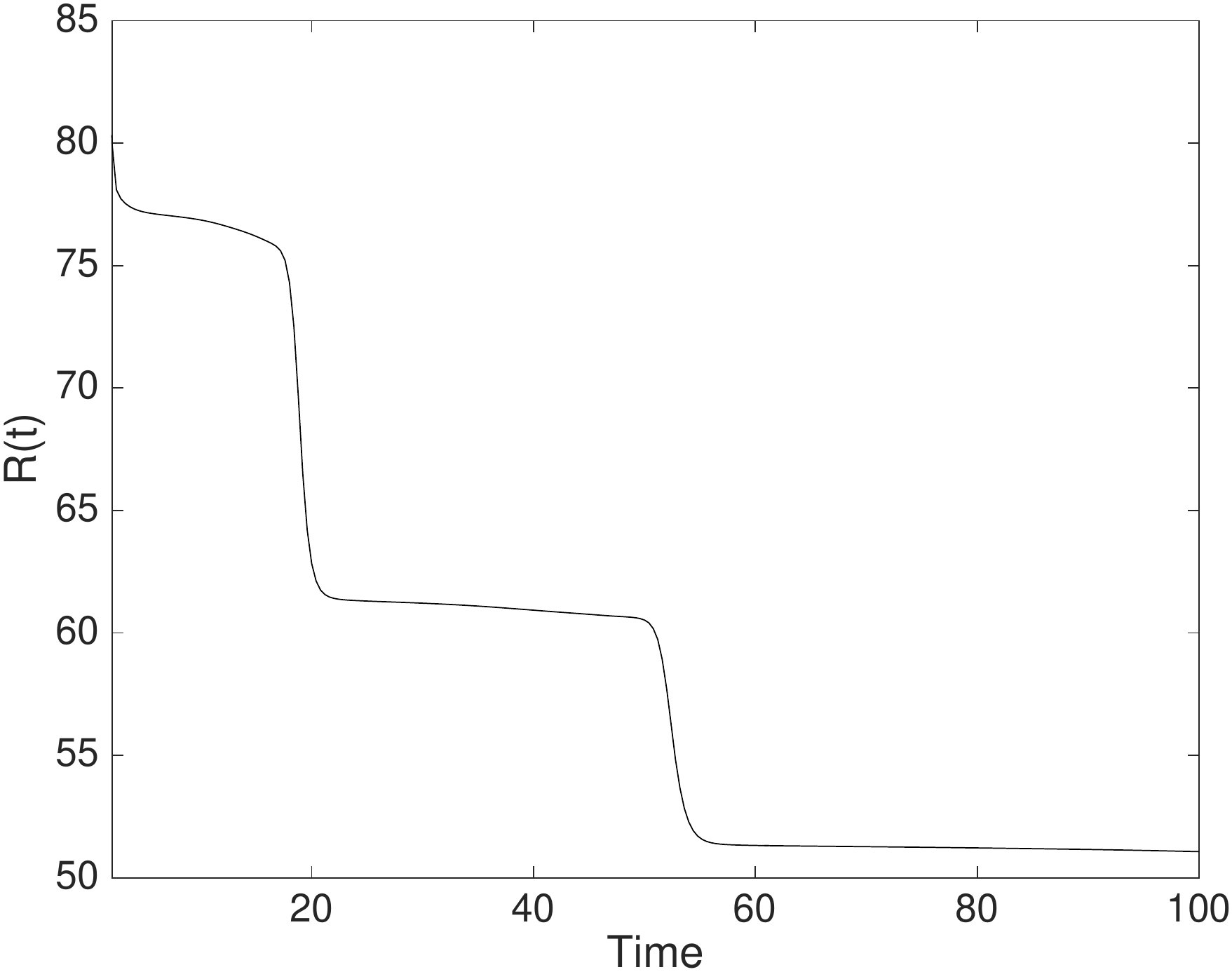}}\;
  \subfigure[history of $R_1(t)$ from the SAV scheme]{ \includegraphics[scale=.25]{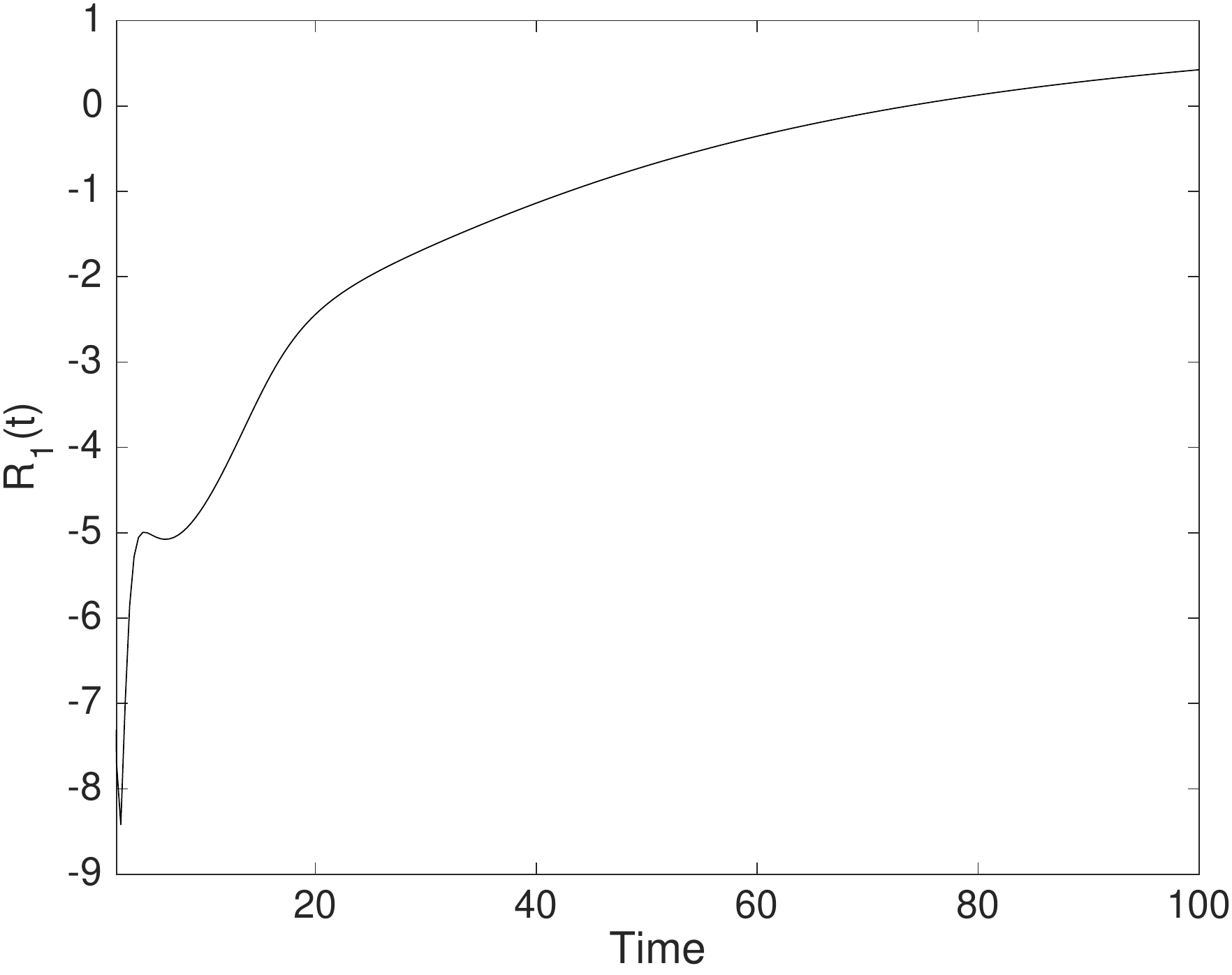}}
  \caption{\small
    Comparison of current gPAV method with other methods (SAV scheme
    and the semi-implicit scheme).
    (a) Time histories of the total energy
    $E(t)$ obtained by the current Scheme 2A and
    the SAV scheme with $\Delta{t}= 10^{-1} - 10^{-4}$, and by the semi-implicit
    scheme from~\cite{DongS2012} with $\Delta{t}=5\times 10^{-3}$.
    (b) Time history of the auxiliary variable $R(t)$ obtained by
    the current Scheme 2A.
    (c) Time history of the auxiliary variable (denoted by $R_1(t)$)
    obtained by the SAV method. 
   \label{fig: CHEEnergy}}
\end{center}
\end{figure}

Finally, we compare the performance of the current schemes
with the SAV method~\cite{ShenXY2018,YangLD2019} and
the semi-implicit scheme~\cite{DongS2012} for the drop evolution
problem. In the SAV method, the auxiliary variable is
defined based on the potential energy only,
$R_1(t)=\sqrt{\int_{\Omega}H(\phi)d\Omega+c_0}$. Here
we use $R_1(t)$ to denote the auxiliary variable in SAV,
in order to distinguish the auxiliary variable $R(t)$
employed here based on the total energy ($R(t)=\sqrt{E(t)}$).
In the semi-implicit scheme~\cite{DongS2012}, the nonlinear term
$h(\phi)$ is simply treated explicitly and the linear terms are
treated implicitly.
Note that in the SAV method, the linear system resulting from
the Cahn-Hilliard equation needs to be solved twice
within a time step~\cite{YangLD2019,ShenXY2018}.
On the other hand, with the schemes proposed here
we only need to solve the linear system  once per time step.
So the operation counts of the current schemes are comparable to
that of the semi-implicit scheme, and are about a half of that
of the SAV method.

Figure \ref{fig: CHEEnergy}(a) shows the time histories of the total energy
$E[\phi]$ obtained using the current Scheme 2A  and
the SAV method with several time step sizes
(ranging from $\Delta{t}=10^{-4}$ to $\Delta{t}=10^{-1}$),
together with the semi-implicit scheme with a time step size
$\Delta{t}=5\times 10^{-3}$.
We note that the history curves corresponding to the relatively small
time step sizes $\Delta{t}=10^{-4}$ and $\Delta{t}=10^{-3}$ all overlap
with one another for both the current scheme and SAV, while some difference between them
can be discerned in the history curves corresponding to
$\Delta{t}=10^{-2}$, suggesting that both the current scheme and SAV
exhibit similar accuracy with relatively small time step sizes.
With the larger step size $\Delta{t}=10^{-1},$ the simulation result is no
longer accurate for both the current scheme and SAV, and indeed we can notice
significant differences when compared with the curves obtained
with smaller time step sizes. Nevertheless, the stability in the computations
with the current scheme and SAV is evident.
This, however, is not the case with the semi-implicit scheme.
Simulation using the semi-implicit scheme blows up after a while
into the computation with a relatively small $\Delta t=5\times 10^{-3}$,
with the energy suddenly growing exponentially. 
These comparisons indicate that the current methods share some
characteristics with SAV in terms of the accuracy and stability.
Note that the current methods require the solution of
the linear system only once within a time step.
So their computational cost is about half of the SAV method.   

The schemes proposed here 
guarantee the unconditional
positivity of the computed $R(t)$ values, irrespective of
the time step size. In Figures \ref{fig: CHEEnergy}(b) and (c) , we
show the time history of the auxiliary variable $R(t)$ computed
using the current Scheme 2A  and the history of
the auxiliary variable $R_1(t)$ obtained by the SAV method, both
with $\Delta t=0.4$.
Note that in the current schemes $R(t)$ is computed by a dynamic
equation stemming from the relation $R(t)=\sqrt{E(t)}$,
while in the SAV method the auxiliary variable $R_1(t)$ is computed
by a dynamic equation stemming from the relation
$R_1(t)=\sqrt{\int_{\Omega}H(\phi)d\Omega+c_0}$.
So both auxiliary variables
should be positive physically.
The current schemes indeed guarantee the positivity of $R(t)$.
In contrast, SAV lacks such a property and
the auxiliary variable computed using SAV can
take negative (unphysical) values, which is evident from Figure \ref{fig: CHEEnergy}(c).

\section{Concluding Remarks}
\label{sec:summary}


In this paper we have presented two first-order
and two second-order unconditionally energy-stable schemes
for numerically solving the Cahn-Hilliard equation.
The stability properties of these schemes have been investigated in relative detail,
and their error analyses are provided.
Besides the discrete unconditional stability,
these schemes have several other attractive properties:
(i) These are linear schemes, and only linear algebraic systems
with a constant coefficient matrix need to be solved.
(ii) The auxiliary variable (scalar-valued number) involved
in each of these schemes is computed by a well-defined explicit form,
and its value is guaranteed to be positive.
(iii) The computational complexity (operation count or computational cost
per time step) of these schemes
is comparable to that of the semi-implicit schemes, and is
about a half of the  gPAV and SAV schemes.

The proposed schemes allow the use of fairly large time step
sizes in dynamic problems and stable computations can be attained. 
These have been demonstrated by numerical examples.
Thanks to the aforementioned  properties, these schemes 
are computationally efficient and simple to implement.
They can be a useful tool for two-phase and multiphase problems
and materials applications.

\section*{Appendix A. Proofs of Several Theorems}

\vspace{0.1in}
\noindent\underline{\bf Theorem \ref{sec3_Lemma3_1B}:}
  Suppose $\phi_{in}\in H^3(\Omega)$ and the condition \eqref{Assume_H2} holds. The following inequality holds with the scheme \eqref{CH1_weak_1B},
  \begin{align*}
  \|\nabla\phi^{n+1}\|_0^2+\frac{\lambda}{2}\|\phi^{n+1}\|_0^2
  +\frac{\Delta t}{2}\|\nabla\Delta\phi^{n+1}\|_0^2+\lambda\Delta t\sum_{k=0}^{n}\|\Delta\phi^{k+1}\|_0^2+\frac{\Delta t}{2}\sum_{k=0}^{n}\|\nabla\mu^{k+1}\|_0^2\leq \widehat{C}_1,
  \end{align*}
  where $\widehat{C}_1$ is the constant
  as given in Theorem \ref{sec3_Lemma3}.        
\begin{proof}

  Taking the inner product of \eqref{CH1_eq1_1B} and \eqref{CH1_eq2_1B} with $\Delta t \mu^{n+1}$ and $\phi^{n+1}-\phi^{n}$, respectively, we have
	\begin{align}\label{sec3_eq6_1B}
	&\frac{1}{2}(\|\nabla\phi^{n+1}\|_0^2-\|\nabla\phi^{n}\|_0^2+\|\nabla\phi^{n+1}-\nabla\phi^{n}\|_0^2)+\frac{\lambda}{2}(\|\phi^{n+1}\|_0^2-\|\phi^{n}\|_0^2+\|\phi^{n+1}-\phi^{n}\|_0^2)
	\nonumber\\
	&\qquad+\Delta t\|\nabla\mu^{n+1}\|_0^2=-|\xi_{_{1B}}^{n}|^2 (h(\phi^{n}),\phi^{n+1}-\phi^{n})\leq\frac{\Delta t}{2}\|\nabla\mu^{n+1}\|_0^2+\frac{|\xi_{_{1B}}^{n}|^4\Delta t}{2}\|\nabla h(\phi^{n})\|_0^2.
	\end{align}
   Take the inner product of  \eqref{CH1_eq1_1B} and \eqref{CH1_eq2_1B}  with $\Delta\phi^{n+1}$ and $\Delta^2\phi^{n+1}$, respectively. Applying the boundary
        condition \eqref{boundary_1B}, we arrive at
	\begin{align}\label{sec3_eq7_1B}
	&\frac{1}{2}(\|\nabla\phi^{n+1}\|_0^2-\|\nabla\phi^{n}\|_0^2+\|\nabla(\phi^{n+1}-\phi^{n})\|_0^2)+\Delta t\|\nabla\Delta\phi^{n+1}\|_0^2+\lambda\Delta t\|\Delta\phi^{n+1}\|_0^2
	\nonumber\\
	&\qquad =|\xi_{_{1B}}^{n}|^2\Delta t(\nabla h(\phi^{n}),\nabla\Delta\phi^{n+1})\leq \frac{\Delta t}{2}\|\nabla\Delta\phi^{n+1}\|_0^2+\frac{|\xi_{_{1B}}^{n}|^4\Delta t}{2}\|\nabla h(\phi^{n})\|_0^2.
	\end{align}
	Summing up \eqref{sec3_eq6_1B} and \eqref{sec3_eq7_1B}, we have
	\begin{align}\label{sec3_eq8_1B}
	&\|\nabla\phi^{n+1}\|_0^2-\|\nabla\phi^{n}\|_0^2+\|\nabla\phi^{n+1}-\nabla\phi^{n}\|_0^2+\frac{\lambda}{2}(\|\phi^{n+1}\|_0^2-\|\phi^{n}\|_0^2+\|\phi^{n+1}-\phi^{n}\|_0^2)
	\nonumber\\
	&\qquad +\frac{\Delta t}{2}\|\nabla\mu^{n+1}\|_0^2+\frac{\Delta t}{2}\|\nabla\Delta\phi^{n+1}\|_0^2+\lambda\Delta t\|\Delta\phi^{n+1}\|_0^2\leq |\xi_{_{1B}}^{n}|^4\Delta t\|\nabla h(\phi^{n})\|_0^2.
	\end{align}
	As is shown in the proof of Theorem \ref{sec3_Lemma3},
        the nonlinear term $\|\nabla h(\phi^{n})\|_0^2$ can be
        estimated by using the positive auxiliary
        variable $\xi_{_{1A}}^{n+1}$, which satisfies the properties
        \eqref{sec3_xi}
        and \eqref{sec3_eq5}. Notice that $\xi_{1B}^n$ satisfies similar
        properties, i.e.~$\xi_{_{1B}}^{n}\leq C(M)$ and
	\[
	\xi_{_{1B}}^{n}=\frac{R^n}{\sqrt{\bm{E}[\phi^{n}]}}\leq\frac{C(M)}{\|\phi^n\|_1}.
	\]
        The rest of the proof parallel those steps in
        the proof of Theorem \ref{sec3_Lemma3}.
\end{proof}

\vspace{0.1in}
\noindent\underline{\bf Theorem \ref{sec4_Lemma4_1B}:}
 Suppose the condition \eqref{regularity_1}, and
 the conditions of Theorems \ref{sec3_Lemma3_1B} and \ref{sec3_Lemma4_1B} hold.
 We have the following result with sufficiently small $\Delta t$, 
 \begin{equation*}
 \frac{1}{2}\|\nabla e_{\phi}^{n+1}\|_0^2+ \frac{\lambda}{2}\|e_{\phi}^{n+1}\|_0^2+|e_{R}^{n+1}|^2+\frac{\Delta t}{2}\sum_{k=0}^n\|\nabla e_{\mu}^{k+1}\|_0^2\leq\widehat{C}_4\Delta t^2,
 \end{equation*}
 where $\widehat{C}_4=C\exp(\Delta t\sum_{k=0}^{n}\frac{r^{k+1}}{1-r^{k+1}\Delta t})\int_{0}^{t^{n+1}}\left(\|\phi_t(s)\|_{1}^4+\|\phi_t(s)\|_{1}^2+\|\phi_{tt}(s)\|_{-1}^2\right)ds$, $r^k=1+\|\nabla\mu^k\|_0^2$ and the constant $C$ depends on $T$, $\phi_{in}$, $\Omega$, $\|\phi\|_{L^{\infty}(0,T;W^{3,\infty}(\Omega))}$ and $\|\mu\|_{L^{\infty}(0,T;H^{1}(\Omega))}$. 
\begin{proof} By subtracting \eqref{truncation1B} from \eqref{CH1_weak_1B}, we have
	\begin{subequations}\label{error_1B}
		\begin{align}
		\label{error_eq1_1B}
		&\frac{e_{\phi}^{n+1}-e_{\phi}^{n}}{\Delta t}= \Delta e_{\mu}^{n+1}-\frac{1}{\Delta t}T_{\phi_{1B}}^{n+1},\\
		\label{error_eq2_1B}
		&e_{\mu}^{n+1}=-\Delta e_{\phi}^{n+1}+ \lambda e_{\phi}^{n+1} + A_3^{n+1},\\
		\label{error_eq3_1B}
		&\frac{e_{R}^{n+1}-e_{R}^{n}}{\Delta t}= -\frac{1}{2}A_4^{n+1}-\frac{1}{\Delta t}T_{R_{1B}}^{n+1},
		\end{align}
	\end{subequations}
	where 
	\begin{align*}
	A_3^{n+1}&=|\xi_{_{1B}}^{n}|^2h(\phi^{n})-\frac{R(t^{n+1})^2}{\bm{E}[\phi(t^{n+1})]}h(\phi(t^{n+1}))\\
	&=\frac{e_R^{n}(R^{n}+R(t^{n}))}{\bm{E}[\phi^{n}]}h(\phi^{n})+\frac{R(t^{n})^2-R(t^{n+1})^2}{\bm{E}[\phi^{n}]}h(\phi^{n})+
	R(t^{n+1})^2\left(\frac{h(\phi^{n})}{\bm{E}[\phi^{n}]}-\frac{h(\phi(t^{n}))}{\bm{E}[\phi(t^{n})]}\right)\\
	&\quad +R(t^{n+1})^2\left(\frac{h(\phi(t^{n}))}{\bm{E}[\phi(t^{n})]}-\frac{h(\phi(t^{n+1}))}{\bm{E}[\phi(t^{n+1})]}\right),\\
	A_4^{n+1}&=\frac{\xi_{_{1B}}^{n+1}}{\sqrt{\bm{E}[\phi^{n+1}]}}\int_{\Omega}|\nabla\mu^{n+1}|^2d\bm{x}-\frac{R(t^{n+1})}{\bm{E}[\phi(t^{n+1})]} \int_{\Omega}|\nabla\mu(t^{n+1})|^2d\bm{x}\\
	&=\frac{e_{R}^{n+1}}{\bm{E}[\phi^{n+1}]}\int_{\Omega}|\nabla\mu^{n+1}|^2d\bm{x}+
	\frac{R(t^{n+1})}{\bm{E}[\phi^{n+1}]}\int_{\Omega}(|\nabla\mu^{n+1}|^2-|\nabla\mu(t^{n+1})|^2)d\bm{x}\\
	& \quad + R(t^{n+1})\left(\frac{1}{\bm{E}[\phi^{n+1}]}-\frac{1}{\bm{E}[\phi(t^{n+1})]}\right)\int_{\Omega}|\nabla\mu(t^{n+1})|^2d\bm{x}.
	\end{align*}
	Taking the inner product of \eqref{error_eq1_1B} with $\Delta te_{\mu}^{n+1}$ and \eqref{error_eq2_1B} with $e_{\phi}^{n+1}-e_{\phi}^{n}$, and multiplying \eqref{error_eq3_1B} with $2\Delta te_{R}^{n+1}$, we get the  following:
	\begin{subequations}\label{sec4_eq7_1B}
		\begin{align}
		\label{sec4_7_eq1_1B}
		&\frac{1}{2}(\|\nabla e_{\phi}^{n+1}\|_0^2-\|\nabla e_{\phi}^{n}\|_0^2+\|\nabla e_{\phi}^{n+1}-\nabla e_{\phi}^{n}\|_0^2)+\frac{\lambda}{2}(\|e_{\phi}^{n+1}\|_0^2-\|e_{\phi}^{n}\|_0^2
		\nonumber\\
		&\qquad+\|e_{\phi}^{n+1}-e_{\phi}^{n}\|_0^2)+\Delta t\|\nabla e_{\mu}^{n+1}\|_0^2=-(T_{\phi_{1B}}^{n+1},e_{\mu}^{n+1})- (A_3^{n+1},e_{\phi}^{n+1}-e_{\phi}^{n}),\\
		\label{sec4_7_eq2_1B}
		&|e_{R}^{n+1}|^2-|e_{R}^{n}|^2+|e_{R}^{n+1}-e_{R}^{n}|^2= -\Delta tA_4^{n+1}e_{R}^{n+1}-2e_{R}^{n+1}T_{R_{1B}}^{n+1}.
		\end{align}
	\end{subequations}
	Now, the right-hand side terms of \eqref{sec4_7_eq1_1B} can be treated as follows. 
	\begin{align*}
	-(T_{\phi_{1B}}^{n+1},e_{\mu}^{n+1})&\leq\frac{\Delta t}{12}\|\nabla e_{\mu}^{n+1}\|_0^2+\frac{3}{\Delta t}\|(-\Delta)^{-1/2}T_{\phi_{1B}}^{n+1}\|_0^2\\
	&\leq\frac{\Delta t}{12}\|\nabla e_{\mu}^{n+1}\|_0^2+C\Delta t^2\int_{t^{n}}^{t^{n+1}}\|\phi_{tt}(s)\|_{-1}^2ds,\\
	-e_R^{n}(R^{n}+R(t^{n}))\left(\frac{h(\phi^{n})}{\bm{E}[\phi^{n}]},e_{\phi}^{n+1}-e_{\phi}^{n}\right)&=e_R^{n}(R^{n}+R(t^{n}))\left(\frac{h(\phi^{n})}{\bm{E}[\phi^{n}]},\Delta t\Delta e_{\mu}^{n+1}-T_{\phi_{1B}}^{n+1}\right)\\
	&\leq Ce_{R}^{n}\left(\Delta t\|\nabla e_{\mu}^{n+1}\|_0+\|(-\Delta)^{-1/2}T_{\phi_{1B}}^{n+1}\|_0\right)\bigg{\|}\frac{\nabla h(\phi^{n})}{\bm{E}[\phi^{n}]}\bigg{\|}_0\\
	&= Ce_{R}^{n}\left(\Delta t\|\nabla e_{\mu}^{n+1}\|_0+\|(-\Delta)^{-1/2}T_{\phi_{1B}}^{n+1}\|_0\right)\bigg{\|}\frac{ h'(\phi^{n})\nabla\phi^{n}}{\bm{E}[\phi^{n}]}\bigg{\|}_0\\
	&\leq\frac{\Delta t}{12}\|\nabla e_{\mu}^{n+1}\|_0^2+C\Delta t|e_{R}^{n}|^2+
	C\Delta t^2\int_{t^{n}}^{t^{n+1}}\|\phi_{tt}(s)\|_{-1}^2ds,
	\end{align*} 
	\begin{align*}
	-(R(t^{n})^2-R(t^{n+1})^2)\left(\frac{h(\phi^{n})}{\bm{E}[\phi^{n}]},e_{\phi}^{n+1}-e_{\phi}^{n}\right)&=(R(t^{n})^2-R(t^{n+1})^2)\left(\frac{h(\phi^{n})}{\bm{E}[\phi^{n}]},\Delta t\Delta e_{\mu}^{n+1}-T_{\phi_{1B}}^{n+1}\right)\\
	&\leq\frac{\Delta t}{12}\|\nabla e_{\mu}^{n+1}\|_0^2+C\Delta t|R(t^{n})-R(t^{n+1})|^2\\
	& \ \ +
	C\Delta t^2\int_{t^{n}}^{t^{n+1}}\|\phi_{tt}(s)\|_{-1}^2ds,\\
	-R(t^{n+1})^2\left(\frac{h(\phi^{n})}{\bm{E}[\phi^{n}]}-\frac{h(\phi(t^{n}))}{\bm{E}[\phi(t^{n})]},e_{\phi}^{n+1}-e_{\phi}^{n}\right)&=R(t^{n+1})^2\left(\frac{h(\phi^{n})}{\bm{E}[\phi^{n}]}-\frac{h(\phi(t^{n}))}{\bm{E}[\phi(t^{n})]},\Delta t\Delta e_{\mu}^{n+1}-T_{\phi_{1B}}^{n+1}\right)
	\\
	&\leq\frac{\Delta t}{12}\|\nabla e_{\mu}^{n+1}\|_0^2+C\Delta t\Big{\|}\frac{\nabla h(\phi^{n})}{\bm{E}[\phi^{n}]}-\frac{\nabla h(\phi(t^{n}))}{\bm{E}[\phi(t^{n})]}\Big{\|}_0^2\\
	& \ \ +
	C\Delta t^2\int_{t^{n}}^{t^{n+1}}\|\phi_{tt}(s)\|_{-1}^2ds,\\
	-R(t^{n+1})^2\left(\frac{h(\phi(t^{n}))}{\bm{E}[\phi(t^{n})]}-\frac{h(\phi(t^{n+1}))}{\bm{E}[\phi(t^{n+1})]},e_{\phi}^{n+1}-e_{\phi}^{n}\right)&=R(t^{n+1})^2\left(\frac{h(\phi(t^{n}))}{\bm{E}[\phi(t^{n})]}-\frac{h(\phi(t^{n+1}))}{\bm{E}[\phi(t^{n+1})]},\Delta t\Delta e_{\mu}^{n+1}-T_{\phi_{1B}}^{n+1}\right)\\
	&\leq\frac{\Delta t}{12}\|\nabla e_{\mu}^{n+1}\|_0^2+C\Delta t\Big{\|}\frac{\nabla h(\phi(t^{n}))}{ \bm{E}[\phi(t^{n})]}-\frac{\nabla h(\phi(t^{n+1}))}{\bm{E}[\phi(t^{n+1})]}\Big{\|}_0^2\\
	& \ \ +
	C\Delta t^2\int_{t^{n}}^{t^{n+1}}\|\phi_{tt}(s)\|_{-1}^2ds.
	\end{align*}
	Next, the right-hand side terms of \eqref{sec4_7_eq2_1B} can be treated as follows. 
	\begin{equation*}
	\left\{
	\begin{split}
	&-\frac{\Delta t|e_{R}^{n+1}|^2}{\bm{E}[\phi^{n+1}]}\int_{\Omega}|\nabla\mu^{n+1}|^2d\bm{x}\leq C\Delta t\|\nabla\mu^{n+1}\|_0^2|e_{R}^{n+1}|^2, \\
	&-\frac{\Delta te_{R}^{n+1}R(t^{n+1})}{\bm{E}[\phi^{n+1}]}\int_{\Omega}(|\nabla\mu^{n+1}|^2-|\nabla\mu(t^{n+1})|^2)d\bm{x}\\
	&\qquad\qquad
	\leq\frac{\Delta te_{R}^{n+1}}{\bm{E}[\phi^{n+1}]}\|\nabla e_\mu^{n+1}\|_0\|\nabla\mu^{n+1}+\nabla\mu(t^{n+1})\|_0\\
	&\qquad\qquad
	\leq C\Delta t(\|\nabla\mu^{n+1}\|_0^2+1)|e_{R}^{n+1}|^2+\frac{\Delta t}{12}\|\nabla e_\mu^{n+1}\|_0^2,\\
	&-\Delta tR(t^{n+1})e_{R}^{n+1}\left(\frac{1}{\bm{E}[\phi^{n+1}]}-\frac{1}{\bm{E}[\phi(t^{n+1})]}\right)\int_{\Omega}|\nabla\mu(t^{n+1})|^2d\bm{x}\\
	&\qquad\qquad
	\leq C\Delta t\big{(}|e_{R}^{n+1}|^2+\left|\bm{E}[\phi^{n+1}]-\bm{E}[\phi(t^{n+1})]\right|^2\big{)}, \\
	&-2e_{R}^{n+1}T_{R_{1B}}^{n+1}\leq\Delta t|e_{R}^{n+1}|^2+\frac{1}{\Delta t}|T_{R_{1B}}^{n+1}|^2\\
	&\qquad\qquad \leq \Delta t|e_{R}^{n+1}|^2+
	C\Delta t^2\int_{t^{n}}^{t^{n+1}}\left|\frac{d^2R(s)}{dt^2}\right|^2ds.
	\end{split}
	\right.
	\end{equation*}
	Combining \eqref{sec4_eq9}, \eqref{sec4_eq11}, \eqref{sec4_eq12} and the above inequalities with \eqref{sec4_7_eq1_1B} and \eqref{sec4_7_eq2_1B},  we have
	\begin{align}\label{sec4_10_1B}
	&\frac{1}{2}(\|\nabla e_{\phi}^{n+1}\|_0^2-\|\nabla e_{\phi}^{n}\|_0^2)+\frac{\lambda}{2}(\|e_{\phi}^{n+1}\|_0^2-\|e_{\phi}^{n}\|_0^2) +|e_{R}^{n+1}|^2-|e_{R}^{n}|^2+\frac{\Delta t}{2}\|\nabla e_{\mu}^{n+1}\|_0^2
	\nonumber\\
	&\qquad+\frac{1}{2}\|\nabla e_{\phi}^{n+1}-\nabla e_{\phi}^{n}\|_0^2+\frac{\lambda}{2}\|e_{\phi}^{n+1}-e_{\phi}^{n}\|_0^2+|e_{R}^{n+1}-e_{R}^{n}|^2
	\nonumber\\
	&\leq C\Delta t(1+\|\nabla\mu^{n+1}\|_0^2)|e_{R}^{n+1}|^2+C\Delta t\big{(}\|\nabla e_{\phi}^{n+1}\|_0^2+\|e_{\phi}^{n+1}\|_0^2+\|\nabla e_{\phi}^n\|_0^2+\|e_{\phi}^n\|_0^2+|e_{R}^{n}|^2\big{)} 
	\nonumber\\
	&\ \ +C\Delta t^2\int_{t^{n}}^{t^{n+1}}(\|\phi_{tt}(s)\|_{-1}^2+\|\phi_t(s)\|_{1}^2)ds +C\Delta t^2\int_{t^{n}}^{t^{n+1}}\left(\left|\frac{d^2R(s)}{dt^2}\right|^2+\left|\frac{dR(s)}{dt}\right|^2\right)ds.
	\end{align}
	Summing up these equations for the
        indices from $0$ to $n$ and using the discrete Gronwall
        lemma \ref{sec2_Lemma3} conclude the proof.	
\end{proof}

\vspace{0.1in}
\noindent\underline{\bf Theorem \ref{sec5_Lemma2}:}
  Suppose $\phi_{in}\in H^3(\Omega)$ and the condition \eqref{Assume_H2} holds. The following inequality
  holds for all $n$ with the scheme \eqref{CH2},
	\begin{align*}
	&\|\nabla\phi^{n+1}\|_0^2+\|\nabla(2\phi^{n+1}-\phi^{n})\|_0^2+\frac{1}{2}\|\nabla(\phi^{n+1}-\phi^{n})\|_0^2+\frac{\lambda}{2}(\|\phi^{n+1}\|_0^2+\|2\phi^{n+1}-\phi^{n}\|_0^2)
	\\
	&\quad +\lambda\Delta t\|\Delta\phi^{n+1}\|_0^2+\frac{3\Delta t}{2}\|\nabla\Delta\phi^{n+1}\|_0^2+\frac{\Delta t}{2}\|\nabla\Delta(\phi^{n+1}-\phi^n)\|_0^2+\Delta t\sum_{k=0}^{n}\|\nabla\mu^{k+1}\|_0^2\leq\widehat{C}_5,
	\end{align*}
	where $\widehat{C}_5=(2\|\nabla\phi^{0}\|_0^2+\lambda\|\phi^{0}\|_0^2+\lambda\Delta t\|\Delta\phi^{0}\|_0^2+\frac{3\Delta t}{2}\|\nabla\Delta\phi^{0}\|_0^2)\exp\left(C(M)T\right)$.

\begin{proof} Taking the inner product of \eqref{CH2_eq1} and \eqref{CH2_eq2} with $2\Delta t \mu^{n+1}$ and $3\phi^{n+1}-4\phi^{n}+\phi^{n-1}$, respectively, we have
	\begin{align}\label{sec5_eq2}
	&\frac{1}{2}(\|\nabla\phi^{n+1}\|_0^2-\|\nabla\phi^{n}\|_0^2+\|\nabla(2\phi^{n+1}-\phi^{n})\|_0^2-\|\nabla(2\phi^{n}-\phi^{n-1})\|_0^2+\|\nabla(\phi^{n+1}-2\phi^{n}+\phi^{n-1})\|_0^2)\nonumber\\
	&\quad+\frac{\lambda}{2}(\|\phi^{n+1}\|_0^2-\|\phi^{n}\|_0^2+\|2\phi^{n+1}-\phi^{n}\|_0^2-\|2\phi^{n}-\phi^{n-1}\|_0^2+\|\phi^{n+1}-2\phi^{n}+\phi^{n-1}\|_0^2)+2\Delta t\|\nabla\mu^{n+1}\|_0^2
	\nonumber\\
	&=-|\xi^{n+1}_{_{2A}}|^2(h(\overline{\phi}^{n}),3\phi^{n+1}-4\phi^{n}+\phi^{n-1})
	\leq\Delta t\|\nabla\mu^{n+1}\|_0^2+|\xi^{n+1}_{_{2A}}|^4\Delta t\|\nabla h(\overline{\phi}^{n})\|_0^2.
	\end{align}
	Taking the inner product of \eqref{CH2_eq1} and \eqref{CH2_eq2}  with $\Delta(2\phi^{n+1}-\phi^{n})$ and $\Delta^2(2\phi^{n+1}-\phi^{n})$, respectively,
        and by using the boundary conditions, we can re-write
        the two resultant equations  into
	\begin{align}\label{sec5_eq3}
	&\frac{1}{2}(\|\nabla\phi^{n+1}\|_0^2-\|\nabla\phi^{n}\|_0^2+\|\nabla(2\phi^{n+1}-\phi^{n})\|_0^2-\|\nabla(2\phi^{n}-\phi^{n-1})\|_0^2+\|\nabla(\phi^{n+1}-2\phi^{n}+\phi^{n-1})\|_0^2)\nonumber\\
	&\quad +2\|\nabla(\phi^{n+1}-\phi^{n})\|_0^2+\frac{1}{2}\left(\|\nabla(\phi^{n+1}-\phi^{n})\|_0^2-\|\nabla(\phi^{n}-\phi^{n-1})\|_0^2+\|\nabla(\phi^{n+1}-2\phi^{n}+\phi^{n-1})\|_0^2\right)
	\nonumber\\
	&\quad+2\Delta t\|\nabla\Delta\phi^{n+1}\|_0^2+\Delta t\left(\|\nabla\Delta\phi^{n+1}\|_0^2-\|\nabla\Delta\phi^{n}\|_0^2+\|\nabla\Delta(\phi^{n+1}-\phi^{n})\|_0^2\right)+2\lambda\Delta t\|\Delta\phi^{n+1}\|_0^2
	\nonumber\\
	&\quad+\lambda\Delta t\left(\|\Delta\phi^{n+1}\|_0^2-\|\Delta\phi^{n}\|_0^2+\|\Delta\phi^{n+1}-\Delta\phi^{n}\|_0^2\right)\nonumber\\
	&= 2|\xi^{n+1}_{_{2A}}|^2\Delta t(\nabla h(\overline{\phi}^{n}),\nabla\Delta\phi^{n+1})+2|\xi^{n+1}_{_{2A}}|^2\Delta t(\nabla h(\overline{\phi}^{n}),\nabla\Delta(\phi^{n+1}-\phi^{n}))\nonumber\\
	&\leq \Delta t\|\nabla\Delta\phi^{n+1}\|_0^2+\frac{\Delta t}{2}\|\nabla\Delta(\phi^{n+1}-\phi^{n}))\|_0^2+3\Delta t|\xi^{n+1}_{_{2A}}|^4\|\nabla h(\overline{\phi}^{n})\|_0^2.
   \end{align}
	Summing up \eqref{sec5_eq2} and \eqref{sec5_eq3}, we have
	\begin{align}\label{sec5_eq4}
	&\|\nabla\phi^{n+1}\|_0^2-\|\nabla\phi^{n}\|_0^2+\|\nabla(2\phi^{n+1}-\phi^{n})\|_0^2-\|\nabla(2\phi^{n}-\phi^{n-1})\|_0^2+2\|\nabla(\phi^{n+1}-\phi^{n})\|_0^2+\Delta t\|\nabla\mu^{n+1}\|_0^2
	\nonumber\\
	&\quad +\frac{1}{2}\left(\|\nabla(\phi^{n+1}-\phi^{n})\|_0^2-\|\nabla(\phi^{n}-\phi^{n-1})\|_0^2\right)+\frac{\lambda}{2}(\|\phi^{n+1}\|_0^2-\|\phi^{n}\|_0^2+\|2\phi^{n+1}-\phi^{n}\|_0^2
	\nonumber\\
	&\quad-\|2\phi^{n}-\phi^{n-1}\|_0^2)+\Delta t\|\nabla\Delta\phi^{n+1}\|_0^2+\Delta t(\|\nabla\Delta\phi^{n+1}\|_0^2-\|\nabla\Delta\phi^{n}\|_0^2)+\frac{\Delta t}{2}\|\nabla\Delta(\phi^{n+1}-\phi^{n})\|_0^2
	\nonumber\\
	&\quad +\lambda\Delta t(\|\Delta\phi^{n+1}\|_0^2-\|\Delta\phi^{n}\|_0^2+\|\Delta\phi^{n+1}-\Delta\phi^{n}\|_0^2)
	\nonumber\\
	&\leq 4|\xi^{n+1}_{_{2A}}|^4\Delta t\|\nabla h(\overline{\phi}^{n})\|_0^2.
	\end{align}
	Using Lemma \ref{sec2_Lemma1}, one finds
	\begin{align}\label{sec5_eq5}
	4|\xi^{n+1}_{_{2A}}|^4\|\nabla h(\overline{\phi}^{n})\|_0^2\leq 4|\xi^{n+1}_{_{2A}}|^4\|h'(\overline{\phi}^{n})\|_{0,\infty}^2\|\nabla\overline{\phi}^{n}\|_0^2\leq 4C|\xi^{n+1}_{_{2A}}|^4\|\nabla\overline{\phi}^{n}\|_0^2(1+\|\overline{\phi}^{n}\|_{0,\infty}^4).
	\end{align}
	By using the same technique as in
        the proof of Theorem \ref{sec3_Lemma3}, the Cauchy Schwarz
        inequality and the triangle inequality, we have
	\begin{eqnarray*}
	4|\xi^{n+1}_{_{2A}}|^4\|\nabla h(\overline{\phi}^{n})\|_0^2\leq\left\{\begin {array}{lll}
	C(M)\|\nabla\overline{\phi}^{n}\|_0^2      & \mbox{for}\ d =1,\\
	C(M)\|\nabla\overline{\phi}^{n}\|_0^2+C(\epsilon_1, M)\|\nabla\overline{\phi}^{n}\|_0^2+\epsilon_1\|\nabla\Delta\overline{\phi}^{n}\|_0^2        & \mbox{for}\ d=2,\\
	C(M)\|\nabla\overline{\phi}^{n}\|_0^2+C(\epsilon_2, M)\|\nabla\overline{\phi}^{n}\|_0^2+\epsilon_2\|\nabla\Delta\overline{\phi}^{n}\|_0^2        & \mbox{for}\ d=3.	\end{array}\right.
	\end{eqnarray*}
  Noting that 
	\[
	\|\nabla\Delta\overline{\phi}^{n}\|_0^2\leq 2\|\nabla\Delta\phi^{n}\|_0^2+2\|\nabla\Delta(\phi^{n}-\phi^{n-1})\|_0^2,
	\]
	and by setting $\epsilon_1=\epsilon_2=\frac{1}{4}$
        and combining the above inequalities with \eqref{sec5_eq4},
        we obtain 
	\begin{align}\label{sec5_eq7}
	&\|\nabla\phi^{n+1}\|_0^2-\|\nabla\phi^{n}\|_0^2+\|\nabla(2\phi^{n+1}-\phi^{n})\|_0^2-\|\nabla(2\phi^{n}-\phi^{n-1})\|_0^2+\Delta t\|\nabla\mu^{n+1}\|_0^2
	\nonumber\\
	&\quad+\frac{1}{2}\left(\|\nabla(\phi^{n+1}-\phi^{n})\|_0^2-\|\nabla(\phi^{n}-\phi^{n-1})\|_0^2\right)+\frac{\lambda}{2}(\|\phi^{n+1}\|_0^2-\|\phi^{n}\|_0^2+\|2\phi^{n+1}-\phi^{n}\|_0^2\nonumber\\
	&\quad-\|2\phi^{n}-\phi^{n-1}\|_0^2)+\lambda\Delta t(\|\Delta\phi^{n+1}\|_0^2-\|\Delta\phi^{n}\|_0^2)+\frac{3\Delta t}{2}(\|\nabla\Delta\phi^{n+1}\|_0^2-\|\nabla\Delta\phi^{n}\|_0^2)
	\nonumber\\
	&\quad +\frac{\Delta t}{2}(\|\nabla\Delta(\phi^{n+1}-\phi^{n})\|_0^2-\|\nabla\Delta(\phi^{n}-\phi^{n-1})\|_0^2)
	\nonumber\\
	&\leq C(M)\Delta t\|\nabla\overline{\phi}^{n}\|_0^2\leq C(M)\Delta t(\|\nabla\phi^{n}\|_0^2+\|\nabla\phi^{n-1}\|_0^2).
	\end{align}
        We conclude the proof by
	taking the sum of \eqref{sec5_eq7} for the indices
         from $0$ to $n$ 
         and using the discrete Gronwall lemma \ref{sec2_Lemma4}.         
\end{proof}

\vspace{0.1in}
\noindent\underline{\bf Theorem \ref{sec5_Lemma3}:}
Suppose $\phi_{in}\in H^4(\Omega)$, and
the conditions for Lemmas \ref{sec2_Lemma1} and \ref{sec2_Lemma2} hold.
The following inequality holds for all $n$ with the scheme \eqref{CH2}, 
\begin{align*}
&\frac{1}{2}\|\Delta\phi^{n+1}\|_0^2+\frac{1}{2}\|\Delta(2\phi^{n+1}-\phi^{n})\|_0^2+\frac{1}{2}\|\Delta(\phi^{n+1}-\phi^{n})\|_0^2+\lambda\Delta t\|\nabla\Delta\phi^{n+1}\|_0^2\\
&\qquad +\frac{3\Delta t}{2}\|\Delta^2\phi^{n+1}\|_0^2+\frac{\Delta t}{2}\|\Delta^2(\phi^{n+1}-\phi^{n})\|_0^2+\lambda\Delta t\sum_{k=0}^{n}\|\nabla\Delta\phi^{k+1}\|_0^2\leq\hat{C}_6,
\end{align*}
where $\hat{C}_6=\|\Delta\phi^{0}\|_0^2+\frac{3\Delta t}{2}\|\Delta^2\phi^{0}\|_0^2+\lambda\Delta t\|\nabla\Delta\phi^0\|_0^2+C(M)T$.

\begin{proof} Multiplying \eqref{CH2_eq1} by
  $2\Delta t\Delta^2(2\phi^{n+1}-\phi^{n})$ and combining \eqref{CH2_eq1} with \eqref{CH2_eq2}, we obtain
	\begin{align}\label{sec5_eq8}
	&\frac{1}{2}\left(\|\Delta\phi^{n+1}\|_0^2-\|\Delta\phi^{n}\|_0^2+\|\Delta(2\phi^{n+1}-\phi^{n})\|_0^2-\|\Delta(2\phi^{n}-\phi^{n-1})\|_0^2+\|\Delta(\phi^{n+1}-2\phi^{n}+\phi^{n-1})\|_0^2\right)
	\nonumber\\
	&\quad  
	+2\|\Delta(\phi^{n+1}-\phi^{n})\|_0^2+\frac{1}{2}\left(\|\Delta(\phi^{n+1}-\phi^{n})\|_0^2-\|\Delta(\phi^{n}-\phi^{n-1})\|_0^2+\|\Delta(\phi^{n+1}-2\phi^{n}+\phi^{n-1})\|_0^2\right)
	\nonumber\\
	&\quad
	+2\Delta t\|\Delta^2\phi^{n+1}\|_0^2+\Delta t\left(\|\Delta^2\phi^{n+1}\|_0^2
	-\|\Delta^2\phi^{n}\|_0^2+\|\Delta^2(\phi^{n+1}-\phi^{n})\|_0^2\right)
	\nonumber\\
	&\quad+2\lambda\Delta t\|\nabla\Delta\phi^{n+1}\|_0^2+\lambda\Delta t\left(\|\nabla\Delta\phi^{n+1}\|_0^2
	-\|\nabla\Delta\phi^{n}\|_0^2+\|\nabla\Delta(\phi^{n+1}-\phi^{n})\|_0^2\right)
	\nonumber\\
	&=2|\xi^{n+1}_{_{2A}}|^2\Delta t(\Delta h(\overline{\phi}^{n}),\Delta^2\phi^{n+1})+2|\xi^{n+1}_{_{2A}}|^2\Delta t(\Delta h(\overline{\phi}^{n}),\Delta^2(\phi^{n+1}-\phi^n))
	\nonumber\\
	&\leq\frac{\Delta t}{2}\|\Delta^2\phi^{n+1}\|_0^2+\frac{\Delta t}{2}\|\Delta^2(\phi^{n+1}-\phi^n)\|_0^2+4|\xi^{n+1}_{_{2A}}|^4\Delta t\|\Delta h(\overline{\phi}^{n})\|_0^2.
	\end{align}
	According to \eqref{Assume_H5}, for any $\epsilon>0$, there exists a constant $C(\epsilon,M)$ depending on $\epsilon$ such that
	\begin{align*}
	4|\xi^{n+1}_{_{2A}}|^4\Delta t\|\Delta h(\overline{\phi}^{n})\|_0^2\leq C(M)\Delta t(1+\|\Delta^2\overline{\phi}^{n}\|_0^{2\sigma})\leq \epsilon\Delta t\|\Delta^2\overline{\phi}^{n}\|_0^2+C(\epsilon,M)\Delta t.
	\end{align*}
	By the triangle inequality, we have
	\[
	\|\Delta^2\overline{\phi}^{n}\|_0^2\leq 2\|\Delta^2\phi^{n}\|_0^2+2\|\Delta^2(\phi^{n}-\phi^{n-1})\|_0^2.
	\]
	Combining the above inequalities with \eqref{sec5_eq8} and choosing $\epsilon=\frac{1}{4}$, we have
	\begin{align}\label{sec5_eq9}
	&\frac{1}{2}\left(\|\Delta\phi^{n+1}\|_0^2-\|\Delta\phi^{n}\|_0^2+\|\Delta(2\phi^{n+1}-\phi^{n})\|_0^2-\|\Delta(2\phi^{n}-\phi^{n-1})\|_0^2+\|\Delta(\phi^{n+1}-2\phi^{n}+\phi^{n-1})\|_0^2\right)
	\nonumber\\
	&\quad
	+2\|\Delta(\phi^{n+1}-\phi^{n})\|_0^2+\frac{1}{2}\left(\|\Delta(\phi^{n+1}-\phi^{n})\|_0^2-\|\Delta(\phi^{n}-\phi^{n-1})\|_0^2+\|\Delta(\phi^{n+1}-2\phi^{n}+\phi^{n-1})\|_0^2\right)
	\nonumber\\
	&\quad
	+\Delta t\|\Delta^2\phi^{n+1}\|_0^2+\frac{3\Delta t}{2}\left(\|\Delta^2\phi^{n+1}\|_0^2
	-\|\Delta^2\phi^{n}\|_0^2\right)+\frac{\Delta t}{2}\left(\|\Delta^2(\phi^{n+1}-\phi^{n})-\|\Delta^2(\phi^{n}-\phi^{n-1})\|_0^2\right)
	\nonumber\\
	&\quad+2\lambda\Delta t\|\nabla\Delta\phi^{n+1}\|_0^2+\lambda\Delta t\left(\|\nabla\Delta\phi^{n+1}\|_0^2
	-\|\nabla\Delta\phi^{n}\|_0^2+\|\nabla\Delta(\phi^{n+1}-\phi^{n})\|_0^2\right)
	\nonumber\\
	&\leq C(M)\Delta t.
	\end{align}
	We conclude the proof by taking the sum of this inequality
        for the indices from $0$ to $n$.
\end{proof}

\vspace{0.1in}
\noindent\underline{\bf Theorem \ref{sec5_Lemma3_2B}:}
Suppose $\phi_{in}\in H^3(\Omega)$ and the condition \eqref{Assume_H2} holds.
The following inequality holds for all n with the scheme \eqref{CH2_2B},
\begin{align*}
&\|\nabla\phi^{n+1}\|_0^2+\|\nabla(2\phi^{n+1}-\phi^{n})\|_0^2+\frac{1}{2}\|\nabla(\phi^{n+1}-\phi^{n})\|_0^2+\frac{\lambda}{2}(\|\phi^{n+1}\|_0^2+\|2\phi^{n+1}-\phi^{n}\|_0^2)\\
&\quad +\lambda\Delta t\|\Delta\phi^{n+1}\|_0^2+\frac{3\Delta t}{2}\|\nabla\Delta\phi^{n+1}\|_0^2+\frac{\Delta t}{2}\|\nabla\Delta(\phi^{n+1}-\phi^{n})\|_0^2+\Delta t\sum_{k=0}^{n}\|\nabla\mu^{k+1}\|_0^2\leq\widehat{C}_5,
\end{align*}
where $\widehat{C}_5$
is given in Theorem \ref{sec5_Lemma2}.

\begin{proof} By Lemma \ref{sec5_Lemma2_2B}, $\widehat{\xi}^n_{_{2B}}$ satisfies the conditions
  $\left|\widehat{\xi}^n_{_{2B}}\right|\leq \frac{3M}{\sqrt{C_0}}$
  and $\left|\widehat{\xi}^n_{_{2B}}\right|\leq \frac{C(M)}{\|\phi\|_1}$.
  In parallel to the proof of Theorem \ref{sec5_Lemma2},
  we take the same steps therein but replace $\xi^{n+1}_{_{2A}}$
  by $\widehat{\xi}^n_{_{2B}}$.
\end{proof}

\vspace{0.1in}
\noindent\underline{\bf Theorem \ref{sec5_Lemma4_2B}:}
Suppose $\phi^0\in H^4(\Omega)$, and the conditions for
Lemmas \ref{sec2_Lemma1} and \ref{sec2_Lemma2} hold.
Then the following inequality holds with the scheme \eqref{CH2_2B},
\begin{align*}
&\frac{1}{2}\|\Delta\phi^{n+1}\|_0^2+\frac{1}{2}\|\Delta(2\phi^{n+1}-\phi^{n})\|_0^2+\frac{1}{2}\|\Delta(\phi^{n+1}-\phi^{n})\|_0^2+\lambda\Delta t\|\nabla\Delta\phi^{n+1}\|_0^2\\
&\qquad +\frac{3\Delta t}{2}\|\Delta^2\phi^{n+1}\|_0^2+\frac{\Delta t}{2}\|\Delta^2(\phi^{n+1}-\phi^{n})\|_0^2+\lambda\Delta t\sum_{k=0}^{n}\|\nabla\Delta\phi^{k+1}\|_0^2\leq\widehat{C}_6,
\end{align*}
where $\widehat{C}_6$
is given in Theorem \ref{sec5_Lemma3}.

\begin{proof} The proof is essentially the same as for Theorem \ref{sec5_Lemma3}. One can refer to the proof of Theorem \ref{sec5_Lemma3}.
\end{proof}

\vspace{0.1in}
\noindent\underline{\bf Theorem \ref{sec5_Lemma5_2B}:}
Suppose the condition \eqref{regularity_2}, and the conditions for Theorems
\ref{sec5_Lemma3_2B} and \ref{sec5_Lemma4_2B} hold.
The following inequality holds for sufficiently small $\Delta t$,
\begin{equation*}
\frac{1}{2}\left(\|\nabla e_{\phi}^{n+1}\|_0^2+\|\nabla(2e_{\phi}^{n+1}-e_{\phi}^{n})\|_0^2\right)+\frac{\lambda}{2}\left(\|e_{\phi}^{n+1}\|_0^2+\|2e_{\phi}^{n+1}-e_{\phi}^{n}\|_0^2\right)+\frac{\Delta t}{2}\|\nabla e_{\mu}^{n+1}\|_0^2+|e_{R}^{n+1}|^2\leq\widehat{C}_7\Delta t^4,
\end{equation*}
where $\widehat{C}_7=C\exp(\Delta t\sum_{k=0}^{n+1}\frac{r^{k+1/2}}{1-r^{k+1/2}\Delta t})\int_{0}^{t^{n+1}}\left(\|\phi_{t}(s)\|_{1}^4+\|\phi_{t}(s)\|_{1}^2+\|\phi_{tt}(s)\|_{1}^2+\|\phi_{ttt}(s)\|_{-1}^2\right)ds$,
$r^{k+1/2}=1+\|\nabla\mu^{k+1/2}\|_0^2$,
and the constant $C$  depends on
$T$, $\phi_{in}$, $\Omega$, $\|\phi\|_{L^{\infty}(0,T;W^{3,\infty}(\Omega))}$, $\|\phi_t\|_{L^{\infty}(0,T;L^{2}(\Omega))}$ and $\|\mu\|_{L^{\infty}(0,T;H^{1}(\Omega))}$. 

\begin{proof} By subtracting \eqref{truncation2_2B} from \eqref{CH2_2B}, we have
	\begin{subequations}\label{error2_eq}
		\begin{align}
		\label{error2_eq1}
		&\frac{3e_{\phi}^{n+1}-4e_{\phi}^{n}+e_{\phi}^{n-1}}{2\Delta t}= \Delta e_{\mu}^{n+1}-\frac{1}{\Delta t}T_{\phi_{2B}}^{n+1},\\
		\label{error2_eq2}
		&e_{\mu}^{n+1}=-\Delta e_{\phi}^{n+1}+ \lambda e_{\phi}^{n+1} + A_5^{n+1},\\
		\label{error2_eq3}
		&\frac{e_{R}^{n+1}-e_{R}^{n}}{\Delta t}= -\frac{1}{2}A_6^{n+1}-\frac{1}{\Delta t}T_{R_{2B}}^{n+1},
		\end{align}
	\end{subequations}
	where 
	\begin{align*}
	A_5^{n+1}&=|\widehat{\xi}^{n}_{_{2B}}|^2h(\overline{\phi}^{n})-\frac{R(t^{n+1})^2}{\bm{E}[\phi(t^{n+1})]}h(\phi(t^{n+1}))\\
	&=\frac{(\overline{R}^{n})^2-\overline{R}(t^{n})^2}{\bm{E}[\overline{\phi}^{n}]}h(\overline{\phi}^{n})+\frac{\overline{R}(t^{n})^2-R(t^{n+1})^2}{\bm{E}[\overline{\phi}^{n}]}h(\overline{\phi}^{n})+
	R(t^{n+1})^2\left(\frac{h(\overline{\phi}^{n})}{\bm{E}[\overline{\phi}^{n}]}-\frac{h(\overline{\phi}(t^{n}))}{\bm{E}[\overline{\phi}(t^{n})]}\right)\\
	&\quad +R(t^{n+1})^2\left(\frac{h(\overline{\phi}(t^{n})))}{\bm{E}[\overline{\phi}(t^{n})]}-\frac{h(\phi(t^{n+1}))}{\bm{E}[\phi(t^{n+1})]}\right),\\
	A_6^{n+1}&=\frac{\xi^{n+1}_{_{2B}}}{\sqrt{\bm{E}[\widetilde{\phi}^{n+1/2}]}}\int_{\Omega}|\nabla\mu^{n+1/2}|^2d\bm{x}-\frac{R(t^{n+1})}{\sqrt{\bm{E}[\phi(t^{n+1})]}}\frac{1}{\sqrt{\bm{E}[\phi(t^{n+1/2})]}} \int_{\Omega}|\nabla\mu(t^{n+1/2})|^2d\bm{x}\\
	&=\frac{e_{R}^{n+1}}{\sqrt{\bm{E}[\widetilde{\phi}^{n+1/2}]\bm{E}[\phi^{n+1}]}}\int_{\Omega}|\nabla\mu^{n+1/2}|^2d\bm{x}+
	\frac{R(t^{n+1})}{\sqrt{\bm{E}[\widetilde{\phi}^{n+1/2}]\bm{E}[\phi^{n+1}]}}\int_{\Omega}(|\nabla\mu^{n+1/2}|^2-|\nabla\mu(t^{n+1/2})|^2)d\bm{x}\\
	& \quad + \frac{R(t^{n+1})}{\sqrt{\bm{E}[\phi^{n+1}]}}\left(\frac{1}{\sqrt{\bm{E}[\widetilde{\phi}^{n+1/2}]}}-\frac{1}{\sqrt{\bm{E}[\phi(t^{n+1/2})]}}\right)\int_{\Omega}|\nabla\mu(t^{n+1/2})|^2d\bm{x}\\
	&\quad +
	 \frac{R(t^{n+1})}{\sqrt{\bm{E}[\phi(t^{n+1/2})]}}\left(\frac{1}{\sqrt{\bm{E}[\phi^{n+1}]}}-\frac{1}{\sqrt{\bm{E}[\phi(t^{n+1})]}}\right)\int_{\Omega}|\nabla\mu(t^{n+1/2})|^2d\bm{x}.
	\end{align*}
	Taking the inner product of \eqref{error2_eq1} with $2\Delta te_{\mu}^{n+1}$ and \eqref{error2_eq2} with $3e_{\phi}^{n+1}-4e_{\phi}^{n}+e_{\phi}^{n-1}$, and multiplying \eqref{error2_eq3} with $2\Delta te_{R}^{n+1}$, we get the following:
	\begin{subequations}\label{sec5_eq15}
		\begin{align}
		\label{sec5_eq15_eq1}
		&\frac{1}{2}(\|\nabla e_{\phi}^{n+1}\|_0^2-\|\nabla e_{\phi}^{n}\|_0^2 +\|\nabla(2e_{\phi}^{n+1}-e_{\phi}^{n})\|_0^2-\|\nabla(2e_{\phi}^{n}-e_{\phi}^{n-1})\|_0^2)
		+2\Delta t\|\nabla e_{\mu}^{n+1}\|_0^2
		+\frac{\lambda}{2}(\|e_{\phi}^{n+1}\|_0^2-\|e_{\phi}^{n}\|_0^2
		\nonumber\\
		&\quad+\|2e_{\phi}^{n+1}-e_{\phi}^{n}\|_0^2-\|2e_{\phi}^{n}-e_{\phi}^{n-1}\|_0^2)\leq- (A_5^{n+1},3e_{\phi}^{n+1}-4e_{\phi}^{n}+e_{\phi}^{n-1})-2(T_{\phi_{2B}}^{n+1},e_{\mu}^{n+1}),\\
		\label{sec5_eq15_eq2}
		&|e_{R}^{n+1}|^2-|e_{R}^{n}|^2+|e_{R}^{n+1}-e_{R}^{n}|^2= -\Delta tA_6^{n+1}e_{R}^{n+1}-2e_{R}^{n+1}T_{R_{2B}}^{n+1}.
		\end{align}
	\end{subequations}
        
	The terms on the right hand side of \eqref{sec5_eq15_eq1} can be treated as follows. 
	\begin{align*}
	  \overline{e}_R^{n}(\overline{R}^{n}+\overline{R}(t^{n}))\left(\frac{h(\overline{\phi}^{n})}{\bm{E}[\overline{\phi}^{n}]},3e_{\phi}^{n+1}-4e_{\phi}^{n}+e_{\phi}^{n-1}\right)
          &=\overline{e}_R^{n}(\overline{R}^{n}+\overline{R}(t^{n}))\left(\frac{h(\overline{\phi}^{n})}{\bm{E}[\overline{\phi}^{n}]},2\Delta t\Delta e_{\mu}^{n+1}-2T_{\phi_{2B}}^{n+1}\right)
	\nonumber\\
	&\leq C\overline{e}_R^{n}\left(\Delta t\|\nabla e_{\mu}^{n+1}\|_0+\|(-\Delta)^{-1/2}T_{\phi_{2B}}^{n+1}\|_0\right)\bigg{\|}\frac{\nabla h(\overline{\phi}^{n})}{\bm{E}[\overline{\phi}^{n}]}\bigg{\|}_0
	\nonumber\\
	&\leq\frac{\Delta t}{5}\|\nabla e_{\mu}^{n+1}\|_0^2+C\Delta t|\overline{e}_R^{n}|^2+
	C\Delta t^4\int_{t^{n-1}}^{t^{n+1}}\|\phi_{ttt}(s)\|_{-1}^2ds,\\
		(\overline{R}(t^{n})^2-R(t^{n+1})^2)\left(\frac{h(\overline{\phi}^{n})}{\bm{E}[\overline{\phi}^{n}]},3e_{\phi}^{n+1}-4e_{\phi}^{n}+e_{\phi}^{n-1}\right)&\leq\frac{\Delta t}{5}\|\nabla e_{\mu}^{n+1}\|_0^2+C\Delta t|\overline{R}(t^{n})-R(t^{n+1})|^2\\
	&\quad+
	C\Delta t^4\int_{t^{n-1}}^{t^{n+1}}\|\phi_{ttt}(s)\|_{-1}^2ds,\\
	R(t^{n+1})^2\left(\frac{h(\overline{\phi}^{n})}{\bm{E}[\overline{\phi}^{n}]}-\frac{h(\overline{\phi}(t^{n}))}{\bm{E}[\overline{\phi}(t^{n})]},3e_{\phi}^{n+1}-4e_{\phi}^{n}+e_{\phi}^{n-1}\right)
	&\leq\frac{\Delta t}{5}\|\nabla e_{\mu}^{n+1}\|_0^2+C\Delta t\Big{\|}\frac{\nabla h(\overline{\phi}^{n})}{\bm{E}[\overline{\phi}^{n}]}-\frac{\nabla h(\overline{\phi}(t^{n}))}{\bm{E}[\overline{\phi}(t^{n})]}\Big{\|}_0^2\\
	& \quad+
	C\Delta t^4\int_{t^{n-1}}^{t^{n+1}}\|\phi_{ttt}(s)\|_{-1}^2ds,
	\end{align*}
   \begin{align*}
	R(t^{n+1})^2\left(\frac{h(\overline{\phi}(t^{n}))}{\bm{E}[\overline{\phi}(t^{n})]}-\frac{h(\phi(t^{n+1}))}{\bm{E}[\phi(t^{n+1})]},3e_{\phi}^{n+1}-4e_{\phi}^{n}+e_{\phi}^{n-1}\right)
	&\leq\frac{\Delta t}{5}\|\nabla e_{\mu}^{n+1}\|_0^2+
	C\Delta t^2\int_{t^{n-1}}^{t^{n+1}}\|\phi_{ttt}(s)\|_{-1}^2ds\\
	&\quad+
	C\Delta t\Big{\|}\frac{\nabla h(\overline{\phi}(t^{n}))}{ \bm{E}[\overline{\phi}(t^{n})]}-\frac{\nabla h(\phi(t^{n+1}))}{\bm{E}[\phi(t^{n+1})]}\Big{\|}_0^2,\\
	-2(T_{\phi_{2B}}^{n+1},e_{\mu}^{n+1})\leq\frac{\Delta t}{5}\|\nabla e_{\mu}^{n+1}\|_0^2+\frac{C}{\Delta t}\|(-\Delta)^{-1/2}T_{\phi_{2B}}^{n+1}\|_0^2
	&\leq\frac{\Delta t}{5}\|\nabla e_{\mu}^{n+1}\|_0^2+C\Delta t^4\int_{t^{n-1}}^{t^{n+1}}\|\phi_{ttt}(s)\|_{-1}^2ds.
	\end{align*}
	Next, the right-hand side terms of \eqref{sec5_eq15_eq2} can be treated as follows. 
	\begin{align*}
	&-2e_{R}^{n+1}T_{R_{2B}}^{n+1}\leq\Delta t|e_{R}^{n+1}|^2+\frac{C}{\Delta t}|T_{R_{2B}}^{n+1}|^2\leq \Delta t|e_{R}^{n+1}|^2+
	C\Delta t^4\int_{t^{n}}^{t^{n+1}}\left|\frac{d^3R(s)}{dt^3}\right|^2ds,\\
	&-\frac{\Delta t|e_{R}^{n+1}|^2}{\sqrt{\bm{E}[\widetilde{\phi}^{n+1/2}]\bm{E}[\phi^{n+1}]}}\int_{\Omega}|\nabla\mu^{n+1/2}|^2d\bm{x}\leq C\Delta t\|\nabla\mu^{n+1/2}\|_0^2|e_{R}^{n+1}|^2,
	\nonumber\\
	&-\frac{\Delta te_{R}^{n+1}R(t^{n+1})}{\sqrt{\bm{E}[\widetilde{\phi}^{n+1/2}]\bm{E}[\phi^{n+1}]}}\int_{\Omega}(|\nabla\mu^{n+1/2}|^2-|\nabla\mu(t^{n+1/2})|^2)d\bm{x}\leq C\Delta te_{R}^{n+1}\|\nabla e_\mu^{n+1/2}\|_0\|\nabla\mu^{n+1/2}+\nabla\mu(t^{n+1/2})\|_0
	\nonumber\\
	&\qquad\leq
	C\Delta t(\|\nabla\mu^{n+1/2}\|_0^2+1)|e_{R}^{n+1}|^2+\frac{\Delta t}{2}(\|\nabla e_\mu^{n+1}\|_0^2+\|\nabla e_\mu^{n}\|_0^2),
	\nonumber\\
	&-\frac{\Delta tR(t^{n+1})e_{R}^{n+1}}{\sqrt{\bm{E}[\phi^{n+1}]}}\left(\frac{1}{\sqrt{\bm{E}[\widetilde{\phi}^{n+1/2}]}}-\frac{1}{\sqrt{\bm{E}[\phi(t^{n+1/2})]}}\right)\int_{\Omega}|\nabla\mu(t^{n+1/2})|^2d\bm{x}
	\nonumber\\
	&\qquad
	\leq C\Delta te_{R}^{n+1}\left|\frac{1}{\sqrt{\bm{E}[\widetilde{\phi}^{n+1/2}]}}-\frac{1}{\sqrt{\bm{E}[\phi(t^{n+1/2})]}}\right|
	\leq C\Delta t|e_{R}^{n+1}|^2+C\Delta t\left|\frac{1}{\sqrt{\bm{E}[\widetilde{\phi}^{n+1/2}]}}-\frac{1}{\sqrt{\bm{E}[\phi(t^{n+1/2})]}}\right|^2,
	\nonumber\\
	&-\frac{\Delta tR(t^{n+1})e_{R}^{n+1}}{\sqrt{\bm{E}[\phi(t^{n+1/2})]}}\left(\frac{1}{\sqrt{\bm{E}[\phi^{n+1}]}}-\frac{1}{\sqrt{\bm{E}[\phi(t^{n+1})]}}\right)\int_{\Omega}|\nabla\mu(t^{n+1/2})|^2d\bm{x}
	\nonumber\\
	&\qquad\leq C\Delta te_{R}^{n+1}\left|\frac{1}{\sqrt{\bm{E}[\phi^{n+1}]}}-\frac{1}{\sqrt{\bm{E}[\phi(t^{n+1})]}}\right|\leq C\Delta t|e_{R}^{n+1}|^2+ C\Delta t\left|\frac{1}{\sqrt{\bm{E}[\phi^{n+1}]}}-\frac{1}{\sqrt{\bm{E}[\phi(t^{n+1})]}}\right|^2.
	\end{align*}	
	For the term $\frac{\nabla h(\overline{\phi}^{n})}{\bm{E}[\overline{\phi}^{n}]}-\frac{\nabla h(\overline{\phi}(t^{n}))}{\bm{E}[\overline{\phi}(t^{n})]}$, we rewrite it into
	\begin{equation*}
	\frac{\nabla h(\overline{\phi}^{n})}{\bm{E}[\overline{\phi}^{n}]}-\frac{\nabla h(\overline{\phi}(t^{n}))}{\bm{E}[\overline{\phi}(t^{n})]}=\frac{\nabla h(\overline{\phi}^{n})-\nabla h(\overline{\phi}(t^{n}))}{\bm{E}[\overline{\phi}^{n}]}
	+\frac{\nabla h(\overline{\phi}(t^{n}))\big{(}\bm{E}[\overline{\phi}(t^{n})]-\bm{E}[\overline{\phi}^{n}]\big{)}}{  \bm{E}[\overline{\phi}^{n}]\bm{E}[\overline{\phi}(t^{n})]}
	\end{equation*}
	In light of the H${\rm\ddot{o}}$lder's inequality and Sobolev
        embedding theorem, we have
	\begin{align*}
	\|\nabla h(\overline{\phi}^{n})-\nabla h(\overline{\phi}(t^{n}))\|_0\leq&\|(h'(\overline{\phi}^{n})-h'(\overline{\phi}(t^{n})))\nabla\overline{\phi}(t^{n})\|_0+\|h'(\overline{\phi}^{n})\nabla \overline{e}_{\phi}^n\|_0
	\nonumber\\
	\leq&C\|\nabla\overline{\phi}(t^{n})\overline{e}_{\phi}^n\|_{0}+C\|\nabla \overline{e}_{\phi}^n\|_0
	\leq C(\|\nabla\overline{\phi}(t^{n})\|_{0,3}\|\overline{e}_{\phi}^n\|_{0,6}+\|\nabla \overline{e}_{\phi}^n\|_0)
	\nonumber\\
	\leq&C\|\overline{\phi}(t^{n})\|_{2}\|\overline{e}_{\phi}^n\|_{1}+C\|\nabla \overline{e}_{\phi}^n\|_0\leq C(\|\nabla \overline{e}_{\phi}^n\|_0+\|\overline{e}_{\phi}^n\|_0).
	\end{align*}
	According to the definition of $\bm{E}[\phi]$, we have
	\begin{align}\label{sec5_eq17}
	\bm{E}[\overline{\phi}^{n}]-\bm{E}[\overline{\phi}(t^{n})]&= \frac{1}{2}\int_{\Omega}(\nabla\overline{\phi}^{n}+\nabla\overline{\phi}(t^{n}))\nabla \overline{e}_\phi^nd\bm{x}+\frac{\lambda}{2}\int_{\Omega}(\overline{\phi}^{n}+\overline{\phi}(t^{n}))\overline{e}_\phi^nd\bm{x}+\int_{\Omega}\big{(}H(\overline{\phi}^{n})-H(\overline{\phi}(t^{n})\big{)}d\bm{x}
	\nonumber\\
	&\leq C \|\nabla \overline{e}_\phi^n\|_0+C\|\overline{e}_\phi^n\|_0+\int_{\Omega}H'\big{(}\theta\overline{\phi}^{n}+(1-\theta)\overline{\phi}(t^{n})\big{)}\big{(}\overline{\phi}^{n}-\overline{\phi}(t^{n})\big{)}d\bm{x}
	\nonumber\\
	&\leq C(\|\nabla \overline{e}_\phi^n\|_0+\|\overline{e}_\phi^n\|_0).
	\end{align}
	Then, we have
	\begin{align}\label{sec5_eq18}
	\Big{\|}\frac{\nabla h(\overline{\phi}^{n})}{\bm{E}[\overline{\phi}^{n}]}-\frac{\nabla h(\overline{\phi}(t^{n}))}{\bm{E}[\overline{\phi}(t^{n})]}\Big{\|}_0^2
	=&\Bigg{\|}\frac{\nabla h(\overline{\phi}^{n})-\nabla h(\overline{\phi}(t^{n}))}{\bm{E}[\overline{\phi}^{n}]}
	+\frac{\nabla h(\overline{\phi}(t^{n}))\big{(}\bm{E}[\overline{\phi}(t^{n})]-\bm{E}[\overline{\phi}^{n}]\big{)}}{  \bm{E}[\overline{\phi}^{n}]\bm{E}[\overline{\phi}(t^{n})]}\Bigg{\|}_0^2
	\nonumber\\
	\leq&C\|\nabla h(\overline{\phi}^{n})-\nabla h(\overline{\phi}(t^{n}))\|_0^2+C\|\nabla h(\overline{\phi}(t^{n}))\|_0^2\big{|}\bm{E}[\overline{\phi}(t^{n})]-\bm{E}[\overline{\phi}^{n}]\big{|}^2
	\nonumber\\
	\leq&C(\|\nabla \overline{e}_{\phi}^n\|_0^2+\|\overline{e}_{\phi}^n\|_0^2).
	\end{align}
	Similarly, 
	\begin{align}\label{sec5_eq19}
	\Big{\|}\frac{\nabla h(\overline{\phi}(t^{n}))}{ \bm{E}[\overline{\phi}(t^{n})]}-\frac{\nabla h(\phi(t^{n+1}))}{\bm{E}[\phi(t^{n+1})]}\Big{\|}_0^2
	\leq&C\|\nabla h(\overline{\phi}(t^{n}))-\nabla h(\phi(t^{n+1}))\|_0^2+C\|\nabla h(\phi(t^{n+1}))\|_0^2\big{|}\bm{E}[\phi(t^{n+1})]-\bm{E}[\overline{\phi}(t^{n})]\big{|}^2
	\nonumber\\
	\leq&C(\|\nabla\phi(t^{n+1})-\nabla\overline{\phi}(t^{n})\|_0^2+\|\phi(t^{n+1})-\overline{\phi}(t^{n})\|_0^2)
	\nonumber\\
	\leq&C\Delta t^3\int_{t^{n-1}}^{t^{n+1}}\|\phi_{tt}(s)\|_{1}^2ds,
	\nonumber\\
	\left|\frac{1}{\sqrt{\bm{E}[\widetilde{\phi}^{n+1/2}]}}-\frac{1}{\sqrt{\bm{E}[\phi(t^{n+1/2})]}}\right|^2
	&=\frac{\left|\bm{E}[\phi(t^{n+1/2})]-\bm{E}[\widetilde{\phi}^{n+1/2}]\right|^2}{\bm{E}[\widetilde{\phi}^{n+1/2}]\bm{E}[\phi(t^{n+1/2})](\sqrt{\bm{E}[\widetilde{\phi}^{n+1/2}]}+\sqrt{\bm{E}[\phi(t^{n+1/2})]})^2}
	\nonumber\\
	&\leq C\left|\bm{E}[\phi(t^{n+1/2})]-\bm{E}[\widetilde{\phi}(t^{n+1/2})]\right|^2+C\left|\bm{E}[\widetilde{\phi}(t^{n+1/2})]-\bm{E}[\widetilde{\phi}^{n+1/2}]\right|^2
	\nonumber\\
	&\leq C(\|\nabla \widetilde{e}_{\phi}^{n+1/2}\|_0^2+\|\widetilde{e}_{\phi}^{n+1/2}\|_0^2)+C\Delta t^3\int_{t^{n-1}}^{t^{n+1}}\|\phi_{tt}(s)\|_{1}^2ds,
	\nonumber\\
	\left|\frac{1}{\sqrt{\bm{E}[\phi^{n+1}]}}-\frac{1}{\sqrt{\bm{E}[\phi(t^{n+1})]}}\right|^2
	&=\frac{\left|\bm{E}[\phi(t^{n+1})]-\bm{E}[\phi^{n+1}]\right|^2}{\bm{E}[\phi^{n+1}]\bm{E}[\phi(t^{n+1})](\sqrt{\bm{E}[\phi^{n+1}]}+\sqrt{\bm{E}[\phi(t^{n+1})]})^2}\nonumber\\
	&\leq C(\|\nabla e_{\phi}^{n+1}\|_0^2+\|e_{\phi}^{n+1}\|_0^2),\nonumber\\
	|\overline{R}(t^{n})-R(t^{n+1})|^2&=|2R(t^{n})-R(t^{n-1})-R(t^{n+1})|^2\leq C\Delta t^3\int_{t^{n-1}}^{t^{n+1}}\left|\frac{d^2R(s)}{dt^2}\right|^2ds.
	\end{align}
	By combining the above inequalities with \eqref{sec5_eq15_eq1} and \eqref{sec5_eq15_eq2},  we have
	\begin{align}\label{sec5_eq20}
	&\frac{1}{2}(\|\nabla e_{\phi}^{n+1}\|_0^2-\|\nabla e_{\phi}^{n}\|_0^2+\|\nabla(2e_{\phi}^{n+1}-e_{\phi}^{n})\|_0^2-\|\nabla(2e_{\phi}^{n}-e_{\phi}^{n-1})\|_0^2)+\frac{\Delta t}{2}(\|\nabla e_{\mu}^{n+1}\|_0^2-\|\nabla e_{\mu}^{n}\|_0^2)
	\nonumber\\
	&\qquad +\frac{\lambda}{2}(\|e_{\phi}^{n+1}\|_0^2-\|e_{\phi}^{n}\|_0^2+\|2e_{\phi}^{n+1}-e_{\phi}^{n}\|_0^2-\|2e_{\phi}^{n}-e_{\phi}^{n-1}\|_0^2)
	+|e_{R}^{n+1}|^2-|e_{R}^{n}|^2+|e_{R}^{n+1}-e_{R}^{n}|^2
	\nonumber\\
    &\leq C\Delta t(1+\|\nabla\mu^{n+1/2}\|_0^2)|e_{R}^{n+1}|^2+C\Delta t\big{(}\|\nabla e_{\phi}^{n+1}\|_0^2+\|e_{\phi}^{n+1}\|_0^2+\|\nabla e_{\phi}^{n}\|_0^2+\|e_{\phi}^{n}\|_0^2+\|\nabla e_{\phi}^{n-1}\|_0^2
	\nonumber\\
	&\ \ +\|e_{\phi}^{n-1}\|_0^2+|\overline{e}_{R}^{n}|^2\big{)}+C\Delta t^4\int_{t^{n-1}}^{t^{n+1}}\left(\|\phi_{ttt}(s)\|_{-1}^2+\|\phi_{tt}(s)\|_{1}^2+\left|\frac{d^3R(s)}{dt^3}\right|^2+\left|\frac{d^2R(s)}{dt^2}\right|^2\right)ds.
	\end{align}
	Summing up \eqref{sec5_eq20} for the indices
        from $0$ to $n$, we get
	\begin{align}\label{sec5_eq21}
	&\frac{1}{2}\left(\|\nabla e_{\phi}^{n+1}\|_0^2+\|\nabla(2e_{\phi}^{n+1}-e_{\phi}^{n})\|_0^2\right)+\frac{\lambda}{2}\left(\|e_{\phi}^{n+1}\|_0^2+\|2e_{\phi}^{n+1}-e_{\phi}^{n}\|_0^2\right)+\frac{\Delta t}{2}\|\nabla e_{\mu}^{n+1}\|_0^2+|e_{R}^{n+1}|^2
	\nonumber\\
	&\leq C\Delta t\sum_{k=0}^{n}(1+\|\nabla\mu^{k+1/2}\|_0^2)|e_{R}^{k+1}|^2+C\Delta t\sum_{k=0}^{n}\big{(}\|\nabla e_{\phi}^{k+1}\|_0^2+\|e_{\phi}^{k+1}\|_0^2+\|\nabla e_{\phi}^{k}\|_0^2+\|e_{\phi}^{k}\|_0^2+\|\nabla e_{\phi}^{k-1}\|_0^2\big{)}
	\nonumber\\
	&\ \ +C\Delta t\sum_{k=0}^{n}(\|e_{\phi}^{k-1}\|_0^2+|e_{R}^{k}|^2+|e_{R}^{k-1}|^2)+C\Delta t^4
	\nonumber\\
	&\leq C\Delta t\sum_{k=0}^{n}(1+\|\nabla\mu^{k+1/2}\|_0^2)\left(|e_{R}^{k+1}|^2+\|\nabla e_{\phi}^{k+1}\|_0^2+\|e_{\phi}^{k+1}\|_0^2\right)+C\Delta t^4.
	\end{align} 
	Based on Theorem \ref{sec5_Lemma3_2B} and the  triangle inequality $\|\nabla\mu^{k+1/2}\|_0^2\leq \frac{1}{2}(\|\nabla\mu^k\|_0^2+\|\nabla\mu^{k+1}\|_0^2)$,
        we have
	\[
	\Delta t\sum_{k=0}^{n+1}\|\nabla\mu^{k+1/2}\|_0^2\leq \widehat{C}_5.
	\]
	We then use the discrete Gronwall lemma \ref{sec2_Lemma3} to finish the
        proof.
\end{proof}

\section*{Acknowledgement}
This work was partially supported by
NSF (DMS-1522537) and the China Scholarship Council (CSC-201906280202).

\bibliographystyle{plain}
\bibliography{engstab,adjstab,yx_ref,obc,mypub,nse,sem,contact_line,interface,multiphase}

\end{document}